\documentclass{amsart}

\usepackage{tikz}
\usepackage{verbatim}
\usepackage{amssymb}
\usepackage{eepic}
\usepackage{color}

\allowdisplaybreaks

\numberwithin{equation}{section}

\newtheorem{theorem}{Theorem}[section]
\newtheorem{lemma}[theorem]{Lemma}
\newtheorem{proposition}[theorem]{Proposition}
\newtheorem{corollary}[theorem]{Corollary}
\newtheorem{remark}[theorem]{Remark}

\newtheorem{TheoA}{Theorem A1}
\newtheorem{TheoAA}{Theorem A2}
\newtheorem{TheoAAA}{Theorem A3}
\newtheorem{TheoB}{Theorem B}

\newcommand{\N}{\mathbb{N}}
\newcommand{\Z}{\mathbb{Z}}
\newcommand{\R}{\mathbb{R}}
\newcommand{\C}{\mathbb{C}}
\newcommand{\F}{\mathbb{F}}

\newcommand{\summ}{\sum\nolimits}

\def\G{\mathrm{G}}
\def\1{\mathbf{1}}
\def\H{\mathcal{H}}

\def\M{\mathcal{M}}

\def\P{\mathcal{P}}
\def\T{\mathcal{T}}
\def\S{\mathcal{S}}

\def\k{\underline{k}} 
\def\l{\underline{\ell}} 
 
\def\x{\underline{\xi}} 
\def\z{\underline{\zeta}} 
\def\d{\underline{d}} 
\def\ibar{\underline{i}}

\renewcommand{\L}{\mathcal{L}}

\newcommand{\dem}{\noindent {\bf Proof. }}

\newcommand{\fin}{\hspace*{\fill} $\square$ \vskip0.2cm}

\addtocounter{tocdepth}{0}

\begin{document}

\addtolength{\parskip}{+1ex}

\title[Hypercontractivity in group von Neumann algebras]{Hypercontractivity \\ in group von Neumann algebras}

\author[Junge, Palazuelos, Parcet and Perrin]
{Marius Junge, Carlos Palazuelos, \\ Javier Parcet and Mathilde Perrin}

\maketitle

\vskip-35pt \null 

\begin{abstract}
In this paper, we provide a combinatorial/numerical method to establish \hskip-1pt new \hskip-1pt hypercontractivity estimates in group \hskip-1pt von Neumann algebras. \hskip-3pt We will illustrate our method with free groups, triangular groups and finite cyclic groups, for which we shall obtain optimal time hypercontractive $L_2 \to L_q$ inequalities with respect to the Markov process given by the word length and with $q$ an even integer. Interpolation and differentiation also yield general $L_p \to L_q$ hypercontrativity for $1 < p \le q < \infty$ via logarithmic Sobolev inequalities. Our method admits further applications to other discrete groups without small loops as far as the numerical part ---which varies from one group to another--- is implemented and tested in a computer. We also develop another combinatorial method which does not rely on computational estimates and provides (non-optimal) $L_p \to L_q$ hypercontractive inequalities for a larger class of groups/lengths, including any finitely generated group equipped with a conditionally negative word length, like infinite Coxeter groups. Our second method also yields hypercontractivity bounds for groups admitting a finite dimensional proper cocycle. Hypercontractivity fails for conditionally negative lengths in groups satisfying Kazhdan property (T). 
\end{abstract}




\section*{{\bf Introduction and main results}}

Given a periodic function $f: \R \to \C$, the decay of its Fourier coefficients is closely related to the integrability properties of $f$. The Riemann-Lebesgue lemma shows that $\widehat{f}(n) \to 0$ as $|n| \to \infty$ for $f \in L_1(\mathbb{T})$, but this convergence could be arbitrarily slow. On the contrary, Plancherel's theorem goes further and gives $\sum_n |\widehat{f}(n)|^2 < \infty$ when $f \in L_2(\mathbb{T})$. Given $1 < p < 2$, a classical problem in harmonic analysis is to determine conditions on a weight $\zeta: \Z \to \R_+$ so that $$f \in L_p(\mathbb{T}) \ \Rightarrow \ \sum_{n \in \Z} \zeta(n) |\widehat{f}(n)|^2 < \infty.$$ Rudin's notion of $\Lambda_p$-set was motivated by this problem for characteristic functions $\zeta$, while Hausdorff-Young inequality immediately provides sufficient conditions for rational functions $\zeta$. The case of exponential functions $\zeta(n) = \exp(-2t|n|)$ leads to norm estimates for the Poisson semigroup on $\Z$. Our goal in this paper is to open a door through similar estimates replacing the frequency group $\Z$ above by other discrete groups $\G$ and the Poisson semigroup by other semigroups (Markovian or not) acting diagonally (Fourier multipliers) on the trigonometric system. Our problem reduces to norm estimates for operators of the form $$\S_{\psi,t}: \sum_{g \in \G} \widehat{f}(g) \lambda(g) \, \mapsto \, \sum_{g \in \G} e^{-t \psi(g)} \widehat{f}(g) \lambda(g),$$ with $t > 0$ and $\G$ a discrete group with left regular representation $\lambda: \G \to \mathcal{B}(\ell_2(\G))$.

The group von Neumann algebra $\mathcal{L}(\G)$ is the weak operator closure of the linear span of $\lambda(\G)$. If $e$ denotes the identity of $\G$, the algebra $\mathcal{L}(\G)$ comes equipped with the standard trace $\tau(f) = \langle \delta_e, f \delta_e \rangle$. Let $L_p(\mathcal{L}(\mathrm{G}))$ be the noncommutative $L_p$ space over the quantum probability space $(\mathcal{L}(\mathrm{G}), \tau)$ with $\|f\|_p^p = \tau |f|^p$. We invite the reader to check that $L_p(\mathcal{L}(\G)) = L_p(\mathbb{T})$ for $\G = \Z$, after identifying $\lambda_{\Z}(k)$ with $e^{2\pi i k \cdot}$. In general, the absolute value and the power $p$ are obtained from functional calculus for the (unbounded) operator $f$ on the Hilbert space $\ell_2(\G)$, we refer to \cite{PX} for further details. Markovian semigroups are composed of self-adjoint, completely positive and unital maps. By Schoenberg's theorem \cite{Sc}, we know that $\S_{\psi,t}$ is Markovian iff $\psi(e)=0$, $\psi(g) = \psi(g^{-1})$ and $\summ_g a_g = 0 \Rightarrow \summ_{g,h} \overline{a_g} a_h \psi(g^{-1}h) \le 0$. Any such $\psi: \G \to \R_+$ is called a conditionally negative length function.  As we shall see, the additional assumptions below play also a crucial role for hypercontractivity estimates in the group algebra $\mathcal{L}(\G)$
\begin{itemize}
\item[$\bullet$] Spectral gap: $\sigma = \inf_{g \neq e} \psi(g) > 0$, 

\vskip2pt

\item[$\bullet$] Subadditivity: $\psi(gh) \le \psi(g) + \psi(h)$,

\vskip2pt

\item[$\bullet$] Absence of 3-loops: $g_1 g_2 g_3 \neq e$ when $\psi(g_j) = \sigma$.
\end{itemize} 

By the Markovian nature of the semigroup, the $\S_{\psi,t}$'s are contractive maps on $L_p(\mathcal{L}(\G))$ for $1 \le p \le \infty$, which become more and more regular for $t$ large. The hypercontractivity problem for $1< p \le q < \infty$ consists in determining the optimal time $t_{p,q} > 0$ above which $$\|\S_{\psi,t} f \|_q \, \le \, \|f\|_p \qquad \mbox{for all} \qquad t \ge t_{p,q}.$$  This expected behavior has been studied since the early 70's in this and other related contexts. Originally, Bonami \cite{B} considered the group $\G = \Z_2$ with the standard length $\psi(g) = \delta_{g \neq e}$ and the cartesian powers $\Z_2^n$ with the Hamming distance $\psi(g_1, g_2, \ldots, g_n) = |\{j : g_j \neq e\}|$. Optimal hypercontractivity for $\Z_2$ is known as the two-point inequality, which was rediscovered by Gross \cite{G1} and used by Beckner \cite{Be} to obtain optimal constants for the Hausdorff-Young inequality. The two-point inequality has also shown a deep impact in both classical and quantum information theory. Shortly after, Weissler obtained in \cite{W} the same optimal time for the Poisson semigroup in the circle group, which is given in our terminology by $\G = \Z$ and $\psi(n) = |n|$. This yields for $1 < p \le 2 < \infty$ $$\Big( \sum_{n \in \Z} e^{-2 t |n|} |\widehat{f}(n)|^2 \Big)^{\frac12} \, \le \, \|f\|_p \quad \Leftrightarrow \quad t \ge \frac{1}{2} \log \frac{1}{p-1} = t_{p,2}.$$ Independently and almost simultaneously, hypercontractivity also emerged from quantum field theory, where Poisson processes are replaced by Orstein-Uhlenbeck like semigroups acting diagonally on generalized gaussians. Optimal estimates in the Bosonic, Fermonic and $q$-deformed cases can be found in \cite{Bi,CL,G2,N1,Se} while closely related results appear in \cite{G1,Ja,K,LR}. Also, similar methods apply for the heat-diffusion semigroup \cite{Be,W2}. Gaussian hypercontractivity bounds rely ultimately on two fundamental results by Gross \cite{G1} and Ball/Carlen/Lieb \cite{BCL}. In spite of this, no further significant results have appeared for Poisson-like processes in the trigonometric setting, which is perhaps explained from the lack of a convexity inequality \`a la Ball/Carlen/Lieb for group von Neumann algebras. Poisson-like hypercontractivity is closely related to Sobolev type theorems and norm estimates for Fourier multipliers in group von Neumann algebras \cite{J1,J2,JM,JMP,OR,PR}.

Let $\psi: \G \to \R_+$ be a conditionally negative length. If $\sigma = \inf_{g \neq e} \psi(g)$, pick $g \in \G$ with $\psi(g) \sim \sigma$ (an identity if the infimum is attained) and $\delta > 0$ small to estimate the $p$-norm of $f = \mathbf{1} + \delta (\lambda(g) + \lambda(g^{-1}))$ in the abelian algebra generated by it. The same estimate for the $q$-norm of $\S_{\psi,t} f \sim \mathbf{1} + \delta e^{-t \sigma} (\lambda(g) + \lambda(g^{-1}))$ yields the universal restriction $$t_{p,q} \, \ge \, \frac{1}{2 \sigma} \log \Big( \frac{q-1}{p-1} \Big) \, =: \, T(p,q,\sigma).$$ In particular, hypercontractivity imposes the existence of a spectral gap and the subadditivity of $\psi$ leads to a Poisson-like process. A far reaching goal is to determine whether $T(p,q,\sigma)$ is the optimal time for  any subadditive conditionally negative length admitting a spectral gap $\sigma$. We will show that $t_{p,q} > T(p,q,\sigma)$ in the presence of 3-loops, see \eqref{Noloop3} for details. Nevertheless, as we shall see, this case seems to be a pathological phenomenon and the problem is still meaningful after removing it. In the first part of this paper, we provide a general method to obtain optimal hypercontractivity $L_2 \to L_q$ estimates on discrete groups equipped with subadditive conditionally negative lengths which admit a spectral gap and satisfy the following conditions
\begin{itemize}
\item[$\bullet$] Growth: There exists $r > 0$ such that $\sum_{g \in \G} r^{\psi(g)}  < \infty$,

\vskip2pt

\item[$\bullet$] Cancellation: $g_1 g_2 \cdots g_m \neq e$ when $\psi(g_j) = \sigma$, $g_j \neq g_{j+1}^{-1}$ and $m$ small.
\end{itemize}
The first condition holds for instance under any exponential growth assumption $| \{ g \in \G : \psi(g) \le R \} | \le \mathrm{C} \rho^R$ for some $\rho > 1$. It also implies that the infimum $\sigma$ is attained, so that our second condition says that $\psi$ does not admit $m$-loops for $m$ small, not just $m=3$. The combinatorial nature of our approach forces $q$ to be an even integer. Standard interpolation and differentiation arguments lead to general $L_p \to L_q$ hypercontractivity estimates for any $1 < p \le q < \infty$ via log--Sobolev type inequalities. Unfortunately, our method fails for a finite number of terms which depend on the pair $(\G,\psi)$ and must be estimated with computer assistance. This forces us to apply it in specific scenarios by estimating each time the complete set of pathologies. We illustrate it with free, triangular and finite cyclic groups.    

Let us start with the finitely generated free groups $\F_n$ equipped with the word length $| \cdot|$ which measures the distance to $e$ in the Cayley graph. In this case, $| \cdot|$ is conditionally negative and the free Poisson semigroup is $$\mathcal{P}_t f \, = \, \sum_{g \in \F_n}^{\null} e^{-t|g|} \widehat{f}(g) \lambda(g).$$ It was introduced by Haagerup \cite{Ha} to prove that the reduced $\mathrm{C}^*$-algebra of the free group $\F_n$ has the metric approximation property. Free groups are morally an endpoint for the behavior of other  discrete groups with lower growth rate but more cancellation. Using an embedding of $\F_n$ into $\Z_2 * \Z_2 * \cdots * \Z_2$ with $2n$ factors, we exploit in \cite{JPPPR} a probabilistic approach which yields $L_p \to L_q$ hypercontractivity in $\F_n$ for $t \ge \log (q-1/p-1)$. This is twice the expected optimal time $T(p,q,1)$, while a more elaborated argument gives $1.173 \cdot T(p,q,1)$. Estimates below that constant turn out to be more challenging. Our first result gives the optimal time $L_2 \to L_q$ inequalities in $\F_2$ for $q \in 2\Z_+$ and reduces the constant $1.173$ above to $\log \, 3 \sim 1.099$ for the general case $1 < p \le q < \infty$ by adapting Gross program \cite{G1,G2}. We also get optimal estimates for $\F_n$ which are less general by computational limits. 

\begin{TheoA}
We find 
\begin{itemize}
\item[i)] If $f \in L_2(\mathcal{L}(\F_2))$ and $q \in 2 \Z_+$ $$\|\mathcal{P}_t f\|_{L_q(\mathcal{L}(\F_2))} \le \|f\|_{L_2(\mathcal{L}(\F_2))} \quad \Leftrightarrow \quad t \ge \frac12 \log (q-1).$$ In particular, $\|\mathcal{P}_t f\|_q \le \|f\|_p$ for $t \ge \frac{\log 3}{2} \log \frac{q-1}{p-1}$ and $1 < p \le q < \infty$.

\vskip5pt

\item[ii)] If $n \ge 3$ and $f \in L_2(\L(\F_n))$, we obtain optimal $L_2 \to L_q$ bounds for the free Poisson semigroup and $q$ a large even integer. Namely, the following holds for every $q \in 2 \Z_+$ greater than or equal to certain index $q(n) > 0$ depending on $n$ $$\|\mathcal{P}_t f\|_{L_q(\mathcal{L}(\F_n))} \le \|f\|_{L_2(\mathcal{L}(\F_n))} \quad \Leftrightarrow \quad t \ge \frac12 \log (q-1).$$
\end{itemize}
\end{TheoA}

\vskip-2pt 

A crucial property in our proof is that the Cayley graph of $\F_n$ is a tree. This suggests that similar ideas might be applicable to hyperbolic groups, which admit a tree-like Cayley graph. To be more precise, we just need that the group behaves locally like that, since for large lengths the decay of the Markov process is fast enough to produce nice estimates. That is why we avoid small loops associated to the length $\psi$ in our general framework. We illustrate this with triangular groups $\Delta_{\alpha\beta\gamma} \, = \, \langle a,b,c \, : \, a^2 = b^2 = c^2 = (ab)^\alpha = (bc)^\beta = (ca)^\gamma = e \rangle$, which are natural examples of hyperbolic groups in the Coxeter family. As finitely generated groups we may consider again the word length $| \cdot |$ ---which is conditionally negative in any infinite Coxeter group by \cite{BJS}--- and we shall keep $\P_t$ for the Poisson-like semigroup associated to the word length in any discrete group. Since relations $a^2 = b^2 = c^2 = e$ are not considered as loops, the smallest loop has length $2 \min(\alpha, \beta, \gamma)$.

\begin{TheoAA}
If $q \in 2\Z_+$, $f \in L_2(\L(\Delta_{\alpha\beta\gamma}))$ and $\min(\alpha,\beta,\gamma) \ge 8$ $$\|\mathcal{P}_t f\|_{L_q(\mathcal{L}(\Delta_{\alpha\beta\gamma}))} \le \|f\|_{L_2(\mathcal{L}(\Delta_{\alpha\beta\gamma}))} \quad \Leftrightarrow \quad t \ge \frac12 \log (q-1).$$ In particular, if $1 < p \le q < \infty$ we obtain $\|\mathcal{P}_t f\|_q \le \|f\|_p$ for $t \ge \frac{\log 3}{2} \log \frac{q-1}{p-1}$.
\end{TheoAA}

Our next family is given by the finite cyclic groups $\Z_n$. Hypercontractivity for these groups with respect to the word length is quite intriguing. Beyond Bonami's two-point inequality for $\Z_2$, only $\Z_4$ and $\Z_5$ are settled and very few is known for other values of $n$. $\Z_3$ is the simplest group with a 3-loop in its Cayley graph and the optimal time is not even conjectured. In fact, only partial results due to Andersson and Diaconis/Saloff-Coste \cite{A2,DS} are known. Andersson got optimal $L_2 \to L_q$ bounds for the $2$-truncations $\sum_{|k| \le 2} e^{-t |k|} \hat{f}(k) \lambda(k)$, with $q \in 2 \Z_+$ for $n$ odd. Since $\mathrm{diam}(\Z_n,|\cdot|) \le 2$ for $n \le 5$, the problem is still open for any $n \ge 6$. 

\begin{TheoAAA}
Given $n \ge 6$, we find $$\|\mathcal{P}_t f\|_{L_q(\mathcal{L}(\Z_n))} \le \|f\|_{L_2(\mathcal{L}(\Z_n))} \quad \Leftrightarrow \quad t \ge \frac12 \log (q-1)$$ whenever one of the following conditions hold
\begin{itemize}
\item[i)] $q \in 2\Z_+$ and $n$ is even,

\item[ii)] $q \in 2\Z_+$ and $n$ is odd with $n \ge q$.
\end{itemize}
Hence, we always have $\|\mathcal{P}_t f\|_q \le \|f\|_p$ for $1 < p \le q < \infty$ and $t \ge \frac{\log 3}{2} \log \frac{q-1}{p-1}$.
\end{TheoAAA}

Theorems A1, A2 and A3 provide by duality a nearly optimal solution to the problem opening this paper for the Markov processes on these groups associated to the word length. Namely, given $1 < p \le 2$ there exists a function $1 \le \beta_{\G}(p) \le \log 3$ with $\beta_{\G}(p) = 1$ for $p = 2, 4/3, 6/5, 8/7, \ldots$ and such that we find $$\sum_{g \in \G} e^{-2t |g|} |\widehat{f}(g)|^2 \le \|f\|_{L_p(\mathcal{L}(\G))}^2 \qquad \mbox{for} \qquad t \ge \frac{\beta_{\G}(p)}{2} \log \frac{1}{p-1}.$$ 

The proof relies on a new combinatorial method which is presented in Section \ref{Method}. It applies a priori to every pair $(\G,\psi)$ satisfying our growth/cancellation conditions provided the numerical part can be tested in a computer and the corresponding super-pathological terms can be estimated by hand. The success depends crucially on the accuracy of our estimates for the growth rate of $(\G,\psi)$ and the cancellation relations given by loops associated to the length $\psi$. The main idea is to provide an algorithm to complete squares in the expression we get for the $q$-norm of the Fourier series defining $\mathcal{P}_t f$. Indeed, the statement for $q \in 2 \Z_+$ can be rewritten for $f \ge 0$ and $t = \frac12 \log(q-1)$ as $$\sum_{\begin{subarray}{c} g_1 g_2 \cdots g_q = e \\ g_j \in \G \end{subarray}} \prod_{j=1}^q \frac{\widehat{f}(g_j)}{\sqrt{q-1}^{\psi(g_j)}} \ \le \, \sum_{h_1, h_2, \ldots, h_{\frac{q}{2}} \in \G}\, \prod_{j=1}^{q/2} |\widehat{f}(h_j)|^2.$$ The combinatorial challenge is to find a nearly optimal way to complete squares in the left hand side ---to obtain an upper bound of the form given in the right hand side--- and therefore complete the proof by a simple comparison of the coefficients in both infinite sums. As a naive starting point, we might expect that the larger the lengths $\psi(g_j)$ are, the easier will be to find an admissible way to complete squares for the term associated to $(g_1, g_2, \ldots, g_q)$. In this paper we construct a critical function $$\mu_q(\G,\psi,\cdot): \big\{ 1,2,\ldots,\frac{q}{2} \big\} \to \N$$ which allows us to decide what is \lq\lq large enough" and what is not. Thus, the terms we obtain after completing squares are divided into regular and (finitely many) pathological ones according to this function. The regular ones will be handled by means of purely combinatorial methods, while the pathological ones additionally require computer assistance to test whether our estimates are fine. Only a small subset (at least for the groups/lengths in Theorems A1-A3) of super-pathological terms fails this test, and demands finer estimates. A crucial property of the critical function to make the problem computationally solvable ---if we want to treat all possible values of $q \in 2 \Z_+$ at once--- is to be uniformly bounded. In fact, we will show that we can keep all the critical functions $\mu_q(\F_n, |\cdot|, m)$ uniformly bounded in the variables $(n,m,q)$ with $q \ge q(n)$ for some index $q(n)$ depending on the number $n$ of generators. This computational limitation is what forces us to restrict the indices $q$ in Theorem A1 ii). Much lower indices can be considered by working with $\F_3, \F_4, \ldots$ isolatedly, but this demands more computations. Also, we may produce optimal time hypercontractivity bounds in $\F_\infty$ equipped with a weighted length function, see \eqref{WeightedLength}. Additionally, as we did for free groups, we can keep the critical functions associated to all Coxeter groups uniformly bounded taking $q$ greater than or equal to some index which depends on the number of generators of the group considered. This leads to potential generalizations of Theorem A2.  

Our second contribution in this paper provides a less accurate but more general combinatorial method. It also has the advantage that it does not rely on numerical procedures. More precisely, let us call $\psi: \G \to \R_+$ a Poisson-like length whenever the following conditions hold
\begin{itemize}
\item[$\bullet$] $\psi$ is conditionally negative,

\vskip2pt

\item[$\bullet$] Spectral gap: $\sigma = \inf_{g \neq e} \psi(g) > 0$, 

\vskip2pt

\item[$\bullet$] Subadditivity: $\psi(gh) \le \psi(g) + \psi(h)$,

\vskip2pt

\item[$\bullet$] Exponential growth: $| \{ g \in \G : \psi(g) \le R \} | \le \mathrm{C} \rho^R$ for some $\rho > 1$.
\end{itemize} 
Note that we admit the existence of small loops. In particular, the word length for any finitely generated group $\G$ satisfies all the properties above automatically except perhaps the conditional negativity. Recall the general map $\S_{\psi,t}$ introduced at the beginning of the paper. Our result can be stated in full generality as follows.

\begin{TheoB}
Let $(\G,\psi)$ be a discrete group equipped with any Poisson-like length. Then, there exists a constant $\beta(\G,\psi) \ge 1$ such that the following hypercontractivity estimate holds for any $1 < p \le q < \infty$ and any $f \in L_p(\mathcal{L}(\G))$ $$\|\S_{\psi,t}f\|_q \, \le \, \|f\|_p \qquad \mbox{for all} \qquad t \ge \frac{\beta(\G,\psi)}{2 \sigma} \log \Big( \frac{q-1}{p-1} \Big).$$ Moreover, if $\psi$ is not conditionally negative, the same holds for $(p,q) = (2,4)$.
\end{TheoB}

The price to avoid the numerical dependence of our previous method is that we may not expect optimal time estimates as in Theorems A1-A3, our argument gives $\beta(\G,\psi) \le \eta \log (\rho)$ for $\rho$ large and $\eta > 1 + \sigma$. The proof of Theorem B is nevertheless much simpler and applies to a large class of groups. This includes many finitely generated groups equipped with the word length, like infinite Coxeter groups. Also groups admitting small loops in its Cayley graph ---like the discrete Heisenberg group--- for which our previous method is not efficient. On the other hand, we construct Poisson-like lengths on any discrete group admitting a finite-dimensional proper cocycle, while hypercontractivity does not occur for conditionally negative lengths in discrete groups satisfying Kazhdan property (T). Theorem B also yields apparently new estimates for non-Markovian semigroups associated to Poisson-like lengths failing to be conditionally negative. In fact, our first combinatorial method towards optimal time $L_2 \to L_q$ estimates for $q \in 2\Z_+$ also applies a priori to non-Markovian semigroups, since conditional negativity is only crucial in Gross extrapolation method for general indices $1 < p \le q < \infty$ via logarithmic Sobolev inequalities. We have avoided this more general formulation for simplicity in the exposition, see Section \ref{RFnCritical} for further details. Using our first combinatorial method we may also prove ultracontractivity bounds for arbitrary lengths which improve the trivial ones when $\psi$ admits a large concentration around $0$, see Corollary \ref{CorollaryC} and the comments after it for further details. 

To conclude, it is worth mentioning yet another approach introduced by Bakry and \'Emery \cite{BaEm} to deduce  hypercontractivity bounds for diffusion semigroups with generators satisfying the $\Gamma_2$-criterion. However, this criterion generally fails for Poisson-like processes. Moreover, a recent counterexample in \cite{JZ} shows that this approach does not necessarily work in the noncommutative setting. Instead, our combinatorial approach relies on more basic tools which still apply for the class of noncommutative Poisson-like semigroups considered in this paper. 


\vskip5pt

\noindent \textbf{Acknowledgements.} We want to thank Javier G\'omez-Serrano for helping us to optimize the computer algorithm. The authors were partially supported by ICMAT Severo Ochoa Grant SEV-2011-0087 (Spain); ERC Grant StG-256997-CZOSQP (EU); NSF Grant DMS-0901457 (USA); MINECO Grants MTM-2010-16518 \& MTM-2011-26912 and  ``Juan de la Cierva'' program (Spain).

\section{{\bf The combinatorial method}} \label{Method}

In this section we present our combinatorial method in the context of pairs $(\G,\psi)$ satisfying the growth/cancellation conditions given in the Introduction. The result below shows that Markovian semigroups associated to conditionally negative lengths have positive maximizers. In particular, it suffices to prove Theorems A1-A3 for elements in the positive cone of $L_p(\mathcal{L}(\G))$. 

\begin{lemma} \label{posmaximizers}
If $\psi$ is a conditionally negative length, we have $$\|\S_{\psi,t}\|_{p \to q} \, = \, \sup \Big\{ \|\S_{\psi,t}f\|_q \, : \ \|f\|_p=1, f \in L_p^+(\mathcal{L}(\G)) \Big\} \quad \mbox{for any} \quad t > 0.$$
\end{lemma}

\dem According to Schoenberg's theorem \cite{Sc}, we know from the properties of $\psi$ that $\S_{\psi,t}$ is a unital c.p. map. By Stinespring's factorization \cite{St}, we may find a Hilbert space $\mathcal{K}_t \supset \ell_2(\G)$ and a  $*$-homomorphism $\pi_t: \mathcal{L}(\G) \to \mathcal{B(K}_t)$, so that $\S_{\psi,t} = \mathcal{E}_t \circ \pi_t$ with $\mathcal{E}_t$ the natural conditional expectation $\mathcal{B(K}_t) \to \mathcal{B}(\ell_2(\G))$. Once this is known, the argument to find positive maximizers follows verbatim Carlen and Lieb in \cite[Theorem 3]{CL} for the CAR algebra. \fin

Moreover, it suffices to assume that $\widehat{f}(g) = \widehat{f}(g^{-1}) \ge 0$ for all $g \in \G$ when we deal with $L_2 \to L_q$ estimates for $q \in 2\Z_+$. Indeed, we may assume from Lemma \ref{posmaximizers} that $f \ge 0$, so that $$\widehat{f}(g) = \overline{\widehat{f}(g^{-1})}.$$ On the other hand, since $\tau(\lambda(g)) = \delta_{g=e}$ and $q \in 2 \Z_+$ we have 
$$\|\S_{\psi,t}f\|_q^q \, = \sum_{g_1 g_2 \cdots g_q = e} \prod_{j=1}^q \widehat{f}(g_j) e^{-t \psi(g_j)} \, \le \sum_{g_1 g_2 \cdots g_q = e} \prod_{j=1}^q |\widehat{f}(g_j)| e^{-t \psi(g_j)} = \|\S_{\psi,t} f'\|_q^q,$$ where $f' = \sum |\widehat{f}(g)| \lambda(g)$. Note that $\widehat{f'}(g) = \widehat{f'}(g^{-1}) \ge 0$ whenever $f \ge 0$ and proving hypercontractivity for $f'$ implies the same estimate for $f$, since both $f$ and $f'$ share the same $L_2$-norm. Therefore, we shall assume along this section that $f \in L_2(\L(\G))$ is a non necessarily positive operator, which admits symmetric positive Fourier coefficients. Note that this assumption is no longer valid ---neither needed--- for Gross extrapolation argument in Appendix A. 

\subsection{Notation} Given a partition $\pi$ of $\{1,2, \ldots, u\}$, define an equivalence relation on $\{1,2, \ldots, u\}$ by setting $i \, \mathcal{R}_\pi \, j$ iff both belong to the same block of $\pi$. Then, if $\k=(k_1, k_2, \ldots, k_u)\in \R^u$ we consider the partition $\pi(\k)$ determined by $i \, \mathcal{R}_{\pi(\k)} \,  j$ iff $k_i=k_j$. It will be useful below to think of $k_1,k_2, \ldots, k_u$ as colored balls which share the same color iff they belong to the same block of the partition $\pi(\k)$. Let us set $j(\k) = (j_1,j_2, \ldots, j_u)$, where $j_i$ is the number of $i$-blocks in the partition $\pi(\k)$ so that $\sum_{i=1}^u j_i i = u$. Consider also the equivalence relation on $u$-tuples $$\k \sim \k' \, \Leftrightarrow \, k_{j} = k'_{\sigma(j)} \quad \mbox{for} \quad 1 \le j \le u$$ and some permutation $\sigma \in \mathfrak{S}_u$, the symmetric group on $\{1,2, \ldots, u\}$. In other words, we will write $\k \sim \k'$ when both vectors share the same non-increasing rearrangement. Given $\k = (k_1,k_2, \ldots, k_u) \in \R^u$, let 
$$M(\k) \, = \, \big| \big\{ \k' \in \R^u : \ \k \sim \k' \big\} \big|$$ the number of all possible ways of ordering the $u$ colored balls $k_1,k_2, \ldots, k_u$. 

\begin{lemma}\label{M}
If $\k \in \R^u$ and $j(\k)=(j_1,j_2, \ldots,j_u)$, we obtain $\displaystyle M(\k) \, = \, u! \prod_{i=1}^u \frac{1}{(i!)^{j_i}}$.
\end{lemma}

\dem The partition $\pi(\k)$ consists of 
\begin{itemize}
\item $j_1$ singletons: $k_1^1, k_1^2, \ldots, k_1^{j_1}$.
\item $j_2$ pairs: $(k_2^1,k_2^1), (k_2^2,k_2^2), \ldots, (k_2^{j_2},k_2^{j_2})$.
\item  $j_i$ blocks of size $i$: $\underbrace{(k_i^1,\cdots, k_i^1)}_{i}, \ldots, \underbrace{(k_i^{j_i},\cdots, k_i^{j_i})}_{i}$.
\end{itemize}
Setting $\{i_1,i_2,\ldots,i_\ell\}=\{i\;:\; j_i \neq 0\}$, $M(\k)$ is given by the product $$\prod_{k=1}^\ell \prod_{s=1}^{j_{i_k}} \underbrace{\binom{u - \sum_{m=1}^{k-1} j_{i_m} i_m - (s-1) i_k}{i_k}}_{\mbox{place } k_{i_k}^s\mbox{'s}}.$$ After simplifying this expression, we obtain $u! / (i_1!)^{j_{i_1}} \cdots (i_\ell!)^{j_{i_\ell}}$ as desired. \fin

Given $\G$ a discrete group and $\psi: \G \to \Z_+$ an integer-valued length function on $\G$ (like the word length) write $W_k = \{g \in \G \, : \, \psi(g)=k\}$ for the set of elements of length $k$. Define $N_k = |W_k|$ (which must be finite according to our growth assumption in the Introduction) and enumerate $W_k$ by $w_k(1), w_k(2), \ldots, w_k(N_k)$. Then we set for a fixed $f \in L_2(\mathcal{L}(\G))$ with symmetric positive Fourier coefficients $$a_0 = \widehat{f}(e) \quad \mbox{and} \quad a_k(i)=\widehat{f}(w_k(i))$$ for $1\le i \le N_k$. Given $k \ge 1$, define the coefficients $$\alpha_0 = \widehat{f}(e)^2 = a_0^2 \quad \mbox{and} \quad \alpha_k = \sum_{g \in W_k} \widehat{f}(g)^2 = \sum_{i=1}^{N_k} a_k(i)^2.$$
 
\subsection{Aim of the method} \label{outline proof}

Consider a discrete group $\G$ which comes equipped with a conditionally negative subadditivite length $\psi: \G \to \R_+$ admitting a spectral gap $\sigma > 0$ and satisfying our growth/cancellation conditions
\begin{itemize}
\item[$\bullet$] There exists $r> 0$ such that $\sum_{g \in \G} r^{\psi(g)} < \infty$,

\vskip2pt

\item[$\bullet$] $g_1 g_2 \cdots g_m \neq e$ whenever $\psi(g_j)=\sigma$, $g_j \neq g_{j+1}^{-1}$ and $m$ small.
\end{itemize}
Assume for simplicity that $\psi: \G \to \Z_+$ and $\sigma = 1$. Take $r = e^{-t}$ and define $$\T_{\psi,r} f  \, = \, \sum_{g \in \G} r^{\psi(g)} \widehat{f}(g) \lambda(g) \, = \, \sum_{g \in \G}^{\null} e^{-t \psi(g)} \widehat{f}(g) \lambda(g) \, = \, \S_{\psi,t} f.$$ The basic idea of our method is the following. As explained above, we may fix $f \in L_2(\mathcal{L}(\G))$ with symmetric positive Fourier coefficients and $q$ an even integer greater than or equal to 4, since $q=2$ follows from Plancherel theorem. Set $q = 2s$, $$L(s) = \Big\{ \k = (k_1, k_2, \ldots, k_s) \in \N^s \; : \; k_1 \ge k_2 \ge \cdots \ge k_s \ge 0 \Big\}$$ and $C_q^{\mathrm{right}}[\k] = M(\k)$ for any $\k \in L(s)$. If $\alpha_{\k} = \alpha_{k_1} \alpha_{k_2} \cdots \alpha_{k_s}$, we have $$\|f\|_2^q \, = \, \Big( \sum_{g \in \G} \widehat{f}(g)^2 \Big)^s \, = \, \Big( \sum_{k \ge 0} \alpha_k \Big)^s \, = \, \sum_{\k \in L(s)} C_q^{\mathrm{right}}[\k] \, \alpha_{\k}.$$ Moreover, since $\tau(\lambda(g)) = \delta_{g=e}$ and $q$ is an even integer, we also find that 
\begin{eqnarray} \label{sumas-su}
\hskip10pt \|\T_{\psi,r}f \|_q^q & = & \sum_{u=0}^q \, \sum_{\substack{g_1 g_2 \cdots g_{2s} = e \\ |\{j\;:\; g_j \neq e\}|=u}} \Big(\prod_{j=1}^{2s} \widehat{f}(g_j)\Big) \, r^{\sum_{j=1}^{2s} \psi(g_j)} \\ \nonumber & = & \widehat{f}(e)^q + \sum_{u=1}^q \binom{q}{u} \widehat{f}(e)^{q-u} \underbrace{\sum_{\substack{g_1 g_2 \cdots g_{u} = e \\ g_j \neq e}} \Big( \prod_{j=1}^u \widehat{f}(g_j) \Big) r^{\sum_{j=1}^u \psi(g_j)}}_{s_u(r)}. 
\end{eqnarray}
We have also used our condition on $f$, which gives $\widehat{f}(g) \! = \! \widehat{f}(g^{-1}) \! \ge \! 0$ for all $g$. Since $\widehat{f}(e)^q = \alpha_0^s$ and $s_1(r)=0$, optimal time hypercontractivity will follow from 
\begin{equation} \label{AIM}
\left\{ \begin{array}{c}
\widehat{f}(e)^{q-u} s_u(r) \ \le \ \displaystyle \sum_{\k \in L(s) \setminus \{\underline{0}\}}^{\null} C_{q,s_u}^{\mathrm{left}}[\k](r) \, \alpha_{\k}, \\ [15pt] C_q^{\mathrm{left}}[\k](r) := \displaystyle \sum_{u=2}^q \binom{q}{u} C_{q,s_u}^{\mathrm{left}}[\k](r) \le C_q^{\mathrm{right}}[\k] \ \, \mbox{with} \, \ \k \in L(s) \setminus \{\underline{0}\}, \end{array} \right\}
\end{equation}
for any $0 \le r \le \frac{1}{\sqrt{q-1}}$. The proof of \eqref{AIM} is outlined in the following paragraphs.

\subsection{Admissible lengths} Note that $$s_u(r) \, = \, \sum_{\l \in L(u)} \sum_{\substack{g_1 g_2 \cdots g_{u} = e \\ g_j \neq e \\ (\psi(g_1), \psi(g_2), \ldots, \psi(g_u)) \sim \l}}^{\null} \Big( \prod_{j=1}^u \widehat{f}(g_j) \Big) \, r^{\sum_{j=1}^u \psi(g_j)} \, =: \, \sum_{\l \in L(u)} s_u[\l](r).$$ The following is a simple consequence of the symmetry and subadditivity of $\psi$.

\begin{lemma}
If $g_1 g_2 \cdots g_u = e$, we find that $\psi(g_i) \le \sum_{j \neq i} \psi(g_j)$ for all $1 \le i \le u$.
\end{lemma}

The previous result allows to refine the set of admissible lengths $\psi(g_j)$ which may appear in the sums $s_u(r)$. This leads naturally to the following set of admissible length $u$-tuples with $2 \le u \le q$ $$\mathrm{Adm}_u \, = \, \Big\{ \l \in L(u) \; : \; \ell_1 \ge \ell_2 \ge \ldots \ge \ell_u \ge 1, \; \ell_1 \le \sum_{j=2}^{u} \ell_j \Big\}.$$ Of course, we clearly have the refined identity $s_u(r) \, = \, \sum_{\l \in \mathrm{Adm}_u}^{\null} s_u[\l](r)$.

\subsection{Completing squares I} As we shall justify later, the sums $s_2(r)$ and $s_3(r)$ must be estimated apart, so we just focus on the general method to estimate the sums $s_u(r)$ for $4 \le u \le q$. The terms $s_u[\l](r)$ for $\l$ in a certain (finite) exceptional set $B_u \subset \mathrm{Adm}_u$ to be defined below will contribute to the so-called super-pathological terms and also require specific methods. For technical reasons that will appear only at the end of the proofs of Theorems A1-A3, we need to introduce here these exceptional sets $B_u$ to be completely rigorous. However, we suggest the reader to assume $B_u = \emptyset$ in a first reading. The remaining terms $s_u[\l](r)$ require two different ways of completing squares, which we now describe. If $u = 2m$ $(2 \le m \le s)$ is even, we set 

\vskip-10pt

\begin{eqnarray*} 
C_{2m} \!\!\! & = & \!\!\! \Big\{ \x=(\xi_1,\xi_2,\ldots,\xi_m)\;:\; \xi_j\in \{2j-1,2j\} \; \ \mbox{for all} \ 1 \le j \le m \Big\}, \\ 
\Lambda_{2m} \!\!\! & = & \!\!\! \Big\{ \lambda = (\x,\l,\d,\ibar) \; : \; \x \in C_{2m}, \; \l \in \mathrm{Adm}_{2m} \!\! \setminus \!\! B_{2m},\; \d \sim \l, \; \prod_{1 \le j \le 2m}^{\rightarrow} w_{d_j}(i_j) = e \Big\}.
\end{eqnarray*}
The arrow in the product means that the order is $j$-increasing, so that we get $w_{d_1}(i_1) w_{d_2}(i_2) \cdots w_{d_{2m}}(i_{2m})$. It is clear that $\ibar = (i_1, i_2, \ldots, i_{2m})$ with $1 \le i_j \le N_{d_j}$ for all $1 \le j \le 2m$. Given any $\d = (d_1, d_2, \ldots, d_{2m})$ such that $\d \sim \l$, we may pick a (non-unique) permutation $\sigma_{\d} \in \mathfrak{S}_{2m}$ such that $d_j = \ell_{\sigma_{\d}(j)}$ for all $1 \le j \le 2m$. Once we have fixed $\sigma_{\d}$ for each $\d$ and $\lambda = (\x,\l,\d,\ibar) \in \Lambda_{2m}$, we set $|\l| = \summ_j \ell_j$ and $$\gamma_\lambda \, = \, \frac{1}{2^m} \prod_{j=1}^m a_{\ell_{\xi_j}}(i_{\sigma_{\d}^{-1}(\xi_j)})^2, \quad \quad \nu_\lambda(r) \, = \, r^{|\l|}.$$

\begin{lemma} \label{EvenLemma}
If $2 \le m \le s$, we find
\begin{eqnarray*}
s_{2m}(r) & \le & \sum_{\lambda \in \Lambda_{2m}} \gamma_\lambda \nu_\lambda(r) + \sum_{\l \in B_{2m}} s_{2m}[\l](r) \\ & = & \sum_{\k \in L_m(s)} \sum_{\lambda \in \Lambda_{2m}[\k]} \gamma_\lambda \nu_\lambda(r) + \sum_{\l \in B_{2m}} s_{2m}[\l](r),
\end{eqnarray*} 
where $L_m(s) = \{ \k \in L(s) \, : \, k_m \neq 0 = k_{m+1} \}$ and $ \Lambda_{2m}[\k] = \{ \lambda \! \in \! \Lambda_{2m} \, : \, \ell_{\xi_j} = k_j \}$.
\end{lemma}

\dem The last identity is clear, while the first inequality follows from
\begin{eqnarray*}
\lefteqn{\hskip-20pt \sum_{\l \in \mathrm{Adm}_{2m} \! \setminus \! B_{2m}}^{\null} s_{2m}[\l](r)} \\ \!\!\! & = & \!\!\! 
\sum_{\l \in \mathrm{Adm}_{2m} \! \setminus \! B_{2m}} \ \sum_{\d \sim \l}^{\null} \sum_{\substack{\ibar \\ w_{d_1}(i_1) \cdots w_{d_{2m}}(i_{2m}) = e}} \Big( \prod_{j=1}^{2m} a_{d_j}(i_j) \Big) \, r^{|\l|} \\ \!\!\! & = & \!\!\! \sum_{\l \in \mathrm{Adm}_{2m} \! \setminus \! B_{2m}} \ \sum_{\d \sim \l}^{\null} \sum_{\substack{\ibar \\ w_{d_1}(i_1) \cdots w_{d_{2m}}(i_{2m}) = e}} \Big( \prod_{j=1}^{2m} a_{\ell_j}(i_{\sigma_{\d}^{-1}(j)}) \Big) \, r^{|\l|}.
\end{eqnarray*}  
This yields the desired estimate $\sum_{\lambda \in \Lambda_{2m}} \gamma_\lambda \nu_\lambda(r)$, since completing squares gives  

\vskip-5pt

\null \hfill \hskip10pt $\begin{array}{rcl}
\displaystyle \prod_{j=1}^{2m} a_{\ell_j}(i_{\sigma_{\d}^{-1}(j)}) & \le & \displaystyle \prod_{j=1}^m  \Big[ \frac{a_{\ell_{2j-1}}(i_{\sigma_{\d}^{-1}(2j-1)})^2 + a_{\ell_{2j}}(i_{\sigma_{\d}^{-1}(2j)})^2}{2} \Big] \\ & = & \displaystyle \sum_{\x \in C_{2m}} \ \prod_{j=1}^m \ \frac{a_{\ell_{\xi_j}}(i_{\sigma_{\d}^{-1}(\xi_j)})^2}{2} \ \ = \ \sum_{\x \in C_{2m}} \gamma_\lambda. \end{array}$ \hfill \hskip10pt $\begin{array}{r} \null \\ [22pt] \hskip10pt \square \end{array}$

\vskip5pt

\noindent To consider the case where $u = 2m-1$ is odd, we set for $3 \le m \le s$
\begin{eqnarray*} 
\Lambda'_{2m} \!\!\! & = & \!\!\! \Big\{ \lambda = (\x,\l,\d,\ibar) \; : \; \x \in C_{2m}, \; \l \in \mathrm{Adm}'_{2m} \!\! \setminus \!\! B'_{2m},\; \d \sim \l, \; \prod_{1 \le j \le 2m}^{\rightarrow} w_{d_j}(i_j) = e \Big\},
\end{eqnarray*}
where $\mathrm{Adm}_{2m}'$ and $B'_{2m}$ are the sets of tuples $(\ell_1, \ell_2, \ldots, \ell_{2m-1},0)$ satisfying that $(\ell_1, \ell_2, \ldots, \ell_{2m-1})$ belongs to $\mathrm{Adm}_{2m-1}$ and $B_{2m-1}$ respectively. Also, given any $\lambda = (\x,\l,\d,\ibar) \in \Lambda_{2m}'$, we define the coefficients $$\gamma_\lambda' \, = \, \frac{1}{2m} \Big( \frac{1}{2^m} \prod_{j=1}^m a_{\ell_{\xi_j}}(i_{\sigma_{\d}^{-1}(\xi_j)})^2 \Big).$$ Now we state and prove the analog of Lemma \ref{EvenLemma} in the odd case.

\begin{lemma} \label{OddLemma}
If $3 \le m \le s$, we find 
\begin{eqnarray*}
\widehat{f}(e) s_{2m-1}(r) & \le & \sum_{\lambda \in \Lambda'_{2m}} \gamma'_\lambda \nu_\lambda(r) + \sum_{\l \in B_{2m-1}} \widehat{f}(e) s_{2m-1}[\l](r) \\ & = & \sum_{\k \in L_m'(s)} \sum_{\lambda \in \Lambda'_{2m}[\k]} \gamma'_\lambda \nu_\lambda(r) + \sum_{\l \in B_{2m-1}} \widehat{f}(e) s_{2m-1}[\l](r),
\end{eqnarray*} 
where $L_m'(s) = L_{m-1}(s) \cup L_m(s)$ and $\Lambda'_{2m}[\k] = \big\{ \lambda \! \in \! \Lambda'_{2m} \, : \, \ell_{\xi_j} = k_j, \, 1 \le j \le m \big\}$.
\end{lemma}

\dem Again, the last identity is straightforward. Namely, the only difference is that for $\lambda = (\x,\l,\d,\ibar) \in \Lambda_{2m}'$ we have that $\l \in \mathrm{Adm}_{2m}' \! \setminus \! B_{2m}'$. This means that $k_m = \ell_{\xi_m} \in \{\ell_{2m-1}, \ell_{2m} \}$ could be $0$ or not, so that $\Lambda_{2m}'$ splits as the disjoint union of $\Lambda'_{2m}[\k]$'s over $L'_m(s) = L_{m-1}(s) \cup L_m(s)$. The  inequality is very similar to the one in Lemma \ref{EvenLemma}. Indeed, it suffices to note that \\ \vskip-20pt  \null \hfill $\displaystyle \sum_{\substack{g_1g_2 \cdots g_{2m-1}=e \\ g_j \neq e}} \widehat{f}(e) \Big( \prod_{j=1}^{2m-1} \widehat{f}(g_j) \, r^{\psi(g_j)} \Big) \, = \, \frac{1}{2m} \sum_{\substack{g_1g_2 \cdots g_{2m} = e\\ \exists ! \; j / g_j = e}} \Big( \prod_{j=1}^{2m} \widehat{f}(g_j) \, r^{\psi(g_j)} \Big).$ \hfill $\square$ 

\subsection{A decomposition of $s_u(r)$} We now estimate $s_u(r)$ for $4 \le u \le q$ by three sums of regular, pathological and super-pathological terms respectively. The crucial decomposition is given by a partition of the set $L_m(s)$ into two subsets, which is determined by a critical function $\mu_q(\G, \psi, \cdot \, ): \big\{1,2, \ldots, \frac{q}{2} \big\} \to \N$. As pointed out in the Introduction, the critical function will quantify what is \lq\lq large enough" for the exponent $\sum_j \psi(g_j)$ in $s_u(r)$ given by \eqref{sumas-su}. Of course, we may not expect to quantify anything before applying the method itself! This means that the critical function arises necessarily a posteriori, and its construction must be done case by case. In the proofs of Theorems A1-A3 we need to begin by providing the explicit critical functions, but only at the end of the argument the reader will be able to fully justify our choice for this function. Let us define 
\begin{eqnarray*}
\mathrm{Reg}_m(s) & = & \Big\{ \k = (k_1, k_2, \ldots, k_m, \underline{0}) \in L_m(s) \; : \; |\k| \ge m + \mu_q(\G,\psi,m) \Big\}, \\ \mathrm{Pat}_m(s) & = & \Big\{ \k = (k_1, k_2, \ldots, k_m, \underline{0}) \in L_m(s) \; : \; |\k| < m + \mu_q(\G,\psi,m) \Big\}. 
\end{eqnarray*}
It is crucial to note that the sets $\mathrm{Pat}_m(s)$ of pathological terms are finite. As in the previous paragraph, it will be instrumental to provide slightly modified sets to deal with the odd case. Namely, we set $\mathrm{Reg}'_m(s) = \mathrm{Reg}_{m-1}(s) \cup \mathrm{Reg}_m(s)$ and $\mathrm{Pat}'_m(s) = \mathrm{Pat}_{m-1}(s) \cup \mathrm{Pat}_m(s)$. Then we may estimate the sums $s_u(r)$ as follows 
\begin{equation} \label{DecompositionRPS1}
\begin{array}{rcl}
s_{2m}(r) & \le & \mathbf{R}_{2m}(\Lambda,r) + \mathbf{P}_{2m}(r) + \mathbf{S}_{2m}(r) \quad (2 \le m \le s), \\ \widehat{f}(e) s_{2m-1}(r) & \le & \mathbf{R}'_{2m}(\Lambda,r) + \mathbf{P}'_{2m}(r) + \mathbf{S}'_{2m}(r) \quad (3 \le m \le s), 
\end{array}
\end{equation}
where 
\begin{eqnarray*}
\big( \mathbf{R}_{2m}(\Lambda,r), \mathbf{R}'_{2m}(\Lambda,r) \big) & = & \Big( \sum^{\null}_{\substack{\k \in \mathrm{Reg}_m(s) \\ \lambda \in \Lambda_{2m}[\k]}} \gamma_\lambda \nu_\lambda(r), \sum^{\null}_{\substack{\k \in \mathrm{Reg}'_m(s) \\ \lambda \in \Lambda'_{2m}[\k]}} \gamma'_\lambda \nu_\lambda(r) \Big), \\ \big( \mathbf{P}_{2m}(r), \mathbf{P}'_{2m}(r) \big) & = & \Big( \sum_{\substack{\k \in \mathrm{Pat}_m(s) \\ \lambda \in \Lambda_{2m}[\k]}} \gamma_\lambda \nu_\lambda(r), \, \sum^{\null}_{\substack{\k \in \mathrm{Pat}'_m(s) \\ \lambda \in \Lambda'_{2m}[\k]}} \, \gamma'_\lambda \nu_\lambda(r) \Big), \\ [3pt] \big( \mathbf{S}_{2m}(r), \mathbf{S}'_{2m}(r) \big) & = & \Big( \sum_{\substack{\l \in B_{2m}}} s_{2m}[\l](r), \sum_{\substack{\l \in B_{2m-1}}} \widehat{f}(e) s_{2m-1}[\l](r) \Big).
\end{eqnarray*}

\subsection{Completing squares II} Our second way of completing squares is perhaps less accurate but definitely more symmetric. This will be crucial in certain estimates below. If $m \ge 2$ and $g_1, g_2, \ldots, g_{2m} \in \G$ (allowing repetitions), we define the equivalence relation $\mathcal{R}_{g_1,\ldots,g_{2m}}$ on $$\mathcal{M}_m \, = \, \Big\{ \z = (\zeta_1, \zeta_2, \ldots, \zeta_m) \in \N^m \, : \ 1 \le \zeta_1 < \zeta_2 < \cdots < \zeta_m \le 2m \Big\}$$ by setting $$\z \ \mathcal{R}_{g_1,\ldots,g_{2m}} \, \z' \, \Leftrightarrow \, g_{\zeta_j} = g_{\zeta'_{\sigma(j)}} \mbox{ for some } \sigma \in \mathfrak{S}_m \, \Leftrightarrow \, (g_{\zeta_j}) \sim (g_{\zeta_j'}).$$ Let $\stackrel{\circ}{\z}$ denote the class of $\z$ in the quotient space $\mathcal{M}_m(g_1,\ldots,g_{2m}) = \mathcal{M}_m/\mathcal{R}_{g_1,\ldots,g_{2m}}$. 

\vskip10pt 

\begin{lemma}\label{complsq2}
If $m \ge 2$ and $g_1,g_2, \ldots,g_{2m} \in \G$ $$\prod_{j=1}^{2m} \widehat{f}(g_j) \, \le \, \sum_{\stackrel{\circ}{\z} \in \mathcal{M}_m(g_1,\ldots,g_{2m})} \prod_{j=1}^m \widehat{f}(g_{\zeta_j})^2.$$
\end{lemma}

\vskip5pt

\dem Completing squares as in Lemma \ref{EvenLemma} gives
$$\prod_{j=1}^{2m} \widehat{f}(g_j) \le \frac{1}{2^m} \sum_{\x \in C_{2m}} \prod_{j=1}^m \widehat{f}(g_{\xi_j})^2 \le \sup_{\z \in \mathcal{M}_m} \prod_{j=1}^m \widehat{f}(g_{\zeta_j})^2 \le \hskip-15pt \sum_{\stackrel{\circ}{\z} \in \mathcal{M}_m(g_1,\ldots,g_{2m})} \prod_{j=1}^m \widehat{f}(g_{\zeta_j})^2. \hskip10pt \square$$


Since we have already completed squares in Lemmas \ref{EvenLemma} and \ref{OddLemma}, it is necessary to show that the way of pairing stated in Lemma \ref{complsq2} is consistent with our previous estimates, which involve some kind of Fubini argument. Namely, our goal is to estimate $\mathbf{R}_{2m}(\Lambda,r)$ and $\mathbf{R}_{2m}'(\Lambda,r)$ by other regular sums. If $2 \le m \le s$, we define in the even case $$\Delta_{2m} = \Big\{ \delta = (g_1,\ldots, g_{2m}, \stackrel{\circ}{\z}) : \, 
g_1 \cdots g_{2m} = e, \, g_j \neq e, \, \stackrel{\circ}{\z} \in \mathcal{M}_m(g_1,\ldots,g_{2m}) \Big\}.$$
If $\delta \in \Delta_{2m}$ and $\k = (k_1, k_2, \ldots, k_m, \underline{0}) \in L_m(s)$, we also set $$\big( \gamma_\delta, \nu_\delta(r) \big) \, = \, \Big( \prod_{j=1}^m \widehat{f}(g_{\zeta_j})^2, r^{\sum_{j=1}^{2m} \psi(g_j)} \Big),$$ $$\Delta_{2m}[\k] = \Big\{ \delta \in \Delta_{2m} \, : \ (\psi(g_{\zeta_1}),\ldots, \psi(g_{\zeta_m})) \sim (k_1,\ldots,k_m) \Big\}.$$ Note that $\gamma_\delta$ and $\Delta_{2m}[\k]$ are well defined since $\z \, \mathcal{R}_{g_1,\ldots,g_{2m}} \, \z'$ iff $(g_{\zeta_j}) \sim (g_{\zeta_j'})$. By the usual modifications, we may also consider the corresponding sets and coefficients in the odd case. Given $3 \le m \le s$, define $$\Delta_{2m}' = \Big\{ \delta = (g_1,\ldots, g_{2m}, \stackrel{\circ}{\z}) : \, g_1 \cdots g_{2m} = e, \, \exists ! \; j \, \mbox{s.t.} \, g_j = e, \, \stackrel{\circ}{\z} \in \mathcal{M}_m(g_1,\ldots,g_{2m}) \Big\}.$$ \vskip5pt \noindent If $\delta \in \Delta_{2m}'$ and $\k = (k_1, k_2, \ldots, k_m, \underline{0}) \in L_m'(s) = L_{m-1}(s) \cup L_m(s)$, set $$\big( \gamma'_\delta, \nu_\delta(r) \big) \, = \, \Big( \frac{1}{2m} \prod_{j=1}^m \widehat{f}(g_{\zeta_j})^2, r^{\sum_{j=1}^{2m} \psi(g_j)} \Big),$$ $$\Delta'_{2m}[\k] \, = \, \Big\{ \delta \in \Delta'_{2m} \; : \; (\psi(g_{\zeta_1}), \ldots, \psi(g_{\zeta_m})) \sim (k_1,\ldots,k_m) \Big\}.$$

\vskip5pt 

\begin{lemma}\label{method2 better than method1}
Given $s \ge 2$, we find 
\begin{itemize}
\item[i)] If $2 \le m \le s$ and $\k \in L_m(s)$, we obtain $$\sum_{\lambda \in \Lambda_{2m}[\k]} \gamma_\lambda \nu_\lambda(r) \, \le \, \sum_{\delta \in \Delta_{2m}[\k]} \gamma_\delta \nu_\delta(r).$$

\item[ii)] If $3 \le m \le s$ and $\k \in L_m'(s) = L_{m-1}(s) \cup L_m(s)$ $$\sum_{\lambda \in \Lambda'_{2m}[\k]} \gamma'_\lambda \nu_\lambda(r) \, \le \, \sum_{\delta \in \Delta'_{2m}[\k]} \gamma'_\delta \nu_\delta(r).$$
\end{itemize}
\end{lemma}

\dem Consider the map $$T_{2m,\k}: (\x,\l,\d,\ibar) \in \Lambda_{2m}[\k] \mapsto (g_1, \ldots, g_{2m}, \stackrel{\circ}{\z}) \in \Delta_{2m}[\k]$$ defined by $g_j = w_{d_j}(i_j)$ and $\z$ the non-decreasing rearrangement of $\underline{\eta} = (\eta_1,\ldots, \eta_m)$ with $\sigma_{\d}(\eta_j) = \xi_j$. Let us show that this map takes values in $\Delta_{2m}[\k]$. The conditions $g_1\cdots g_{2m} = w_{d_1}(i_1)\cdots w_{d_{2m}}(i_{2m})=e$ with $g_j \neq e$ and $\z \in \M_m$ are clear. The identity $(\psi(g_{\zeta_1}),\ldots,\psi(g_{\zeta_m})) \sim (k_1,\ldots,k_m)$ follows from 
\begin{eqnarray*}
(\psi(g_{\zeta_j}))_{1\leq j\leq m} & \sim & (\psi(g_{\eta_j}))_{1\leq j\leq m} \ = \ (d_{\eta_j})_{1 \le j \le m} \\ & = & (d_{\sigma_{\d}^{-1}(\xi_j)})_{1 \le j \le m} \ = \ (\ell_{\xi_j})_{1 \le j \le m} \ = \ (k_j)_{1 \le j \le m}.
\end{eqnarray*}
Define now an equivalence relation $\mathcal{R}_{2m,\k}$ on $\Lambda_{2m}[\k]$ by $$\lambda_1 \, \mathcal{R}_{2m,\k} \, \lambda_2 \ \Leftrightarrow \ T_{2m,\k}(\lambda_1) = T_{2m,\k} (\lambda_2).$$ Set $\widetilde{\Lambda}_{2m}[\k] = \Lambda_{2m}[\k] / \mathcal{R}_{2m,\k}$ and $\widetilde{T}_{2m,\k}: \widetilde{\lambda} \in \widetilde{\Lambda}_{2m}[\k] \mapsto T_{2m,\k}(\lambda) \in \Delta_{2m}[\k]$, which is clearly an injective map. Then, our first assertion (even case) follows from the following two claims
\begin{itemize} 
\item[a)] $\displaystyle \big| \big\{ \lambda \in \Lambda_{2m}[\k] \, : \ \lambda \in \widetilde{\lambda} \big\} \big| \, \le \, 2^m$ for any $\widetilde{\lambda} \in \widetilde{\Lambda}_{2m}[\k]$,

\item[b)] $\displaystyle \big( \gamma_\lambda, \nu_\lambda(r) \big) \, = \, \Big( \frac{1}{2^m} \gamma_{T_{2m,\k}(\lambda)} ,\nu_{T_{2m,\k}(\lambda)}(r) \Big)$. 
\end{itemize}
Indeed, assuming these assertions we obtain 
\begin{eqnarray*}
\sum_{\lambda\in \Lambda_{2m}[\k]} \gamma_\lambda \nu_\lambda(r) & = & \sum_{\widetilde{\lambda} \in \widetilde{\Lambda}_{2m}[\k]} \sum_{\lambda \in \widetilde{\lambda}} \gamma_{\lambda} \nu_{\lambda}(r) \\ & = & \frac{1}{2^m} \sum_{\widetilde{\lambda} \in \widetilde{\Lambda}_{2m}[\k]} \sum_{\lambda \in \widetilde{\lambda}} \gamma_{\widetilde{T}_{2m,\k}(\widetilde{\lambda})} \nu_{\widetilde{T}_{2m,\k}(\widetilde{\lambda})}(r) \\ & \le & \sum_{\widetilde{\lambda} \in \widetilde{\Lambda}_{2m}[\k]} \gamma_{\widetilde{T}_{2m,\k}(\widetilde{\lambda})} \nu_{\widetilde{T}_{2m,\k}(\widetilde{\lambda})}(r) \ \ \le \ \sum_{\delta \in \Delta_{2m}[\k]} \gamma_\delta \nu_\delta(r).
\end{eqnarray*}
Since the odd case is similar, it suffices to prove our claims. To prove a), note that $\lambda_1 \, \mathcal{R}_{2m,\k} \, \lambda_2$ implies $(\l_1,\d_1,\ibar_1) = (\l_2, \d_2, \ibar_2)$. Hence, the cardinal of each equivalence class is dominated by $2^m$ possible choices of $\x$. The proof of b) follows from
\begin{eqnarray*}
\gamma_{T_{2m,\k}(\lambda)} \!\!\! & = & \!\!\! \prod_{j=1}^m \widehat{f}(g_{\zeta_j})^2 = \prod_{j=1}^m \widehat{f}(g_{\eta_j})^2 = \prod_{j=1}^m \widehat{f}(w_{d_{\eta_j}}(i_{\eta_j}))^2 \\ \!\!\! & = & \!\!\! \prod_{j=1}^m \widehat{f}(w_{d_{\sigma_{\d}^{-1}(\xi_j)}}(i_{\sigma_{\d}^{-1}(\xi_j)}))^2 = \prod_{j=1}^m \widehat{f}(w_{\ell_{\xi_j}}(i_{\sigma_{\d}^{-1}(\xi_j)}))^2 = 2^m \gamma_\lambda
\end{eqnarray*}
and $\nu_{T_{2m,\k}(\lambda)}(r) = r^{\sum_{j=1}^{2m} \psi(g_j)}=r^{|\d|}=r^{|\l|} = \nu_\lambda(r)$. This completes the proof. \fin

\begin{remark}
\emph{If we set} $$\big( \mathbf{R}_{2m}(\Delta,r), \mathbf{R}'_{2m}(\Delta,r) \big) \, = \, \Big( \sum^{\null}_{\substack{\k \in \mathrm{Reg}_m(s) \\ \delta \in \Delta_{2m}[\k]}} \gamma_\delta \nu_\delta(r), \sum^{\null}_{\substack{\k \in \mathrm{Reg}'_m(s) \\ \delta \in \Delta'_{2m}[\k]}} \gamma'_\delta \nu_\delta(r) \Big),$$ \emph{then Lemma \ref{method2 better than method1} together with \eqref{DecompositionRPS1} yields}
\begin{equation} \label{DecompositionRPS2}
\begin{array}{rcl}
s_{2m}(r) & \le & \mathbf{R}_{2m}(\Delta,r) + \mathbf{P}_{2m}(r) + \mathbf{S}_{2m}(r) \quad (2 \le m \le s), \\ \widehat{f}(e) s_{2m-1}(r) & \le & \mathbf{R}'_{2m}(\Delta,r) + \mathbf{P}'_{2m}(r) + \mathbf{S}'_{2m}(r) \quad (3 \le m \le s). 
\end{array}
\end{equation}
\end{remark}

\subsection{Analysis of both approaches} \label{LambdavsDelta}

Needing two ways of completing squares requires an explanation. Our first approach is more accurate, but definitely less symmetric and the price is that the critical function $\mu_q(\G,\psi,\cdot \, )$ that we would obtain with this method is not even bounded, since the estimates needed to construct the critical function are less precise in the absence of symmetry. Thus, we could not use a computer to identify the super-pathological terms among the pathological ones. Hence, our first method will be reserved for isolated cases. The second method is more symmetric and yields a uniformly bounded critical function $\mu_q(\G,\psi,\cdot \, )$, but it can only be used under the growth condition $\sum r^{\psi(g)} < \infty$ for all $r$ strictly smaller than some index $R(\G,\psi)$. Since the expected optimal $r$ for $L_2 \to L_q$ hypercontractivity is given by $(q-1)r^2=1$ we need to have $$q \, > \, 1 + \frac{1}{R(\G,\psi)^2} \, =: \, q(\G,\psi).$$ We will use our second approach to estimate the $\Delta$-regular sums for these $q$'s. On the other hand, our first (more accurate) way of completing squares will be used for the finitely many pathological terms and also for the regular terms associated to the finitely many $q$'s below the critical index $q(\G,\psi)$. We refer to Paragraph \ref{RFnCritical} for an illustration in $\F_n$ of the behavior of the critical function with both approaches. 

\subsection{$\Lambda$-estimates} \label{Lambdaest}

Our $\Lambda$-estimates will be used for the sums $\mathbf{R}_{2m}(\Lambda,r)$/$ \mathbf{R}_{2m}'(\Lambda,r)$ and $q$ less than or equal to the critical index. We will also use the same estimates below for the pathological sums $\mathbf{P}_{2m}(r)$. Recalling the definition of $C_{2m}$ as the set of $\x = (\xi_1, \xi_2, \ldots, \xi_m)$ such that $\xi_j \in \{2j-1,2j\}$, we define $\x^{\star}  = (\xi_1^\star, \xi_2^\star, \ldots, \xi_m^\star)$ by $$\big\{ \xi_j, \xi_j^\star \big\} \, = \, \big\{ 2j-1,2j \big\} \quad \mbox{for all} \quad 1 \le j \le m.$$ 

\begin{proposition} \label{Lambda-est}
Given $0 \le r \le 1$, we find 
\begin{itemize}
\item[i)] If $2 \le m \le s$ and $\k \in L_m(s)$, we obtain 
$$\hskip35pt \alpha_0^{s-m} \sum_{\lambda \in \Lambda_{2m}[\k]} \gamma_\lambda \nu_\lambda(r)  \le  \frac{1}{2^m} \Big[ \sum_{\x \in C_{2m}} \sum_{\begin{subarray}{c} \l \in \mathrm{Adm}_{2m} \setminus B_{2m} \\ \ell_{\xi_j} = k_j \end{subarray}}  M(\l)\Big( \prod_{j=2}^{m} N_{\ell_{\xi_j^\star}}\Big)r^{|\l|} \Big]   \alpha_{\k}.$$

\item[ii)] If $3\leq m \leq s$ and $\k \in L_m(s)\cup L_{m-1}(s)$, we obtain 
$$\hskip20pt \alpha_0^{s-m} \sum_{\lambda \in \Lambda'_{2m}[\k]} \gamma'_\lambda \nu_\lambda(r)  \le  \frac{1}{2^m2m} \Big[ \sum_{\x \in C_{2m}} \sum_{\begin{subarray}{c} \l \in \mathrm{Adm}_{2m}' \setminus B_{2m}' \\ \ell_{\xi_j} = k_j \end{subarray}}  M(\l)\Big( \prod_{j=2}^{m} N_{\ell_{\xi_j^\star}}\Big)r^{|\l|} \Big] \alpha_{\k}.$$
\end{itemize}
\end{proposition}

\dem We start with the even case i). Given $\k \in L_m(s)$, we may consider any triple $(\x,\l,\d)$ with $\x \in C_{2m}$, $\l \in \mathrm{Adm}_{2m} \setminus B_{2m}$ satisfying $\ell_{\xi_j} = k_j$ and $\d \sim \l$. Then we define the sets
$\Lambda_{2m}[\k,\x,\l,\d] \, = \, \{ \lambda = (\x,\l,\d,\ibar) \in \Lambda_{2m}[\k] \}$ and claim 
\begin{equation}\label{claim-patho-even}
\alpha_0^{s-m} \sum_{\lambda \in \Lambda_{2m}[\k,\x,\l,\d]} \gamma_\lambda \, \le \, \frac{1}{2^m} \Big(\prod_{j=2}^{m}N_{\ell_{\xi_j^\star}} \Big)\alpha_{\k}.
\end{equation}
Recalling that $\gamma_\lambda = 2^{-m} \prod_j a_{\ell_{\xi_j}}(i_{\sigma_{\d}^{-1}(\xi_j)})^2$, we can make vary the $2m-1$ indices $$\begin{array}{rcl} i_{\sigma_{\d}^{-1}(\xi_j)} \in \big\{ 1,2, \ldots, N_{\ell_{\xi_j}} \big\} & \mbox{for} & 1 \le j \le m, \\ i_{\sigma_{\d}^{-1}(\xi_j^\star)} \in \big\{ 1,2, \ldots, N_{\ell_{\xi_j^\star}} \big\} & \mbox{for} & 2 \le j \le m, \end{array}$$ in the sum above. Namely, the last summation index is entirely determined by the constraint $w_{d_1}(i_1) w_{d_2}(i_2) \cdots w_{d_{2m}}(i_{2m})=e$. Such an index could not exist, but in that case we get zero and the estimate below still holds true $$\sum_{\lambda \in \Lambda_{2m}[\k,\x,\l,\d]} \gamma_\lambda \ \le \ \frac{1}{2^m} \sum_{\substack{1\le i_{\sigma_{\d}^{-1}(\xi_j^\star)} \le N_{\ell_{\xi_j^\star}}\\ 2 \le j \le m}} \sum_{\substack{1\le i_{\sigma_{\d}^{-1}(\xi_j)} \le N_{\ell_{\xi_j}} \\ 1 \le j\le m}} \prod_{j=1}^m a_{\ell_{\xi_j}}(i_{\sigma_{\d}^{-1}(\xi_j)})^2.$$
Once we know that \eqref{claim-patho-even} holds, we may complete the argument as follows

\vskip-10pt

\begin{eqnarray*}
\alpha_0^{s-m} \sum_{\lambda \in \Lambda_{2m}[\k]}^{\null} \gamma_\lambda \nu_\lambda(r) & = & \alpha_0^{s-m} \sum_{\x \in C_{2m}} \sum_{\substack{\l \in \mathrm{Adm}_{2m} \setminus B_{2m} \\ \ell_{\xi_j} = k_j}} \Big( \sum_{\d \sim \l} 
\sum_{\lambda \in \Lambda_{2m}[\k,\x,\l,\d]} \gamma_\lambda \Big) r^{|\ell|} \\ & \le & \frac{1}{2^m} \, \sum_{\x \in C_{2m}} \sum_{\substack{\l \in \mathrm{Adm}_{2m} \setminus B_{2m} \\ \ell_{\xi_j} = k_j}} \Big[ \sum_{\d \sim \l} 
\Big( \prod_{j=2}^{m}N_{\ell_{\xi_j^\star}} \Big) \Big] \,r^{|\l|} \, \alpha_{\k}.
\end{eqnarray*}
Since $|\{\d \; : \; \d \sim \l\}| = M(\l)$, this completes the proof of i). Assertion ii) is proved similarly. Taking $\Lambda_{2m}'[\k,\x,\l,\d] \, = \, \{ \lambda = (\x,\l,\d,\ibar) \in \Lambda_{2m}'[\k] \}$, we may estimate $\sum_\lambda \gamma_\lambda'$ over this set as in \eqref{claim-patho-even} and obtain an extra factor $(2m)^{-1}$ in our previous bound, then we proceed as above. The proof is complete. \fin  

\subsection{$\Delta$-estimates} \label{Deltaest}

Define the sum$$\mathcal{G}(\G,\psi,r) = \sum_{g \in \G \setminus \{e\}} r^{\psi(g)},$$ whose radius of convergence $R(\G,\psi)$ is used to define the critical index $q(\G,\psi)$ as above. We now estimate $\mathbf{R}_{2m}(\Delta,r)$ and $\mathbf{R}_{2m}'(\Delta,r)$ which represent all but a finite number of terms in our decomposition \eqref{DecompositionRPS2}. These estimates are valid for $q > q(\G,\psi)$, which again represent all $q \in 2 \Z_+$ except for the isolated family considered in Paragraph \ref{Lambdaest}. The following is the core of our method. 

\begin{proposition}\label{method1}
Given $s > \frac12 q(\G,\psi)$ and $0 \le r < R(\G,\psi)$, we find 
\begin{itemize}
\item[i)] If $2 \le m \le s$ and $\k \in L_m(s)$, we obtain 
$$\hskip50pt \alpha_0^{s-m} \sum_{\delta \in \Delta_{2m}[\k]} \gamma_\delta \nu_\delta(r) \le \frac{(2m)!(s-m)!}{m! s!} \mathcal{G}(\G,\psi,r)^{m-1} r^{|\k|+1} M(\k) \alpha_{\k}.$$

\item[ii)] If $3 \le m \leq s$ and $\k \in L_m(s)$, we obtain
$$\hskip35pt \alpha_0^{s-m} \sum_{\delta\in \Delta'_{2m}[\k]} \gamma'_\delta \nu_\delta(r) \le \frac{(2m-1)!(s-m)!}{(m-1)! s!} \mathcal{G}(\G,\psi,r)^{m-2} r^{|\k|+1} M(\k) \alpha_{\k}.$$

\item[iii)] If $3 \le m \leq s$ and $\k \in L_{m-1}(s)$, we obtain 
$$\hskip20pt \alpha_0^{s-m} \sum_{\delta\in \Delta'_{2m}[\k]} \gamma'_\delta \nu_\delta(r) \le \frac{(2m-1)!(s-m+1)!}{m!s!} \mathcal{G}(\G,\psi,r)^{m-1} r^{|\k|+1} M(\k) \alpha_{\k}.$$
\end{itemize}
\end{proposition}

\dem Given $\k \in L_m(s)$, consider $$\underline{\tilde{k}} = (k_1, k_2, \ldots, k_m, -1, -1, \ldots, -1) \in \Z^{2m}.$$ For any $\d \sim \underline{\tilde{k}}$, fix a permutation $\sigma_{\d} \in \mathfrak{S}_{2m}$ satisfying $d_j = \tilde{k}_{\sigma_{\d}(j)}$ and set  $$\z(\d) = \big( \zeta_1(\d), \zeta_2(\d), \ldots, \zeta_m(\d) \big) \in \M_m$$ for the positions $1 \le \zeta_1(\d) < \zeta_2(\d) < \cdots < \zeta_m(\d) \le 2m$ where $\d \neq -1$. Note that $Z(\d) = \{\zeta_j(\d)\}_{1 \le j \le m} = \sigma_{\d}^{-1}(\{1,2,\ldots,m\})$. Once we have fixed the positions $\d$ and $\ibar = (i_1, \ldots, i_m)$ with $1 \le i_j \le N_{k_j}$, we introduce the $\Delta$-sets 
\begin{eqnarray*}
\Delta_{2m}[\k,\d] \!\!\!\! & = & \!\!\!\! \Big\{ \delta = (g_1, \ldots, g_{2m}, \stackrel{\circ}{\z(\d)}) \in \Delta_{2m}[\k] \; : \; \psi(g_j) = d_j \mbox{ for } j \in Z(\d) \Big\}, \\  \Delta_{2m}[\k,\d,\ibar] \!\!\!\! & = & \!\!\!\! \Big\{ \delta = (g_1, \ldots, g_{2m}, \stackrel{\circ}{\z(\d)})\in \Delta_{2m}[\k,\d] \; : \; g_j = \omega_{k_{\sigma_{\d}(j)}}(i_{\sigma_{\d}(j)}), \, j \in Z(\d) \Big\}.
\end{eqnarray*}
Note that $d_j = \tilde{k}_{\sigma_{\d}(j)} = k_{\sigma_{\d}(j)}$ for $j \in Z(\d)$. Then $$\Delta_{2m}[\k] \ = \ \bigcup_{\d \sim  \tilde{\k}} \ \bigcup_{\begin{subarray}{c} \ibar = (i_1, \ldots, i_m) \\ 1 \le i_j \le N_{k_j} \end{subarray}} \Delta_{2m}[\k,\d,\ibar],$$ where the union is not necessarily disjoint. We claim that
\begin{eqnarray} 
\label{gamma_delta} \null \hskip10pt \gamma_\delta \, = \, \prod_{j=1}^m a_{k_j}(i_j)^2 \ \, \mbox{for} \, \ \delta \in \Delta_{2m}[\k,\d,\ibar], \\ 
\label{sum-nu_delta} \sum_{\delta \in \Delta_{2m}[\k,\d,\ibar]} \nu_\delta(r) \ \le \ \mathcal{G}(\G,\psi,r)^{m-1} r^{|\k|+1},
\end{eqnarray}
for $\k,\d,\ibar$ fixed. Assuming our claim, we immediately obtain
\begin{eqnarray*}
\alpha_0^{s-m} \sum_{\delta \in \Delta_{2m}[\k]} \gamma_\delta \nu_\delta(r) & \le & \alpha_0^{s-m} \sum_{\d \sim  \tilde{\k}} \, \sum_{\substack{1 \le i_j \le N_{k_j} \\ 1 \le j \le m}} \sum_{\delta \in \Delta_{2m}[\k,\d,\ibar]} \gamma_\delta \nu_\delta(r) \\ & \le & \alpha_0^{s-m} \sum_{\d \sim \tilde{\k}} \, \sum_{\substack{1 \le i_j \le N_{k_j} \\ 1 \le j \le m}} 
\prod_{j=1}^m a_{k_j}(i_j)^2 \mathcal{G}(\G,\psi,r)^{m-1} r^{|\k|+1}.
\end{eqnarray*}
This gives $$\alpha_0^{s-m} \sum_{\delta \in \Delta_{2m}[\k]} \gamma_\delta \nu_\delta(r) \, \le \, M(\tilde{\k}) \mathcal{G}(\G,\psi,r)^{m-1} r^{|\k|+1} \alpha_{\k},$$ and the assertion will follow from $$M(\tilde{\k}) = \binom{2m}{m} M(k_1,k_2,\ldots,k_m) = \Big[ \binom{2m}{m} \Big/ \binom{s}{m} \Big] M(\k).$$ Claim \eqref{gamma_delta} for $\delta = (g_1,g_2,\ldots,g_{2m},\stackrel{\circ}{\z(\d)}) \in \Delta_{2m}[\k,\d,\ibar]$ follows from $$\gamma_\delta \ = \ \prod_{j=1}^m \widehat{f}(g_{\zeta_j(\d)})^2 \, = \prod_{j \in Z(\d)} \widehat{f}(g_j)^2 \, = \prod_{j \in Z(\d)} a_{k_{\sigma_{\d}(j)}}(i_{\sigma_{\d}(j)})^2 \, = \, \prod_{j=1}^m a_{k_j}(i_j)^2.$$ We now turn to \eqref{sum-nu_delta}. Write $\{\eta_1(\d),\eta_2(\d),\ldots,\eta_m(\d)\} = \{1,2,\ldots,2m\} \setminus Z(\d)$, with $1\le \eta_1(\d) < \eta_2(\d) < \cdots < \eta_m(\d) \le 2m$. Now, given $h_1,h_2,\ldots, h_{m-1} \in \G \setminus \{e\}$ we define the sets
$$\Delta_{2m}[\k,\d,\ibar,h_1,\ldots,h_{m-1}] \, = \, \Big\{ \delta \in \Delta_{2m}[\k,\d,\ibar] \; : \; g_{\eta_j(\d)} = h_j \mbox{ for } 1 \le j \le m-1 \Big\}.$$
Then, it is easy to prove that
\begin{enumerate}
\item[a)] $\nu_\delta(r) \le r^{\sum_{j=1}^{m-1}\psi(h_j) + |\k| + 1}$,
\item[b)] $\big| \Delta_{2m}[\k,\d,\ibar,h_1,\ldots,h_{m-1}] \big| \le 1$. 
\end{enumerate}
Indeed, a) follows from $\psi(g_{\zeta_j(\d)}) = k_j$, $\psi(g_{\eta_j(\d)}) = \psi(h_j)$ and $\psi(g_{\eta_m(\d)}) \ge 1$. On the other hand, b) follows since the only entry in $\delta \in \Delta_{2m}[\k,\d,\ibar,h_1,\ldots,h_{m-1}]$ which is not determined a priori is $g_{\eta_m(\d)}$. Hence, the restriction $g_1 g_2 \cdots g_{2m} = e$ yields the second assertion. Since we have a disjoint union $$\Delta_{2m}[\k,\d,\ibar] \, = \, \bigsqcup_{\substack{h_j \neq e \\ 1 \le j \le m-1}} \Delta_{2m}[\k,\d,\ibar,h_1,\cdots,h_{m-1}],$$ 
we obtain
\begin{eqnarray*}
\sum_{\delta \in \Delta_{2m}[\k,\d,\ibar]} \nu_\delta(r) & = & \sum_{\substack{h_j \neq e \\ 1 \le j \le m-1}} \sum_{\delta \in \Delta_{2m}[\k,\d,\ibar,h_1,\ldots,h_{m-1}]} \nu_\delta(r) \\ & \le & \sum_{\substack{h_j \neq e \\ 1 \le j \le m-1}} r^{\sum_{j=1}^{m-1} \psi(h_j) + |\k| + 1} \ = \ \mathcal{G}(\G,\psi,r)^{m-1} r^{|\k|+1}.
\end{eqnarray*}
This proves \eqref{sum-nu_delta} and concludes the proof of i). The proof of ii) and iii) is very similar. Given $\k \in L_{m-1}(s) \cup L_m(s)$, we may construct $\tilde{\k}$ as above and define $\Delta'_{2m}[\k,\d,\ibar] \subset \Delta'_{2m}[\k]$ for $\d \sim \tilde{\k}$ and $\ibar \in \prod_{j=1}^m \{1,\ldots,N_{k_j}\}$ accordingly, so that $$\Delta'_{2m}[\k] \, = \, \bigcup_{\d \sim \tilde{\k}} \ \bigcup_{\begin{subarray}{c} \ibar = (i_1, \ldots, i_m) \\ 1 \le i_j \le N_{k_j} \end{subarray}} \Delta'_{2m}[\k,\d,\ibar].$$
In that case, we claim that for $\k,\d,\ibar$ fixed we have 
\begin{eqnarray}
\label{gamma'_delta} \hskip25pt \gamma'_\delta \, = \, \frac{1}{2m} \prod_{j=1}^m a_{k_j}(i_j)^2 & \mbox{ for } & \delta \in \Delta_{2m}'[\k,\d,\ibar], \\ \label{sum'-nu_delta} \hskip25pt  \sum_{\delta \in \Delta'_{2m}[\k,\d,\ibar]} \nu_\delta(r) \ \le \ 
\mathcal{G}(\G,\psi,r)^{m-1} r^{|\k|+1} & \mbox{ for } & \k \in L_{m-1}(s), \\ \label{sum'-nu_delta2} \hskip25pt \sum_{\delta \in \Delta'_{2m}[\k,\d,\ibar]} \nu_\delta(r) \ \le \ m \mathcal{G}(\G,\psi,r)^{m-2} r^{|\k|+1} & \mbox{ for } & \k \in L_m(s).
\end{eqnarray}
The proof of \eqref{gamma'_delta} and \eqref{sum'-nu_delta} follows verbatim \eqref{gamma_delta} and \eqref{sum-nu_delta}, while the argument for \eqref{sum'-nu_delta2} requires slight modifications. Given $\delta \in \Delta'_{2m}[\k,\d,\ibar]$, in this situation we have $g_{\zeta_j(\d)} \neq e$ for all $1 \le j \le m$ and $g_{\eta_j(\d)} = e$ holds for one and only one $1 \le j \le m$. 
This leads to define $$\Delta'_{2m}[\k,\d,\ibar,j_0] \, = \, \Big\{ \delta \in \Delta'_{2m}[\k,\d,\ibar] \; : \;  g_{\eta_{j_0}(\d)} = e \Big\}$$ for $1 \le j_0 \le m$. Then we set for $h_1, h_2, \ldots,h_{m-2} \in \G \setminus \{e\}$ 
$$\Delta'_{2m}[\k,\d,\ibar,j_0,h_1,\ldots,h_{m-2}] = \Big\{ \delta \in \Delta'_{2m}[\k,\d,\ibar,j_0] \; : \; g_{\tilde{\eta}_{j}(\d,j_0)} = h_j \mbox{ for } j \le m-2 \Big\},$$
where we write $$\tilde{\eta}_{j}(\d,j_0) \, = \, \left\{\begin{array}{ll} \eta_j(\d) & \quad \mbox{ for } \ 1 \le j \le j_0-1, \\ \eta_{j+1}(\d) & \quad \mbox{ for } j_0 \le j \le m-1. \end{array} \right.$$ The word $g_{\tilde{\eta}_{m-1}(\d,j_0)}$ is uniquely determined in $\G$ by $g_1 g_2 \cdots g_{2m} = e$ and the fact that the other $g_j$'s are determined a priori in $\Delta'_{2m}[\k,\d,\ibar,j_0,h_1,\ldots,h_{m-2}]$. Hence we also deduce that $$\big| \Delta'_{2m}[\k,\d,\ibar,j_0,h_1,\ldots,h_{m-2}] \big| \, \le \,  1.$$ Since $\nu_\delta(r) \le r^{\sum_{j=1}^{m-2} \psi(h_j) + |\k| + 1}$ for $\delta \in \Delta'_{2m}[\k,\d,\ibar,j_0,h_1,\ldots,h_{m-2}]$ and $$\Delta'_{2m}[\k,\d,\ibar] \, = \, \bigsqcup_{1 \le j_0 \le m} \bigsqcup_{\substack{h_j \neq e \\ 1 \le j \le m-2}} \Delta'_{2m}[\k,\d,\ibar,j_0,h_1,\cdots,h_{m-2}],$$ 
we obtain
\begin{eqnarray*}
\sum_{\delta \in \Delta'_{2m}[\k,\d,\ibar]} \nu_\delta(r) \!\!\!\! & = & \!\!\!\! \sum_{j_0 = 1}^m \sum_{\substack{h_j \neq e \\ 1 \le j \le m-2}}\sum_{\delta \in \Delta'_{2m}[\k,\d,\ibar,j_0,h_1,\ldots,h_{m-2}]} \nu_\delta(r) \\ \!\!\!\! & \le & \!\!\!\! \sum_{j_0 = 1}^m \sum_{\substack{h_j \neq e \\ 1 \le j \le m-2}} r^{\sum_{j=1}^{m-2} \psi(h_j) + |\k| + 1} \, = \, m \, \mathcal{G}(\G,\psi,r)^{m-2} r^{|\k|+1}.
\end{eqnarray*}
This proves \eqref{sum'-nu_delta2}. The assertion follows by joining the pieces as we did for i). \fin

\begin{remark} \label{ConstMu}
\emph{If $|\k|$ is large enough and $q > q(\G,\psi)$, Proposition \ref{method1} provides the estimates for \eqref{AIM}. The optimal size $|\k|$ when $\k \in L_m(s)$ is what determines the critical function $\mu_q(\G,\psi,m)$, see the Strategy below for further details. In view of this, it is clear that the critical function is very much affected by the growth of $(\G,\psi)$ through the size of $\mathcal{G}(\G,\psi,r)$.}
\end{remark}

\subsection{Strategy} \label{Strategy}

\begin{itemize}
\item[{\bf i)}] \textbf{Admissible lengths.} Given any pair $(\G,\psi)$ satisfying our growth and cancellation conditions, we will begin by refining the set $\mathrm{Adm}_u$ of admissible lengths, which depends as we will see on the concrete group $\G$. 

\vskip5pt 

\item[{\bf ii)}] \textbf{Estimates for $s_2(r)$ and $s_3(r)$.} When $u=2$ we have 
$$s_2(r) \, = \, \sum_{g \neq e}^{\null} \widehat{f}(g)^2 r^{2 \psi(g)} \, = \, \sum_{k \ge 1}^{\null} r^{2k} \alpha_k.$$ When $u = 3$, the estimate $$\widehat{f}(e)^{q-3} s_3(r) \ \le \ \displaystyle \sum_{\k \in L(s) \setminus \{\underline{0}\}}^{\null} C_{q,s_3}^{\mathrm{left}}[\k](r) \, \alpha_{\k}$$ in the line of \eqref{AIM} requires a different way of completing squares. The reason is that $C_q^{\mathrm{left}}[\k](r) = \binom{q}{2} C_{q,s_2}^{\mathrm{left}}[\k](r) = C_q^{\mathrm{right}}[\k] = q/2$ is attained at optimal time $r = 1/\sqrt{q-1}$ for the singular term $\k = (1,\underline{0})$. This means that $\alpha_0^{s-1}\alpha_1$ can only appear in the estimate of $s_{2}(r)$, which forces us to be specially careful in the case $u=3$ to avoid this term. On the other hand our assumption on the absence of small loops is crucial in our estimate of the sum $s_3(r)$, which varies from one group to another. All our estimates for $C_q^{\mathrm{left}}[\k](r)$ arising from $s_2$ and $s_3$ sums are collected in what we call $(\alpha)$-estimates in the proofs of Theorems A1-A3.

\vskip5pt

\item[{\bf iii)}] \textbf{General goal}. According to \eqref{DecompositionRPS1} and \eqref{DecompositionRPS2}, our aim is to show that $$\qquad \widehat{f}(e)^{q-{2m}} \Big[ \mathbf{R}_{2m}(r) + \mathbf{P}_{2m}(r) + \mathbf{S}_{2m}(r) \Big] \le \sum_{\k \in L_m(s) \setminus \{\underline{0}\}}^{\null} C_{q,s_{2m}}^{\mathrm{left}}[\k](r) \, \alpha_{\k},$$ and similarly in the odd case. Given $1 \le m \le s$ and $\k \in L_m(s)$, we have 
$$\hskip30pt C_q^{\mathrm{left}}[\k](r) \, = \, \sum_{u=2 \vee (2m-1)}^{q \wedge (2m+1)} \binom{q}{u} C_{q,s_{u}}^{\mathrm{left}}[\k](r),$$ which refines our original definition in \eqref{AIM}. Indeed, taking into account how we complete squares to obtain the regular and pathological terms and also the way we will complete squares for the super-pathological terms, it will become clear that only the sums $s_{2m-1}(r)$, $s_{2m}(r)$ and $s_{2m+1}(r)$ may have a contribution to $\alpha_{\k} = \alpha_{k_1} \alpha_{k_2} \cdots \alpha_{k_m} \alpha_0^{s-m}$ when $\k \in L_m(s)$. Given $q \in 2\Z_+$ and $r=1/\sqrt{q-1}$, the aim is to prove  $$C_q^{\mathrm{left}}[\k](r) \le C_q^{\mathrm{right}}[\k] \quad \mbox{for every} \quad \k \in L(s)\setminus \{0\}.$$

\vskip5pt

\item[{\bf iv)}] \textbf{Algorithm.} Assume first that $q>q(\G,\psi)$. Then, we use Proposition \ref{method1} to construct the coefficients $C_q^{\mathrm{left}}[\k](r)$ as explained in \textbf{iii)}. These will be the $(\gamma)$-estimates in our proof of Theorems A1-A3. Define $\mathcal{A}$ to be the set of $\k \in L(s)\setminus \{0\}$ verifying $C_q^{\mathrm{left}}[\k](r) \le C_q^{\mathrm{right}}[\k]$ with this choice of left coefficients. According to Proposition \ref{method1}, the validity of this inequality just depends on the length $|\k|$. We then construct the critical function as follows $$\mu_q(\G,\psi,m) \, = \, \inf \Big\{ |\k|-m \,: \, \k \in \mathcal{A} \cap L_m(s) \Big\}$$ for $1\le m\le \frac{q}{2}$, so that the set $\mathrm{Reg}_m(s)$ coincides with  $\mathcal{A}\cap L_m(s)$. This gives $C_q^{\mathrm{left}}[\k](r) \le C_q^{\mathrm{right}}[\k]$ for all $\k \in \mathrm{Reg}_m(s)$ and $1 \le m \le s$. We write $\mathrm{Pat}_m(s)$ for $L_m(s) \setminus \mathrm{Reg}_m(s)$. Given  $\k \in \mathrm{Pat}_m(s)$, we now use (a slightly modified version of) Proposition \ref{Lambda-est}, see the $(\delta)$-estimates in the proofs of Theorems A1-A3. This allows us to construct a refined estimate of $C_q^{\mathrm{left}}[\k](r)$ for pathological terms $\k \in \mathrm{Pat}_m(s)$. Then, we use computer assistance to show that most of these terms satisfy $C_q^{\mathrm{left}}[\k](r) \le C_q^{\mathrm{right}}[\k]$ with this choice of left coefficients. When this is not the case, we say that $\k$ is a super-pathological term and denote by $ \mbox{S-Pat}_m(s)$ the set of all these $\k$'s. In order to fix these terms, we must refine the estimates for the sums $s_u[\l](r)$ which may contribute to them. This may require again the computer to find exact expressions for $s_u[\l](r)$ and then complete squares matching the expected inequalities. This happens particularly when the admissible length $\l$ belongs to our exceptional sets $B_u$. At the end, we may construct finer coefficients $C_q^{\mathrm{left}}[\k](r)$ for $\k \in  \mbox{S-Pat}_m(s)$ satisfying the expected estimate $$C_q^{\mathrm{left}}[\k](r) \le C_q^{\mathrm{right}}[\k],$$ which completes the argument for $q>q(\G,\psi)$. If $q \le q(\G,\psi)$, we reapply the same method. However, since Proposition \ref{method1} is no longer valid in that case, we replace it by a bound coming from Proposition \ref{Lambda-est}, see the $(\beta)$-estimates in the proofs of Theorems A1-A3 for further details. 

\vskip5pt

\item[{\bf v)}] \textbf{Extrapolation.} General hypercontractivity for $1 < p \le q < \infty$ follows by adapting Gross technique \cite{G1,G2} to our setting. By interpolation we may prove $L_2 \to L_{2 + \varepsilon}$ hypercontractive inequalities loosing a factor $\log 3$ with respect to the expected optimal time. Differentiating the dual inequality at time $0$ yields the corresponding logarithmic Sobolev inequality. Combining this with a Dirichlet form inequality for the infinitesimal generator  gives rise to the general result. More details can be found in Appendix A below.   
\end{itemize}

\begin{remark}
\emph{Let us briefly review our conditions on $\psi$. The spectral gap is just necessary for hypercontractivity, as shown in the Introduction. The subadditivity produces a Poisson-like semigroup and it is also crucial to determine the admissible lengths. The conditional negativity ---Markovianity of the semigroup--- has been used to reduce the problem to the positive cone and it is also crucial to apply Gross extrapolation argument in Appendix A. The growth condition is needed to apply our second way to complete squares, while the cancellation condition ---absence of small loops--- is needed for our careful estimate of the sum $s_3(r)$ and to deal with super-pathological terms.}
\end{remark}

\section{{\bf Optimal time estimates}} \label{Optimalt}

In this section we will use our combinatorial method to prove Theorems A1-A3 up to some numerical analysis and technical inequalities, which have been postponed to Appendices C and D for clarity in the exposition.

\subsection{Free groups} \label{Sect-ThA}

Let us write $c_1, c_2, \ldots, c_n$ to denote the generators of the free group $\F_n$ equipped with the associated word length $|\hskip0.5pt \cdot|$. In that situation we have $$ N_k = 2n(2n-1)^{k-1}, \quad \mathcal{G}(\F_n,|\cdot|,r) = \frac{2nr}{1-(2n-1)r}
\quad \mbox{and} \quad q(\F_n,|\cdot|)=4n^2-4n+2.$$
Recall that we enumerated the set $W_k$ by $w_k(1), w_k(2), \ldots, w_k(N_k)$. Additionally, we order them by imposing that $w_k(i) = w_k(\frac{N_k}{2}+i)^{-1}$ for all $1 \le i \le \frac{N_k}{2}$ and $w_1(i) = c_i$ for $1 \le i \le n$, which will be helpful for some of our estimates below. Note also that $$\widehat{f}(g) = \widehat{f}(g^{-1}) \ge 0 \ \Rightarrow \ a_k(i) = a_k(\mbox{$\frac{N_k}{2}$}+i) \quad \mbox{so that} \quad \alpha_k = 2 \sum_{i=1}^{N_k/2} a_k(i)^2.$$ 

\subsubsection{Admissible lengths} 

Let us refine the set of admissible lengths for free groups.

\begin{lemma} \label{admlength-Fn}
Given $g_1, g_2, \ldots, g_u \in \F_n$, we find
\begin{itemize}
\item[i)] There exists $0 \le m \le \frac12 \sum_{j=1}^u |g_j|$ such that $$|g_1 g_2 \cdots g_u| = \sum_{j=1}^u |g_j| - 2m.$$

\item[ii)] If $g_1 \cdots g_u = e$, then $\sum_j |g_j|$ is even and $\sum_{j \neq i} |g_j| \ge |g_i|$ for any $1 \le i \le u$.
\end{itemize}
\end{lemma}

\dem The second assertion follows easily from the first one, which in turn can be proved by induction. Indeed, it is clear for $u=1$, while for $u=2$ the identity holds for $0 \le m \le |g_1| \wedge |g_2| \le \frac{1}{2}(|g_1|+|g_2|)$. Namely, in case of cancellation in $g_1g_2$ there is necessarily an even number of letters which will disappear, but no more than $2|g_1| \wedge 2|g_2|$. The general case easily follows from this one by induction. \fin

This lemma refines even more the set of admissible lengths $|g_j|$ which may appear in the sums $s_u(r)$ for free groups. Namely, along the proof of Theorem A1 we will redefine the set $\mathrm{Adm}_u$ for $2 \le u \le q$ by adding a parity condition $$\mathrm{Adm}_u \, = \, \Big\{ \l \in L(u) \; : \; \ell_1 \ge \ell_2 \ge \ldots \ge \ell_u \ge 1, \; \sum_{j=1}^{u} \ell_j \ \mbox{even}, \; \ell_1 \le \sum_{j=2}^{u} \ell_j \Big\}.$$ With this new definition, we still have the identity $s_u(r) \, = \, \sum_{\l \in \mathrm{Adm}_u}^{\null} s_u[\l](r)$.

\subsubsection{Estimates for $s_3(r)$} \label{S2S3}

As we already pointed out in the strategy, our estimate for the sum $s_3(r)$ must be treated apart to avoid the term $\k = (1,\underline{0}) \in L(s)$, which was already \lq\lq saturated" by the sum $s_2(r)$. 
 
\begin{proposition} \label{s23}
We have
\begin{eqnarray*}
\widehat{f}(e) s_3(r) & \le & \frac{3}{4}r^4 \alpha_1^2 + 6 \sum_{1\le i<j \le n}a_1(i)^2a_1(j)^2 r^4 + \frac32 r^4 \alpha_0 \alpha_2
\\ [5pt] & + & \sum_{\begin{subarray}{c} \ell_2 \ge \ell_3 \ge 1 \\ (\ell_2, \ell_3) \neq (1,1) \end{subarray}} \sum_{m=0}^{\lfloor \ell_3/2 \rfloor} A_m(\ell_2, \ell_3,r) \alpha_{\ell_2} \alpha_{\ell_3} + B_m(\ell_2,\ell_3,r) \alpha_0 \alpha_{\ell_2+\ell_3-2m} \\ 
[5pt] & \le & \sum_{ \ell_2 \ge \ell_3 \ge 1} \sum_{m=0}^{\lfloor \ell_3/2 \rfloor} \big[A_m(\ell_2, \ell_3,r) \alpha_{\ell_2} \alpha_{\ell_3}  + B_m(\ell_2,\ell_3,r) \alpha_0 \alpha_{\ell_2+\ell_3-2m} \big],
\end{eqnarray*}
where $A_0(\ell_2, \ell_3,r) = B_0(\ell_2,\ell_3,r)$, 
$$\begin{array}{lcl} \displaystyle A_m(\ell_2,\ell_3,r) = \frac12 M(\ell_2+\ell_3-2m,\ell_2,\ell_3) r^{2(\ell_2+\ell_3-m)} & m \ge 0, \\ [8pt] \displaystyle B_m(\ell_2, \ell_3,r) = \frac12 (2n-2)(2n-1)^{m-1} M(\ell_2+\ell_3-2m,\ell_2,\ell_3) r^{2(\ell_2+\ell_3-m)} & m \ge 1.
\end{array}$$
\end{proposition}

\dem 
We start by decomposing 
$$s_3(r) =  \sum_{\l \in \mathrm{Adm}_3 \setminus \{(2,1,1)\}} s_3[\l](r) \, + \, s_3[(2,1,1)](r).$$
Given $\l \in \mathrm{Adm}_3$ and $\d \sim \l$, we set $$\Lambda_3(\d) \, = \, \Big\{ (g_1, g_2, g_3) \in \F_n^3 \, : \ g_1 g_2 g_3 = e \ \mbox{and} \ |g_j| = d_j \Big\}.$$ If $\d \sim \l$ is fixed, pick a permutation $\sigma_{\d} \in \mathfrak{S}_3$ such that $d_j = \ell_{\sigma_{\d}(j)}$ and set $$\Lambda_3(\d,h) \, = \, \Big\{ (g_1, g_2, g_3) \in \Lambda_3(\d) \, : \ g_{\sigma_{\d}^{-1}(1)} = h \Big\}$$ for any $h \in W_{\ell_1}$. Observe that $\l = (\ell_1, \ell_2, \ell_3) \in \mathrm{Adm}_3$ if and only if $\ell_2 \ge \ell_3 \ge 1$ and $\ell_1 = \ell_2 + \ell_3 - 2m$ for some integer $0 \le m \le \ell_3/2$. Our estimate for $s_3(r)$ relies on the following properties:
\begin{itemize}
\item[a)] $\Lambda_3(\d,h_1) \cap \Lambda_3(\d,h_2) = \emptyset$ for $h_1 \neq h_2$,

\vskip10pt

\item[b)] $\displaystyle \sum_{(g_1, g_2, g_3) \in \Lambda_3(\d)} \widehat{f}(g_{\sigma_{\d}^{-1}(2)})^2 \widehat{f}(g_{\sigma_{\d}^{-1}(3)})^2 \le \alpha_{\ell_2} \alpha_{\ell_3}$,

\item[c)] If $|h| = \ell_2 + \ell_3 - 2m$, $|\Lambda_3(\d,h)| = \delta_{m=0} + (2n-2)(2n-1)^{m-1} \delta_{m>0}$.
\end{itemize}
The first one is clear and the second one follows from the fact that the sum has only two degrees of freedom due to the constraint $g_1g_2g_3=e$. To justify c), we assume for simplicity that $\d=\l$ and $\sigma_{\d}$ is the identity map, the other cases being similar. Let $h \in W_{\ell_1}$ and write $h^{-1} = b_1 b_2 \cdots b_{\ell_1}$ with $|b_j|=1$ and $b_j \neq b_{j+1}^{-1}$. 
Then $(g_1,g_2,g_3) \in \Lambda_3(\d,h)$ if and only if $g_1 = h$ and $$g_2g_3 = h^{-1} = \underbrace{b_1 \cdots b_{\ell_2-m} \gamma}_{g_2} \underbrace{\gamma^{-1} b_{\ell_2-m+1} \cdots b_{\ell_1}}_{g_3}$$ for some word $\gamma$ satisfying $|\gamma| = m$ and $|b_{\ell_2-m}\gamma| = |\gamma^{-1}b_{\ell_2-m+1}|=m+1$. Since $b_{\ell_2-m} b_{\ell_2-m+1} \neq e$, this means that we have $(2n-2)(2n-1)^{m-1}$ possible words $\gamma$ for $m > 0$ and $1$ possible $\gamma$ for $m=0$. The assertion for $\l \neq (2,1,1)$ now follows using these properties and completing squares. Namely, we have 
\begin{eqnarray*}
\lefteqn{\hskip-10pt \widehat{f}(e) \sum_{\l \in \mathrm{Adm}_3 \setminus \{(2,1,1)\}} s_3[\l](r)} \\ & \le & \sum_{\begin{subarray}{c} \l \in \mathrm{Adm}_3 \setminus \{(2,1,1)\} \\ \d \sim \l \end{subarray}} \sum_{\begin{subarray}{c} h \in W_{\ell_1} \\ (g_1,g_2,g_3) \in \Lambda_3(\d,h) \end{subarray}} \widehat{f}(e) \prod_{j=1}^3 \widehat{f}(g_j) \, r^{|\l|} \\ & = & \sum_{\begin{subarray}{c} \ell_2 \ge \ell_3 \ge 1 \\ (\ell_2, \ell_3) \neq (1,1) \end{subarray}} \sum_{m=0}^{\lfloor \ell_3/2 \rfloor} \hskip8pt \sum_{\d \sim \l_m} \hskip5pt \Big( \!\!\!\! \sum_{\begin{subarray}{c} h \in W_{\ell_2 + \ell_3 - 2m} \\ (g_1,g_2,g_3) \in \Lambda_3(\d,h) \end{subarray}} \!\!\!\widehat{f}(e) \prod_{j=1}^3 \widehat{f}(g_j) \, \Big) \ r^{2(\ell_2 + \ell_3 - m)},
\end{eqnarray*}
with $\l_m = (\ell_2 + \ell_3 - 2m, \ell_2, \ell_3)$. Now, completing squares in the last sum we obtain 
\begin{eqnarray*}
\lefteqn{\hskip-10pt \sum_{\begin{subarray}{c} h \in W_{\ell_2 + \ell_3 -2m} \\ (g_1,g_2,g_3) \in \Lambda_3(\d,h) \end{subarray}} \!\!\!\!\!\widehat{f}(e) \prod_{j=1}^3 \widehat{f}(g_j)} \\ & \le & \sum_{\begin{subarray}{c} h \in W_{\ell_2 + \ell_3 - 2m} \\ (g_1,g_2,g_3) \in \Lambda_3(\d,h) \end{subarray}} \!\!\!\!\!\frac{\widehat{f}(e)^2 \widehat{f}(h)^2 + \widehat{f}(g_{\sigma_{\d}^{-1}(2)})^2 \widehat{f}(g_{\sigma_{\d}^{-1}(3)})^2}{2} \\ & \le & \hskip15pt \frac12 \Big( \underbrace{\delta_{m=0} + (2n-2)(2n-1)^{m-1} \delta_{m > 0}}_{\beta_m} \Big) \alpha_0 \alpha_{\ell_2 + \ell_3 -2m} \ + \ \frac12 \alpha_{\ell_2} \alpha_{\ell_3}.
\end{eqnarray*}
Combining both estimates and since there are $M(\l)$ choices for $\d \sim \l$, we find
\begin{eqnarray*}
\lefteqn{\hskip-15pt \widehat{f}(e) \sum_{\l \in \mathrm{Adm}_3 \setminus \{(2,1,1)\}} s_3[\l](r)} \\ & \le &  
\sum_{\begin{subarray}{c} \ell_2 \ge \ell_3 \ge 1 \\ (\ell_2, \ell_3) \neq (1,1) \end{subarray}} \sum_{m=0}^{\lfloor \ell_3/2 \rfloor} \frac12 M(\l_m) \, r^{2(\ell_2 + \ell_3 - m)} \, \alpha_{\ell_2} \alpha_{\ell_3} \\ 
& + &  \sum_{\begin{subarray}{c} \ell_2 \ge \ell_3 \ge 1 \\ (\ell_2, \ell_3) \neq (1,1) \end{subarray}} \sum_{m=0}^{\lfloor \ell_3/2 \rfloor} \frac12 M(\l_m) \, r^{2(\ell_2 + \ell_3 - m)} \beta_m \alpha_0 \alpha_{\ell_2 + \ell_3 -2m}. 
\end{eqnarray*}
It remains to prove the following inequality  
$$\widehat{f}(e) s_3[(2,1,1)](r) \le \frac{3}{4}r^4 (\alpha_1^2+2\alpha_0\alpha_2) +6\sum_{1\le i<j \le n}a_1(i)^2a_1(j)^2 r^4  \le \frac32 r^4 (\alpha_1^2 + \alpha_0 \alpha_2)$$
to conclude the proof. Assuming that $w_2(i)=c_i^2$, we have 
\begin{eqnarray*}
\widehat{f}(e) s_3[(2,1,1)](r) \!\!\!\! & = & \!\!\!\! \sum_{\begin{subarray}{c} g_1g_2g_3=e \\ (|g_1|,|g_2|,|g_3|)\sim (2,1,1) \end{subarray}} \widehat{f}(e) \widehat{f}(g_1)\widehat{f}(g_2)\widehat{f}(g_3) r^4 \\ \!\!\!\! & = & \!\!\!\! 6a_0\Big[\sum_{i=1}^n a_1(i)^2a_2(i) + \sum_{1\le i<j\le n} a_1(i)a_1(j)\Big(\sum_{k\in K(i,j)}a_2(k)\Big)\Big]r^4,
\end{eqnarray*}
where $K(i,j)$ denotes the subset of $\{1,2,\ldots,N_2/2\}$ of cardinal $4$ verifying that $\{w_2(k) \, : \, k \in K(i,j)\}$ is constituted of words with the first letter in $\{c_i,c_i^{-1}\}$ and the last one in $\{c_j,c_j^{-1}\}$. Completing squares as before we get
\begin{eqnarray*}
\null \hskip20pt \lefteqn{\hskip-15pt \widehat{f}(e) s_3[(2,1,1)](r)} \\ &\le &  
3\Big[\sum_{i=1}^n a_1(i)^4 + \sum_{i=1}^{N_2/2}a_0^2a_2(i)^2 + 4 \sum_{1\le i<j\le n} a_1(i)^2a_1(j)^2\Big]r^4 \\
&\le &  \frac{3}{4}r^4(\alpha_1^2+2\alpha_0 \alpha_2)+6\sum_{1\le i<j \le n}a_1(i)^2a_1(j)^2 r^4 \le  \frac{3}{2}r^4 (\alpha_1^2 + \alpha_0 \alpha_2). \hskip20pt \square
\end{eqnarray*}

\subsubsection{Numerical estimates for $\F_n$} \label{numericalF2}

Our goal in \eqref{AIM} was to show that $$\widehat{f}(e)^{q-u} s_u(r) \ \le \ \sum_{\k \in L(s) \setminus \{\underline{0}\}}^{\null} C_{q,s_u}^{\mathrm{left}}[\k](r) \, \alpha_{\k},$$ for all $2\leq u \leq q$ and certain coefficients $C_{q,s_u}^{\mathrm{left}}[\k](r)$. In this section we will identify these coefficients for $\F_n$ equipped with the word length. The proofs of these results are quite simple but tedious, so we have decided to collect them in Appendix C. This will make the core of the argument more transparent for the reader. Let us recall that $R(\F_n,|\cdot|)=\frac{1}{2n-1}$, so that we decompose the sums $s_u(r)$ for $4\leq u \leq q$ following \eqref{DecompositionRPS1} or \eqref{DecompositionRPS2} for $r \ge \frac{1}{2n-1}$ or $r < \frac{1}{2n-1}$ respectively. Hence, we need to define three kinds of left-coefficients as follows. Take $* = \Lambda$ for $r \ge \frac{1}{2n-1}$ and $*=\Delta$ otherwise. Then, the coefficients for $4 \le u \le 
 q$ are given by 
\begin{equation} \label{Cleft3}
C_{q,s_u}^{\mathrm{left}}[\k](r) \, = \, C_{q,s_u}^{\mathrm{left},\mathbf{R}}[\k](r)+C_{q,s_u}^{\mathrm{left},\mathbf{P}}[\k](r)+C_{q,s_u}^{\mathrm{left},\mathbf{S}}[\k](r),
\end{equation}
where the right hand side is determined by the inequalities  
\begin{eqnarray*}
\displaystyle \alpha_0^{s-m} \mathbf{R}_{2m}(*,r) & \le & \sum_{\k \in L_m(s)}^{\null} C_{q,s_{2m}}^{\mathrm{left},\mathbf{R}}[\k](r) \alpha_{\k}, \\ \alpha_0^{s-m} \mathbf{P}_{2m}(r) & \le & \sum_{\k \in L_m(s)}^{\null} C_{q,s_{2m}}^{\mathrm{left},\mathbf{P}}[\k](r) \alpha_{\k}, \\ \alpha_0^{s-m} \, \mathbf{S}_{2m}(r) & \le & \sum_{\k \in L_m(s)}^{\null} C_{q,s_{2m}}^{\mathrm{left},\mathbf{S}}[\k](r) \alpha_{\k},
\end{eqnarray*}
for $u=2m$ and similarly replacing $(\mathbf{R}_{2m}(*,r), \mathbf{P}_{2m}(r), \mathbf{S}_{2m}(r), L_m(s))$ by the sums $(\mathbf{R}_{2m}'(*,r), \mathbf{P}_{2m}'(r), \mathbf{S}_{2m}'(r))$ with the summation index $L_m'(s)$ when $u=2m-1$. By construction, regular and pathological terms do not contribute to each other 
\begin{equation} \label{Reg-Pat}
\begin{array}{rcl}
\displaystyle \k \in \mathrm{Reg}(s) := \bigcup_{m =2}^s \mathrm{Reg}_m(s) & \Rightarrow & \displaystyle C_{q,s_{u}}^{\mathrm{left}, \mathbf{P}}[\k](r)=0 \quad \mbox{for} \quad 4 \le u \le q, \\ \displaystyle \k \in \, \mathrm{Pat}(s) := \bigcup_{m =2}^s \, \mathrm{Pat}_m(s) & \Rightarrow & \displaystyle C_{q,s_{u}}^{\mathrm{left},\mathbf{R}}[\k](r)=0 \quad \mbox{for} \quad 4 \le u \le q. 
\end{array}
\end{equation}
Moreover, we will see that $\mathbf{S}_{2m}(r)$ and $\mathbf{S}'_{2m}(r)$ do not affect the regular terms $\k \in \mathrm{Reg}(s) \Rightarrow C_{q,s_{u}}^{\mathrm{left},\mathbf{S}}[\k](r)=0$. Hence, $C_{q,s_u}^{\mathrm{left}}[\k](r)= C_{q,s_{u}}^{\mathrm{left},\mathbf{R}}[\k](r)$ for regular $\k$'s. It is crucial for the numerical part of the method that the critical function is uniformly bounded. This holds true for $\F_2$ equipped with the length function. Its critical function is given by the graph below.  

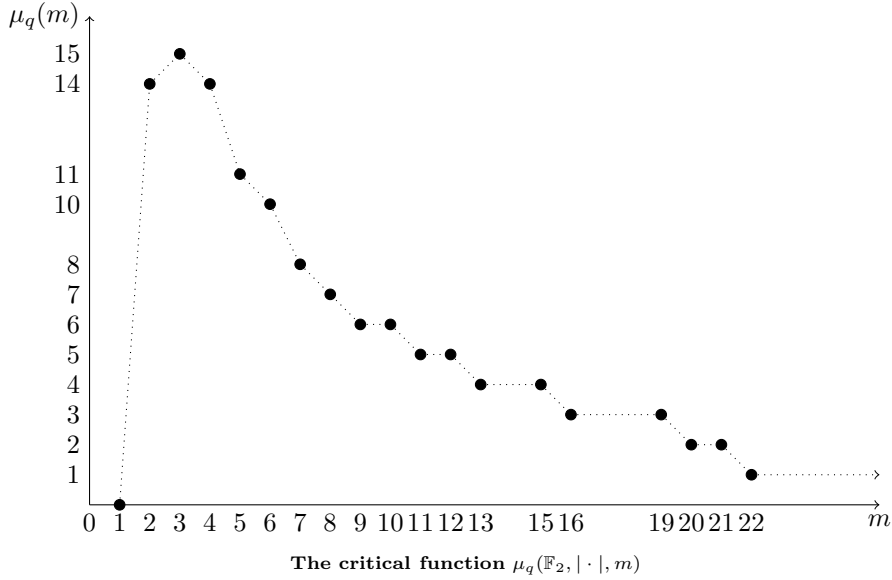
\begin{tikzpicture}

\draw[->] (-0.5,0) -- (10,0) node[anchor=north] {$m$};
\draw	(-0.5,0) node[anchor=north] {0}
		(-0.1,0) node[anchor=north] {1}
		(0.3,0) node[anchor=north] {2}
		(0.7,0) node[anchor=north] {3}
		(1.1,0) node[anchor=north] {4}
		(1.5,0) node[anchor=north] {5}
		(1.9,0) node[anchor=north] {6}
		(2.3,0) node[anchor=north] {7}
		(2.7,0) node[anchor=north] {8}
		(3.1,0) node[anchor=north] {9}
		(3.5,0) node[anchor=north] {10}
		(3.9,0) node[anchor=north] {11}		
		(4.3,0) node[anchor=north] {12}
		(4.7,0) node[anchor=north] {13}
		(5.5,0) node[anchor=north] {15}
		(5.9,0) node[anchor=north] {16}
		(7.1,0) node[anchor=north] {19}
		(7.5,0) node[anchor=north] {20}
		(7.9,0) node[anchor=north] {21}
		(8.3,0) node[anchor=north] {22};
\draw	(4.5,-0.8) node{{\scriptsize {\bf The critical function} $\mu_q(\F_2, |\cdot|, m)$}};

\draw[->] (-0.5,0) -- (-0.5,6.5) node[anchor=east] {$\mu_q(m)$};

\draw 	(-0.5,0.4) node[anchor=east] {1}
		(-0.5,0.8) node[anchor=east] {2}
		(-0.5,1.2) node[anchor=east] {3}
		(-0.5,1.6) node[anchor=east] {4}
		(-0.5,2.0) node[anchor=east] {5}
		(-0.5,2.4) node[anchor=east] {6}
		(-0.5,2.8) node[anchor=east] {7}
		(-0.5,3.2) node[anchor=east] {8}
		(-0.5,4.0) node[anchor=east] {10}
		(-0.5,4.4) node[anchor=east] {11}
		(-0.5,5.6) node[anchor=east] {14}
		(-0.5,6.0) node[anchor=east] {15};
		
\draw[circle,fill=black] (-0.1,0);	

\filldraw [fill=black] (-0.1,0) circle (2pt);
\filldraw [fill=black] (0.3,5.6) circle (2pt);
\filldraw [fill=black] (0.7,6.0) circle (2pt);
\filldraw [fill=black] (1.1,5.6) circle (2pt);
\filldraw [fill=black] (1.5,4.4) circle (2pt);
\filldraw [fill=black] (1.9,4.0) circle (2pt);
\filldraw [fill=black] (2.3,3.2) circle (2pt);
\filldraw [fill=black] (2.7,2.8) circle (2pt);
\filldraw [fill=black] (3.1,2.4) circle (2pt);
\filldraw [fill=black] (3.5,2.4) circle (2pt);
\filldraw [fill=black] (3.9,2.0) circle (2pt);
\filldraw [fill=black] (4.3,2.0) circle (2pt);
\filldraw [fill=black] (4.7,1.6) circle (2pt);
\filldraw [fill=black] (5.5,1.6) circle (2pt);
\filldraw [fill=black] (5.9,1.2) circle (2pt);
\filldraw [fill=black] (7.1,1.2) circle (2pt);
\filldraw [fill=black] (7.5,0.8) circle (2pt);
\filldraw [fill=black] (7.9,0.8) circle (2pt);
\filldraw [fill=black] (8.3,0.4) circle (2pt);
    	
\draw[dotted] (-0.1,0) -- (0.3,5.6);
\draw[dotted] (0.3,5.6) -- (0.7,6.0);
\draw[dotted] (0.7,6.0) -- (1.1,5.6);
\draw[dotted] (1.1,5.6) -- (1.5,4.4);
\draw[dotted] (1.5,4.4) -- (1.9,4.0);
\draw[dotted] (1.9,4.0) -- (2.3,3.2);
\draw[dotted] (2.3,3.2) -- (3.1,2.4);
\draw[dotted] (3.1,2.4) -- (3.5,2.4);
\draw[dotted] (3.5,2.4) -- (3.9,2.0);
\draw[dotted] (3.9,2.0) -- (4.3,2.0);
\draw[dotted] (4.3,2.0) -- (4.7,1.6);
\draw[dotted] (4.7,1.6) -- (5.5,1.6);
\draw[dotted] (5.5,1.6) -- (5.9,1.2);
\draw[dotted] (5.9,1.2) -- (7.1,1.2);
\draw[dotted] (7.1,1.2) -- (7.5,0.8);
\draw[dotted] (7.5,0.8) -- (7.9,0.8);
\draw[dotted] (7.9,0.8) -- (8.3,0.4);
\draw[->][dotted] (8.3,0.4) -- (10,0.4);
\end{tikzpicture}

\noindent The critical function above is thought for $q$ large, since $\mu_q(\F_2, |\cdot|,m)$ is just defined for $m \le q/2$. Although we have not shown yet how to construct this function, we could give better bounds for it at a fixed value of $q$. However, it will be crucial to work with a uniform bound in $q$, that is why the function above does not depend on this parameter. The same comments apply to the critical functions we shall use for triangular and cyclic groups below. In Theorem A1 ii), our aim is to treat all the $\F_n$'s together and get an uniform result in $n$. To do that, we need the critical function $\mu_q(\F_n,|\cdot|,m)$ to be uniformly bounded in $(n,m,q)$, so that we find finitely many forms of pathologies which can be fixed by hand. This is possible if we allow $q$ to be large enough. More precisely, for $q\geq q(n)=(22n)^{44n}+2$ the critical function is given by
$$\mu_q(\F_n,|\cdot|,m)=\left\{
\begin{array}{cl}
0 & \quad \mbox{ for } m=1, \\
2 & \quad \mbox{ for } 2 \leq m\leq s.
\end{array}\right.$$
The index $q(n)$ is very large, this is crucial to obtain a uniform result in $n$. On the contrary, for $n$ fixed, a more careful approach in terms of numerical computations could improve $q(n)$. We might even extend Theorem A1 i) to $\F_3,\F_4, \ldots$ and obtain $q(n)=4$ for these values of $n$. For the uniform result in Theorem A1 ii) we just consider $q\geq q(n)$, so that $\Lambda$-estimates are not needed to control regular terms when $q \le q(\F_n,|\cdot|)$. This is why $(\beta)$-estimates below are stated only for $\F_2$. Moreover, since the set of super-pathological terms is different for $\F_2$ when $q\geq 4$ and for $\F_n$ when $q\geq q(n)$, we will establish an $(\varepsilon)$-estimate for each case.

\vskip5pt

\noindent ${\bf (\alpha)}$ {\bf Coefficients for $u=2,3$.} We start with estimates of $C_{q,s_u}^{\mathrm{left}}[\k](r)$ for $u=2,3$, which we deduce from Proposition \ref{s23}. 
We establish these estimates in the general case $\F_n$. 

\begin{itemize}
\item[$\alpha$i)] $C_{q,s_{2}}^{\mathrm{left}}[(k,\underline{0})](r) \, = \, r^{2k}$.

\vskip8pt

\item[$\alpha$ii)] $\displaystyle C_{q,s_{3}}^{\mathrm{left}}[(k_1,k_2,\underline{0})](r) \, = \, \displaystyle \frac{3}{s(s-1)}M(k_1,k_2,\underline{0}) \Big( \sum_{m=0}^{\lfloor k_2/2 \rfloor} r^{2(k_1+k_2-m)} \Big)$.

\vskip8pt

\item[$\alpha$iii)] If $k  \le 4 \delta_{n=2} + 2 \delta_{n \ge 3}$, we obtain $$\hskip30pt C_{q,s_{3}}^{\mathrm{left}}[(k, \underline{0})](r) \, = \, \left\{\begin{array}{ll} 0 &\mbox{if} \ k = 1, \\ \frac32 r^4+(n-1)r^6 & \mbox{if} \ k = 2, \\ 3 \hskip1pt r^6 + 3(n-1) r^8 & \mbox{if} \ k = 3, \\ \frac92 r^8 \hskip-1pt + 6(n-1) r^{10} + (n-1)(2n-1) r^{12} & \mbox{if} \ k = 4. \end{array}\right.$$

\vskip5pt

\item[$\alpha$iv)] If $k \ge 5 \delta_{n=2} + 3 \delta_{n \ge 3}$, we obtain $$\hskip30pt C_{q,s_{3}}^{\mathrm{left}}[(k, \underline{0})](r) \, = \, \left\{\begin{array}{ll}
\frac{3 k (1-r^2) r^{2k}}{2 (1-(2n-1)r^2)} + \frac{6(n-1)r^{2k}}{(2n-1)(1-(2n-1)r^2)^2} & \mbox{when} \ r < \frac{1}{\sqrt{2n-1}}, \\ [6pt]  \frac{3(n-1)k^2+(30n-24)k+24(n-1)}{4(2n-1)^{k+1}} & \mbox{when} \ r= \frac{1}{\sqrt{2n-1}}.
\end{array}\right.$$
\end{itemize}

\vskip5pt

\noindent ${\bf (\beta)}$ {\bf Coefficients for $\Lambda$-regular terms -- $\F_2$.} Given $4\leq u \leq q$ and $\k \in  \mathrm{Reg}(s)$ we now compute those $C_{q,s_{u}}^{\mathrm{left},\mathbf{R}}[\k](r)$ which arise from the $\Lambda$-estimates established in Paragraph \ref{Lambdaest}. Recall that we use these estimates only for $4 \le q \le q(\F_2,|\cdot|)=10$. Since it suffices to consider the optimal $r = 1/ \sqrt{q-1}$, we need to estimate the left-coefficients in the range $\frac{1}{3}\leq  r \leq \frac{1}{\sqrt{3}}$. For $\k \in \mathrm{Reg}_m(s)$, we find   
\begin{itemize}
\item[$\beta$i)] If $2 \le m \le s$, we obtain
$$\hskip30pt C_{q,s_{2m}}^{\mathrm{left},\mathbf{R}}[\k](r) = \left\{\begin{array}{l} \displaystyle\frac{(2m)!(s-m)!}{s!}M(\k)2^{m-2}3^{-|\k|} \hfill \mbox{ if }  r = \mbox{$\frac{1}{3}$}, \\ \displaystyle\frac{(2m)!(s-m)!}{s!}M(\k)\displaystyle\frac{(3r^2)^{|\k|}}{1-r}\Big(\displaystyle\frac{4r}{3r-1}\Big)^{m-1} \hfill \mbox{ if } \mbox{$\frac{1}{3}$} < r < \mbox{$\frac{1}{\sqrt{3}}$}, \\ 16M(\k)\displaystyle\frac{|\k|3^{(1-|\k|)/2}}{\sqrt{3}-1} \hfill \mbox{ if }  r =\mbox{$\frac{1}{\sqrt{3}}$} \mbox{ and } s=m=2. \end{array}\right.$$

\vskip5pt

\item[$\beta$ii)] If $3 \le m \le s$, we obtain
$$C_{q,s_{2m-1}}^{\mathrm{left},\mathbf{R}}[\k](r) = \left\{\begin{array}{l} \displaystyle\frac{(2m-1)!(s-m)!}{s!}M(\k)2^{m-4}3^{-|\k|} \hfill \mbox{ if }  r = \mbox{$\frac{1}{3}$}, \\ \displaystyle\frac{(2m-1)!(s-m)!}{2(s!)}M(\k)\displaystyle\frac{(3r^2)^{|\k|}}{1-r}\Big(\displaystyle\frac{4r}{3r-1}\Big)^{m-2} \hfill \mbox{ if } \mbox{$\frac{1}{3}$} < r < \mbox{$\frac{1}{\sqrt{3}}$}. \end{array}\right.$$

\vskip5pt

\item[$\beta$iii)] If $3 \le m+1 \le s$, we obtain
$$C_{q,s_{2m+1}}^{\mathrm{left},\mathbf{R}}[\k](r) = \left\{\begin{array}{l} \displaystyle\frac{(2m+1)!(s-m)!}{s!}M(\k)2^{m-3}3^{-|\k|}
\hfill \mbox{ if }  r = \mbox{$\frac{1}{3}$}, \\ \displaystyle\frac{(2m+1)!(s-m)!}{2(s!)}M(\k)\displaystyle\frac{(3r^2)^{|\k|}}{1-r}\Big(\displaystyle\frac{4r}{3r-1}\Big)^{m} \hskip10pt \hfill \mbox{ if } \mbox{$\frac{1}{3}$} < r < \mbox{$\frac{1}{\sqrt{3}}$}. \end{array}\right.$$
\end{itemize}

\vskip5pt

\noindent ${\bf (\gamma)}$ {\bf Coefficients for $\Delta$-regular terms.} If $q > q(\F_n, |\cdot|) = 4n^2-4n+2$ and $\k \in \mathrm{Reg}(s)$, we use the estimates for regular terms proved in Paragraph \ref{Deltaest}. Let us be more precise. Given $s \ge 2n^2-2n+2$ and $0 \le r < \frac{1}{2n-1}$, we may rewrite Proposition \ref{method1} for $\k \in \mathrm{Reg}_m(s)$ as follows
\begin{itemize}
\item[$\gamma$i)] If $2 \le m \le s$, we obtain 
$$\hskip20pt C_{q,s_{2m}}^{\mathrm{left},\mathbf{R}}[\k](r) \, = \, \frac{(2m)!(s-m)!}{m! s!} \Big( \frac{2nr}{1-(2n-1)r} \Big)^{m-1} r^{|\k|+1} M(\k).$$

\item[$\gamma$ii)] If $3 \le m \leq s$, we obtain
$$\hskip30pt C_{q,s_{2m-1}}^{\mathrm{left},\mathbf{R}}[\k](r) \, = \, \frac{(2m-1)!(s-m)!}{(m-1)! s!} \Big( \frac{2nr}{1-(2n-1)r} \Big)^{m-2} r^{|\k|+1} M(\k).$$

\item[$\gamma$iii)] If $3 \le m+1 \le s$, we obtain 
$$\hskip18pt C_{q,s_{2m+1}}^{\mathrm{left},\mathbf{R}}[\k](r) \, = \, \frac{(2m+1)!(s-m)!}{(m+1)! s!} \Big( \frac{2nr}{1-(2n-1)r} \Big)^{m} r^{|\k|+1} M(\k).$$
\end{itemize}

\vskip5pt

\noindent ${\bf (\delta)}$ {\bf Coefficients for pathological terms.} To estimate the coefficients $C_{q,s_u}^{\mathrm{left},\mathbf{P}}[\k](r)$ for pathological terms, we will use again the $\Lambda$-estimates in Paragraph \ref{Lambdaest} in a slightly modified form. Let us introduce some terminology. Given $\k \in \mathrm{Pat}_m(s)$ set $$A_{2m}[\k] \, = \, \Big\{ \l \in \mathrm{Adm}_{2m} \; : \;  k_j \in \{\ell_{2j-1},\ell_{2j}\} \mbox{ for } 1\le j \le m \Big\} \setminus B_{2m}.$$ Similarly, for $\k \in \mathrm{Pat}_m'(s)$ we define $$A'_{2m}[\k] \, = \, \Big\{ \l \in \mathrm{Adm}'_{2m} \; : \;  k_j \in \{\ell_{2j-1},\ell_{2j}\} \mbox{ for } 1\le j \le m  \Big\} \setminus B'_{2m}.$$ It is crucial to observe that, thanks to our definition of the admissible lengths $\mathrm{Adm}_{u}$, the sets $A_{2m}[\k]$ and $A_{2m}'[\k]$ are finite. Given $\l \in \mathrm{Adm}_{2m}\cup \mathrm{Adm}'_{2m}$, we also consider the sets $P(\l)=|\{1\leq j\leq m\; : \; \ell_{2j-1}=\ell_{2j}\}|$.
  Given $0 \le r \le 1$, we find for $\k \in \mathrm{Pat}_m(s)$ the following estimates
\begin{itemize}
\item[$\delta$i)] If $2 \le m \le s$, we obtain 
$$\hskip18pt C_{q,s_{2m}}^{\mathrm{left},\mathbf{P}}[\k](r) \, = \, \sum_{\l \in A_{2m}[\k]}2^{P(\l)-m}\Big(\prod_{j=2}^{m}N_{\ell_{2j-1}}\Big)M(\l)r^{|\l|}.$$

\item[$\delta$ii)] If $3 \le m \leq s$, we obtain
$$\hskip15pt C_{q,s_{2m-1}}^{\mathrm{left},\mathbf{P}}[\k](r) \,  = \,  \sum_{\l \in A_{2m}'[\k]}\frac{2^{P(\l)-m}}{2m}\Big(\prod_{j=2}^{m-1}N_{\ell_{2j-1}}\Big)M(\l)r^{|\l|}.$$

\item[$\delta$iii)] If $3 \le m+1 \leq s$, we obtain 
$$\hskip25pt C_{q,s_{2m+1}}^{\mathrm{left},\mathbf{P}}[\k](r) \, = \,  \sum_{\l \in A_{2m+2}'[\k]}\frac{2^{P(\l)-m}}{4(m+1)}\Big(\prod_{j=2}^{m+1}N_{\ell_{2j-1}}\Big)M(\l)r^{|\l|}.$$
\end{itemize}

\vskip5pt

\noindent ${\bf (\varepsilon)}$ {\bf Coefficients for super-pathological terms -- $\F_2$.} Most of our estimates above for pathological terms will serve for our purposes. However, a few of them are not fine enough. Using a computer, we identify those $\k \in \mathrm{Pat}(s)$ which fail this test and call them super-pathological. These terms have one of the forms $(1,1,\underline{0}),\;(2,1,\underline{0}),\;(3,1,\underline{0}),\;(2,2,\underline{0}),\;(1,1,1,\underline{0}),\;(1,1,1,1,\underline{0})$ adding at the end as many zeros as needed to form $s$-tuples. In order to improve our estimates, we need better bounds for the sums $s_u[\l](r)$ with $\l$ belonging to the exceptional set $B_u$ where 
\begin{eqnarray*}
B_4 & = & \big\{ (1,1,1,1),(3,1,1,1),(2,2,1,1),(3,2,2,1), \\ &  & \hskip5pt (5,2,2,1),(5,3,3,1),(4,3,2,1),(2,2,2,2) \big\},\\
B_5 & = & \big\{ (2,1,1,1,1),(3,2,1,1,1),(2,2,2,1,1) \big\}, \ B_6 \ = \ \big\{ (1,1,1,1,1,1) \big\}, \\ B_7 & = & \big\{ (2,1,1,1,1,1,1) \big\}, \ B_8 \ = \ \big\{ (1,1,1,1,1,1,1,1) \big\},
\end{eqnarray*}
We have $B_u=\emptyset$ otherwise. Given $0\leq r \leq 1$, we find
$$\begin{array}{rclcl}
\varepsilon1)&&s_4[(1,1,1,1)](r)&\leq & 2r^4\alpha_1^2, \\ [2pt] 
\varepsilon2)&&s_4[(3,1,1,1)](r)&\leq & 9r^6\alpha_1^2+r^6 \alpha_1\alpha_3, \\ [2pt] 
\varepsilon3)&&s_4[(2,2,1,1)](r)&\leq &13 r^6 \alpha_1\alpha_2, \\ [2pt] 
\varepsilon4)&&s_4[(3,2,2,1)](r)&\leq &72 r^8 \alpha_1\alpha_3 +24 r^8 \alpha_2^2, \\ [2pt] 
\varepsilon5)&&s_4[(5,2,2,1)](r)&\leq &72 r^{10} \alpha_1\alpha_5 +24 r^{10} \alpha_2^2, \\ [2pt] 
\varepsilon6)&&s_4[(5,3,3,1)](r)&\leq &216 r^{12} \alpha_1\alpha_5 +24 r^{12} \alpha_3^2, \\ [2pt] 
\varepsilon7)&&s_4[(4,3,2,1)](r)&\leq &144 r^{10} \alpha_1\alpha_4 +48 r^{10} \alpha_2\alpha_3, \\ [2pt] 
\varepsilon8)&&s_4[(2,2,2,2)](r)&\leq &12 r^8 \alpha_2^2, \\ [2pt] \varepsilon9)&\widehat{f}(e)&s_5[(2,1,1,1,1)](r)&\leq &\frac{15}{2}r^6 \alpha_0\alpha_1\alpha_2 +\frac{65}{12}r^6 \alpha_1^3, \\ [2pt] \varepsilon10)&\widehat{f}(e)&s_5[(3,2,1,1,1)](r)&\leq & 160 r^8 \alpha_0\alpha_2\alpha_3 + 120 r^8 \alpha_1^3, \\ [2pt] \varepsilon11)&\widehat{f}(e)&s_5[(2,2,2,1,1)](r)&\leq &240 r^8 \alpha_0\alpha_1\alpha_2 + 20 r^8 \alpha_1\alpha_2^2, \\ [2pt] \varepsilon12)&&s_6[(1,1,1,1,1,1)](r)&\leq &4r^6 \alpha_1^3, \\ [2pt] \varepsilon13)&\widehat{f}(e)&s_7[(2,1,1,1,1,1,1)](r)&\leq &\frac{735}{16}r^8\alpha_0\alpha_1^3 + \frac{525}{32} r^8 \alpha_1^3 \alpha_2, \\ [2pt] \varepsilon14)&&s_8[(1,1,1,1,1,1,1,1)](r)&\leq &\frac{29}{3}r^8 \alpha_1^4.
\end{array}$$

\vskip5pt

\noindent ${\bf (\varepsilon')}$ {\bf Coefficients for super-pathological terms -- $\F_n$.} In the general case $\F_n$ the $(\delta)$-estimates above for pathological terms are not fine enough to treat all $n$'s together. Hence we need to be careful with all the pathological terms. Since $\mu_q(\F_n, |\cdot|,m) = 2$ for $2 \le m \le s$, the pathological terms have the form $(1,1,\cdots,1,\underline{0})$ and  $(2,1,\cdots,1,\underline{0})$ in $L_{m}(s)$ for all $2\leq m \leq s$. 
We need to estimate better the sums $s_u[\l](r)$ with $\l$ belonging to the set $B_u$, where
$$\begin{array}{rll}
B_{2m}=\{(1,1,\ldots,1)\} &\in \mathrm{Adm_{2m}} & \mbox{ for all } 2\leq m \leq s, \\
B_{2m+1}=\{(2,1,\ldots,1)\} &\in \mathrm{Adm_{2m+1}} & \mbox{ for all } 2\leq m \leq s-1. 
\end{array}$$
Given $0\leq r \leq 1$, we find for $2\leq m \leq s$
$$\begin{array}{rclcl}
\varepsilon'1) && s_{2m}[(1,\ldots,1)](r)&\leq & \displaystyle\frac{(2m)!}{m! (m+1)!}r^{2m}\alpha_1^m, \\ [8pt] 
\varepsilon'2) & \widehat{f}(e) & s_{2m+1}[(2,1,\ldots,1)](r)&\leq &  
\displaystyle\frac{3(2m+1)!}{2(m-1)! (m+2)!}r^{2m+2}\alpha_0\alpha_1^{m-1}\alpha_2 \\ [8pt] 
&&&+&n(2n-1)\displaystyle\frac{3(2m+1)!}{4(m-1)! (m+2)!}r^{2m+2}\alpha_1^{m+1}. 
\end{array}$$

\begin{remark}\label{RkCleftReg}
\emph{The super-pathological sums in $(\varepsilon)$ and $(\varepsilon')$ do not affect the regular terms $\k \in \mathrm{Reg}(s)$. 
Indeed, we have $|\k| \le 6$ and $|\k| \le m+1$ for the $\alpha_{\k}$'s in the right hand side of the $(\varepsilon)$ and $(\varepsilon')$ estimates respectively. Also, recall that $\k \in \mathrm{Reg}_m(s)$ in $\F_2$ when $$|\k| \ge m + \mu_q(\F_2, |\cdot|, m)\geq m+14$$ since $2 \le m \le 4$ for the $\k$'s in $(\varepsilon)$. Similarly, $\k \in \mathrm{Reg}_m(s)$ in $\F_n$ when $$|\k| \ge m + \mu_q(\F_n, |\cdot|, m)\geq m+2.$$ This means that for these terms we have in both cases $C_{q,s_{u}}^{\mathrm{left},\mathbf{S}}[\k](r)=0$.}
\end{remark}


\subsubsection{Proof of Theorem \emph{A1 i)}} \label{proofA1}

According to \eqref{AIM}, it remains to prove 
\begin{eqnarray*} 
\label{Cleft<Cright} C_q^{\mathrm{left}}[\k](r) & := & \sum_{u=2}^q \binom{q}{u} C_{q,s_{u}}^{\mathrm{left}}[\k](r) \\ \nonumber & = & \sum_{u=2 \vee (2m-1)}^{q \wedge (2m+1)} \binom{q}{u} C_{q,s_{u}}^{\mathrm{left}}[\k](r) \ \le \ C_q^{\mathrm{right}}[\k] \ = \ M(\k),
\end{eqnarray*}
for $\k \in L_m(s) \setminus \{\underline{0}\}$ and $0 \le r \le 1/\sqrt{q-1}$. The second identity follows easily from the process of completing squares, the only sums $s_u(r)$ which may have a contribution to $\alpha_{\k}$ with $\k \in L_m(s)$ are those with $u \in \{2m-1,2m,2m+1\}$. We divide the proof into regular and pathological cases.

\noindent {\bf The regular case $\k \in \mathrm{Reg}_m(s)$.} We have $$C_q^{\mathrm{left}}[\k](r) = C_q^{\mathrm{left},\mathbf{R}}[\k](r)$$ for any $\k \in \mathrm{Reg}(s)$. Indeed, it follows from \eqref{Cleft3}, \eqref{Reg-Pat} and Remark \ref{RkCleftReg}. When $m=1$, we have $\mathrm{Reg}_1(s) = L_1(s)$ since $\mu_q(\F_2,|\cdot|,1)=0$. Therefore, $\k$ is of the form $(k,\underline{0})$ for some $k \ge 1$ and we have to prove $$\binom{q}{2} C_{q,s_2}^{\mathrm{left}, \mathbf{R}}[\k](r) + \binom{q}{3} C_{q,s_3}^{\mathrm{left}, \mathbf{R}}[\k](r) \, \le \, s$$ for any $k \ge 1$ and $s \ge 2$. The case $k=1$ yields $\binom{q}{2} r^2 = s$ for $r=1/\sqrt{2s-1}$, which is sharp as announced in Paragraph \ref{Strategy}. If $k \ge 2$, according to estimates $\alpha$i), $\alpha$iii) and $\alpha$iv), the proof reduces to inequalities \eqref{Reg1-1}-\eqref{Reg1-4} in Appendix D1 for $n=2$, except for $k \ge 5$ and $s=2$ which follows by direct substitution. When $m=2$, regular terms $\k = (k_1, k_2, \underline{0})$ must satisfy $|\k|  \geq 2+\mu_q(\F_2,|\cdot|,2)=16$, so that we need the inequality $$\binom{q}{3} C_{q,s_3}^{\mathrm{left}, \mathbf{R}}[\k](r) + \binom{q}{4} C_{q,s_4}^{\mathrm{left}, \mathbf{R}}[\k](r) + \binom{q}{5} C_{q,s_5}^{\mathrm{left}, \mathbf{R}}[\k](r) \, \le \, M(\k),$$ with the usual convention that $\binom{q}{5}=0$ for $q < 5$. If $2 \leq s \leq 5$, we use $\alpha$ii), $\beta$i) and $\beta$iii) and the expected inequalities follow by direct substitution. If $s \ge 6$ we use $\alpha$ii), $\gamma$i) and $\gamma$iii), which reduces to inequality \eqref{Reg2} in Appendix D1 for $n=2$. When $3 \leq m \leq s$ and $\k \in \mathrm{Reg}_m(s)$ we need to prove $$\binom{q}{2m-1} C_{q,s_{2m-1}}^{\mathrm{left}, \mathbf{R}}[\k](r) + \binom{q}{2m} C_{q,s_{2m}}^{\mathrm{left}, \mathbf{R}}[\k](r) + \binom{q}{2m+1} C_{q,s_{2m+1}}^{\mathrm{left}, \mathbf{R}}[\k](r) \, \le \, M(\k).$$ If $3 \le m \le s \le 5$, we use $\beta$i), $\beta$ii) and $\beta$iii) to obtain the desired inequality by direct substitution noting that $|\k| \ge m + \mu_q(\F_2, |\cdot|,m)$. If $s\geq 6$ we use $\gamma$i), $\gamma$ii) and $\gamma$iii) to reduce the desired inequality to \eqref{Regm} in Appendix D1 for $n=2$. 

\noindent {\bf The pathological case $\k \in \mathrm{Pat}_m(s)$.} The pathological sets $\mathrm{Pat}_m(s)$ are all finite for fixed value of $s$ and $m$, but these values range over an infinite set since $q = 2s$ can be arbitrarily large. We will estimate the terms with $m \ge 22$ by hand. The other cases present finitely many classes of pathologies and will be estimated with computer assistance later on. Since $\mu_q(\F_2, |\cdot|, m) = 1$ when $22 \leq m \leq s$, we have $\mathrm{Pat}_m(s)=\{(1,1,\ldots,1,\underline{0})\}$. Since $B_u=\emptyset$ for $u\geq 9$, the super-pathological sums do not affect these $\k$'s. This means that $C_q^{\mathrm{left}}[\k](r)= C_q^{\mathrm{left},\mathbf{P}}[\k](r)$ is given by our $(\delta)$-estimates. By definition of the admissible lengths $\mathrm{Adm}_u$, we obtain
\begin{eqnarray*}
A'_{2m}[\k] & = & \big\{ (\ell,1,1,\ldots,1,0) \, : \, 2 \leq \ell \leq 2m-2, \, \ell \mbox{ even} \big\} \ \subset \ \mathrm{Adm}'_{2m}, \\
A_{2m}[\k] & = & \big\{ (\ell,1,1,\ldots,1,1) \, : \, 1 \leq \ell \leq 2m-1, \, \ell \mbox{ odd} \big\} \hskip1pt \ \subset \ \hskip2pt  \mathrm{Adm}_{2m}, \\ A'_{2m+2}[\k] & = & \big\{ (\ell,1,1,\ldots,1,0) \, : \, 2 \leq \ell \leq 2m, \, \ell \mbox{ even} \big\} \hskip4pt \ \subset \ \hskip4pt  \mathrm{Adm}'_{2m+2}.
\end{eqnarray*}
Using $\delta$i), $\delta$ii) and $\delta$iii) we get
\begin{eqnarray*}
\lefteqn{C_q^{\mathrm{left}}[\k](r) \, = \, \binom{2s}{2m}\Big(4^{m-1}r^{2m}+\sum_{\begin{subarray}{c} \ell=3 \\ \ell \mathrm{ odd} \end{subarray}}^{2m-1} m 4^{m-1}r^{2m-1+\ell} \Big)} \\ \!\!\!\! & + & \!\!\!\! \binom{2s}{2m-1} \sum_{\begin{subarray}{c} \ell=2 \\ \ell \mathrm{ even} \end{subarray}}^{2m-2} (2m-1)4^{m-3}r^{2m-2+\ell} \, + \, \binom{2s}{2m+1}\sum_{\begin{subarray}{c} \ell=2 \\ \ell \mathrm{ even} \end{subarray}}^{2m} (2m+1)4^{m-1}r^{2m+\ell} \\ \!\!\!\! & \leq & \!\!\!\! \Big[ \binom{2s}{2m-1} \frac{2m-1}{4^3} +\binom{2s}{2m} \frac{m}{4} + \binom{2s}{2m+1} \frac{2m+1}{4} r^2 \Big] \frac{4^m r^{2m}}{1-r^2}. 
\end{eqnarray*}
Since $C_q^{\mathrm{right}}[\k]=\binom{s}{m}$, we are reduced to prove \eqref{Pat} in Appendix D1. It remains to analyze the terms in $\mathrm{Pat}_m(s)$ with $2 \leq m \leq 21$. It is crucial to note that this set is finite. Indeed, using our $(\alpha)$, $(\delta)$ and $(\varepsilon)$ estimates we can express the left coefficients of the pathological terms in $\bigcup_{2\leq m \leq 21}\mathrm{Pat}_m(s)$. Thanks to the parity condition established in Lemma \ref{admlength-Fn}, in these left coefficients $r$ always appears with an even exponent. Hence, we are reduced to prove a finite number of inequalities for functions in the variable $r^2 = \frac{1}{2s-1}$, which become rational functions in $s$. Equivalently, by rearranging we have to prove that a finite number of polynomials in $s$ are positive. Since $s$ only takes integer values, this can be easily done by a computer via the positivity test explained in Appendix D. Only a small number of terms fail this test, the terms $\alpha_0^{s-2}\alpha_1^2$ for $2 \le s \le 6$. 

In this case, we split the coefficient $\alpha_1$ and go back to the coefficients $a_1(1)$ and $a_1(2)$ which we found in the process of completing squares. More precisely, the sums potentially contributing to $\k=(1,1,\underline{0}) \in \mathrm{Pat}_2(s)$ are $s_3(r)$ and ---according to the sets $A_{2m}[\k]$ and $A_{2m+2}'[\k]$--- $s_4[(1,1,1,1)](r)$, $s_4[(3,1,1,1)](r)$, $s_5[(2,1,1,1,1)](r)$, $s_5[(4,1,1,1,1)](r)$. According to Proposition \ref{s23}, the contribution of $s_3(r)$ is given by $$\Lambda_3 \, = \, \binom{2s}{3}\Big(\frac{3}{4}r^4\alpha_1^2+6 r^4a_1(1)^2a_1(2)^2 \Big).$$ On the other hand, according to Appendix C, we show with our estimates of $\varepsilon$1) and $\varepsilon$2) ---not the final bound in terms of $\alpha_1^2$, but the previous ones--- that the contribution of the $s_4$-sums above to $\alpha_1^2$ is dominated by $$\Lambda_4 \, = \, \binom{2s}{4}\Big( \frac{3}{2}r^4\alpha_1^2+4r^4 a_1(1)^2a_1(2)^2+ 36 r^6 (a_1(1)^4+a_1(2)^4)\Big).$$ Finally, our way of completing squares for $\varepsilon$9) does not give any contribution for $s_5[(2,1,1,1,1)](r)$. According to $\delta$iii), the last sum $s_5[(4,1,1,1,1)](r)$ contributes as $\Lambda_5 = \binom{2s}{5}20 r^8 \alpha_1^2$. The goal is to show that $$\Lambda_3 + \Lambda_4 + \Lambda_5 \, \le \, \frac{s(s-1)}{2}\alpha_1^2,$$ where $\alpha_1=2(a_1(1)^2+a_1(2)^2)$. 
We may rewrite $\Lambda_3 + \Lambda_4 + \Lambda_5$ as follows
$$A(s)\alpha_1^2+4B(s)(a_1(1)^4+a_1(2)^4)+8C(s)a_1(1)^2a_1(2)^2 \leq [A(s)+\max\{B(s),C(s)\}]\alpha_1^2,$$
with 
\begin{eqnarray*}
A(s) & = & \frac{s(s-1)(24s^3-4s^2-94s+93)}{12(2s-1)^3}, \\ 
B(s) & = & \frac{3s(s-1)(2s-3)}{2(2s-1)^2}, \\ C(s) & = & \frac{s(s-1)(2s+3)}{12(2s-1)}.
\end{eqnarray*}
The assertion follows from $A(s)+\max\{B(s),C(s)\}\leq \frac{s(s-1)}{2}$ for $2\leq s\leq 6$. \fin

\subsubsection{Proof of Theorem \emph{A1 ii)}}

The argument will follow the same steps as the proof of Theorem A1 i) detailed in Paragraph \ref{proofA1} by using the estimates $(\alpha)$, $(\gamma)$, $(\delta)$ and $(\varepsilon')$ for the free group $\F_n$ with $n$ generators. 

\noindent {\bf The regular case $\k \in \mathrm{Reg}_m(s)$.} 
This is similar to the regular case for $\F_2$. It is even simpler since in our situation we assume $q\geq q(n)$ and we do not need the $\Lambda$-estimates. When $\k \in \mathrm{Reg}_1(s)=L_1(s)$, it reduces to prove inequalities \eqref{Reg1-1} and \eqref{Reg1-4} in Appendix D, the case $k=1$ being trivial. When $\k \in \mathrm{Reg}_2(s)$ with $|\k| \geq 2+\mu_q(\F_n,|\cdot|,2)=4$ we use $\alpha$ii), $\gamma$i) and $\gamma$iii) to end up with inequality \eqref{Reg2} in Appendix D. 
When $3\leq m \leq s$ and $\k \in \mathrm{Reg}_m(s)$, we reduce the desired inequality to \eqref{Regm} in Appendix D by means of $\gamma$i), $\gamma$ii) and $\gamma$iii). 

\noindent {\bf The pathological case $\k \in \mathrm{Pat}_m(s)$.} In the general case $\F_n$ we will fix by hand all the pathological terms 
$\k \in \mathrm{Pat}_m(s)=\{(1,1,\ldots,1,\underline{0}),(2,1,\ldots,1,\underline{0})\}$ for $2\leq m \leq s$. 
We start with $\k=(1,1,\ldots,1,\underline{0}) \in \mathrm{Pat}_m(s)$ for $3\leq m \leq s$. 
By using the description of the sets $A'_{2m}[\k],A_{2m}[\k]$ and $A'_{2m+2}[\k]$ given in Paragraph \ref{proofA1} together with the estimates $(\delta)$ and $\varepsilon'1)$ + $\varepsilon'2)$, we obtain 
\begin{eqnarray*}
\lefteqn{C_q^{\mathrm{left}}[\k](r) }\\ 
& = & 
\binom{2s}{2m-1} \Big(\frac{3n(2n-1)(2m-1)!}{4(m-2)!(m+1)!}r^{2m}+ \sum_{\begin{subarray}{c} \ell=4 \\ \ell \ \mathrm{even} \end{subarray}}^{2m-2} (2m-1)\frac{(2n)^{m-2}}{4}r^{2m-2+\ell}\Big) \\ 
& + & \binom{2s}{2m}\Big(\frac{(2m)!}{m!(m+1)!}r^{2m}+\sum_{\begin{subarray}{c} \ell=3 \\ \ell \ \mathrm{odd} \end{subarray}}^{2m-1} m (2n)^{m-1}r^{2m-1+\ell} \Big) \\ & + & 
\binom{2s}{2m+1}\sum_{\begin{subarray}{c} \ell=4 \\ \ell \ \mathrm{even} \end{subarray}}^{2m} (2m+1)\frac{(2n)^{m}}{4}r^{2m+\ell}. 
\end{eqnarray*}
The term for $\ell=2$ in the last sum vanishes through the use of $\varepsilon'2)$. Recalling that $C_q^{\mathrm{right}}[\k]=\binom{s}{m}$, we are reduced to \eqref{Pat-1-Fn} in Appendix D. For $\k=(1,1,\underline{0})$ in $\mathrm{Pat}_2(s)$, putting together similarly $\alpha$ii) + $\varepsilon'1)$ \& $\delta$i) + $\varepsilon'2)$ \& $\delta$iii) yields \eqref{Pat-2-Fn} in Appendix D. We now turn to $\k=(2,1,\cdots,1,\underline{0}) \in \mathrm{Pat}_m(s)$ for $3\leq m \leq s$. In that case, the definition of admissible lengths yields the following descriptions
\begin{eqnarray*}
A'_{2m}[\k] & = & \big\{ (2,1,\ldots,1,0)\big\} 
\\ & & \cup \big\{ (\ell,2,1,\ldots,1,0) \, : \, 3 \leq \ell \leq 2m-1, \, \ell \mbox{ odd} \big\} \\
&& \cup \big\{ (\ell,2,2,1,1,\ldots,1,0) \, : \, 2 \leq \ell \leq 2m, \, \ell \mbox{ even} \big\}
\hskip15pt \subset \ \mathrm{Adm}'_{2m}, \\
A_{2m}[\k] & = & \big\{ (\ell,2,1,\ldots,1) \, : \, 2 \leq \ell \leq 2m, \, \ell \mbox{ even} \big\} \\
&& \cup \big\{ (\ell,2,2,1,\ldots,1) \, : \, 3 \leq \ell \leq 2m+1, \, \ell \mbox{ odd} \big\} \hskip20pt  \subset \ \mathrm{Adm}_{2m}, \\ 
A'_{2m+2}[\k] & = & \big\{ (2,1,\ldots,1,0)\big\} 
\\ & & \cup \big\{ (\ell,2,1,\ldots,1,0) \, : \, 3 \leq \ell \leq 2m+1, \, \ell \mbox{ odd} \big\} \\
&& \cup \big\{ (\ell,2,2,1,\ldots,1,0) \, : \, 2 \leq \ell \leq 2m+2, \, \ell \mbox{ even} \big\} \hskip8pt \subset \ \mathrm{Adm}'_{2m+2}.
\end{eqnarray*}
By using $(\delta)$ and $\varepsilon'2)$ we find 
$$C_q^{\mathrm{left}}[\k](r) = \binom{2s}{2m-1} A_m(r,s) + \binom{2s}{2m} B_m(r,s)+ \binom{2s}{2m+1} C_m(r,s),$$ where
\begin{eqnarray*}
A_m(r,s) & = &  \sum_{\begin{subarray}{c} \ell=3 \\ \ell \mathrm{ odd } \end{subarray}}^{2m-1} \frac{(2n)^{m-2}}{4}(2m-1)(2m-2)r^{2m-1+\ell} \\ & + & \sum_{\begin{subarray}{c} \ell=2 \\ \ell \mathrm{ even } \end{subarray}}^{2m} \frac{(2n)^{m-2}(2n-1)}{8}(2m-1)(2m-2)(2m-3)\Big(\frac{1}{3}\delta_{\ell=2}+\frac{1}{2}\delta_{\ell\neq 2}\Big)r^{2m+\ell} \\ B_m(r,s) & = & 
\sum_{\begin{subarray}{c} \ell=2 \\ \ell \mathrm{ even } \end{subarray}}^{2m} (2n)^{m-1}m(2m-1)r^{2m+\ell} \\
& + & \sum_{\begin{subarray}{c} \ell=3 \\ \ell \mathrm{ odd } \end{subarray}}^{2m+1}  \frac{(2n)^{m-1}(2n-1)}{4}m(2m-1)(2m-2)r^{2m+1+\ell} \\ C_m(r,s) & = & \frac{3(2m+1)!}{2(m-1)!(m+2)!}r^{2m+2} + \sum_{\begin{subarray}{c} \ell=3 \\ \ell \mathrm{ odd } \end{subarray}}^{2m+1} \frac{(2n)^{m}}{2}m(2m+1)r^{2m+1+\ell} \\ & + & \sum_{\begin{subarray}{c} \ell=2 \\ \ell \mathrm{ even } \end{subarray}}^{2m+2} \frac{(2n)^{m}(2n-1)}{8}(2m+1)(2m)(2m-1)\Big(\frac{1}{3}\delta_{\ell=2}+\frac{1}{2}\delta_{\ell\neq 2}\Big)r^{2m+2+\ell}.
\end{eqnarray*}
Since  $C_q^{\mathrm{right}}[\k]=\frac{s!}{(m-1)!(s-m)!} = \frac{1}{m} \binom{s}{m}$, this yields \eqref{Pat-3-Fn} in Appendix D. 
For $\k=(2,1,\underline{0}) \in \mathrm{Pat}_2(s)$, using $\alpha$ii) and \eqref{Pat-3-Fn} for $m=2$ we get \eqref{Pat-4-Fn}. 
This ends the proof of the pathological case and completes the proof of Theorem A1 ii). \fin

\subsection{Triangular groups}

We now apply the combinatorial method presented in Section \ref{Method} to another natural example of finitely generated group equipped with the words length $|\cdot|$, the triangular groups
$$\Delta_{\alpha\beta\gamma} \, = \, \Big\langle a,b,c \, \big| \, a^2 = b^2 = c^2 = (ab)^\alpha = (bc)^\beta = (ca)^\gamma = e \Big\rangle.$$ 
In the spirit of this paper, we will only highlight the main steps in the proof of Theorem A2 and collect all the technical computations in Appendices C and D. 
We denote by $L:=2\min(\alpha,\beta,\gamma)$ the length of the smallest loop in $\Delta_{\alpha\beta\gamma}$, and we will see that for technical reasons we need to avoid small values of $L$ (namely $L\geq 15$) to make our argument work. 
In that situation we have 
$$ N_k \leq 3\cdot 2^{k-1}, \quad \mathcal{G}(\Delta_{\alpha\beta\gamma},|\cdot|,r) \leq \frac{3r}{1-2r}
\quad \mbox{and} \quad q(\Delta_{\alpha\beta\gamma},|\cdot|)=5.$$

\subsubsection{Admissible lengths}

Observe that all possible loops in the group have even length. Therefore, Lemma \ref{admlength-Fn} still holds true in the case of triangular groups and we may also refine the set of admissible lengths in that situation by adding a parity condition and setting for $2\leq u \leq q$ $$\mathrm{Adm}_u \, = \, \Big\{ \l \in L(u) \; : \; \ell_1 \ge \ell_2 \ge \ldots \ge \ell_u \ge 1, \; \sum_{j=1}^{u} \ell_j \ \mbox{even}, \; \ell_1 \le \sum_{j=2}^{u} \ell_j \Big\}.$$ 

\subsubsection{Estimates for $s_3(r)$} 

The necessity of avoiding small loops appears in the estimate of the $s_3$-sum below. Indeed, when we are below the smallest loop around the origin, the group behaves like a free group and we can obtain a finer estimate for $s_3(r)$. However, when we are above that smallest loop, some cancellations may appear and it becomes more difficult to estimate this sum. Fortunately, this may be compensated by the decay of the semigroup, which allows us to use somehow brutal estimates for lengths beyond certain relatively small quantity which we impose to be the smallest loop length. More precisely, we will choose $L$ such that the following estimates imply $\mu_q(\Delta_{\alpha\beta\gamma},|\cdot|,1)=0$.  

\begin{proposition} \label{s23-triangle}
We have $$\widehat{f}(e) s_3(r)  \le  \sum_{ \ell_2 \ge \ell_3 \ge 1} \sum_{m=0}^{\lfloor \ell_3/2 \rfloor} \big[A_m(\ell_2, \ell_3,r) \alpha_{\ell_2} \alpha_{\ell_3} + B_m(\ell_2,\ell_3,r) \alpha_0 \alpha_{\ell_2+\ell_3-2m} \big],$$
where we have for $K = K(\ell_2, \ell_3,L) = \frac{3(\ell_2+\ell_3) - L}{6}$
\begin{eqnarray*}
\displaystyle A_m(\ell_2,\ell_3,r) & = & \frac12 M(\ell_2+\ell_3-2m,\ell_2,\ell_3) r^{2(\ell_2+\ell_3-m)}, \\ [8pt] 
\displaystyle B_m(\ell_2, \ell_3,r) & = & \frac12 N_{\ell_3} M(\ell_2+\ell_3-2m,\ell_2,\ell_3) r^{2(\ell_2+\ell_3-m)} \delta_{m \leq K},\\ [8pt] & + & \frac12( \delta_{m=0} + 2^{m-1} \delta_{m>0}) M(\ell_2+\ell_3-2m,\ell_2,\ell_3) r^{2(\ell_2+\ell_3-m)} \delta_{m > K}.
\end{eqnarray*}
\end{proposition}

\dem
We follow the argument in the proof of Proposition \ref{s23}, by modifying only the estimate c) of $|\Lambda_3(\d,h)|$.  
Recall that for $\l \in \mathrm{Adm}_3$, $\d \sim \l$ and $h \in W_{\ell_1}$ we define $$\Lambda_3(\d) \, = \, \Big\{ (g_1, g_2, g_3) \in \Delta_{\alpha\beta\gamma}^3 \, : \ g_1 g_2 g_3 = e \ \mbox{and} \ |g_j| = d_j \Big\},$$ $$\Lambda_3(\d,h) \, = \, \Big\{ (g_1, g_2, g_3) \in \Lambda_3(\d) \, : \ g_{\sigma_{\d}^{-1}(1)} = h \Big\}.$$ When $\ell_1 <L/3$, the elements $(g_1,g_2,g_3)\in \Delta_{\alpha\beta\gamma}^3$ considered satisfy $g_1g_2g_3=e$ and $|g_j| <L/3$ for all $1\leq j \leq 3$. 
Hence no loops could appear, and the only possible cancellations are the usual ones, as in the free case. In that situation c) becomes 
\begin{itemize}
\item[$\mathrm{c}'$)] If $|h| = \ell_2 + \ell_3 - 2m$, $|\Lambda_3(\d,h)| = \delta_{m=0} + 2^{m-1} \delta_{m>0}$,
\end{itemize}
whenever $\ell_1 =\ell_2+\ell_3-2m<L/3 \Leftrightarrow m > K(\ell_2,\ell_3,L)$. If $\ell_1\geq L/3$, we may have extra cancellations coming from the loops in the group, and we cannot estimate $|\Lambda_3(\d,h)|$ precisely. We will use the trivial inequality $|\Lambda_3(\d,h)| \leq N_{\ell_3}$ whenever $m \leq K(\ell_2,\ell_3,L)$. Putting all together yields the required result. \fin

\subsubsection{Numerical estimates for $\Delta_{\alpha\beta\gamma}$}

We keep the notation $$C_{q,s_u}^{\mathrm{left},\mathbf{R}}[\k](r), \quad C_{q,s_u}^{\mathrm{left},\mathbf{P}}[\k](r), \quad C_{q,s_u}^{\mathrm{left},\mathbf{S}}[\k](r)$$ introduced in Paragraph \ref{numericalF2}. Then \eqref{Reg-Pat} still holds true, and we will see that the super-pathological sums do not affect the regular terms. We collect below the formulas for these left coefficients arising from the results detailed in the general method in Section \ref{Method}. We refer to Appendix C for the details of the proofs. In the case of a triangular group $\Delta_{\alpha\beta\gamma}$ with $L\geq 15$ equipped with the word length, the critical function is given by the graph

\begin{tikzpicture}

\draw[->] (-0.5,0) -- (9,0) node[anchor=north] {$m$};
\draw	(-0.5,0) node[anchor=north] {0}
		(0,0) node[anchor=north] {1}
		(0.5,0) node[anchor=north] {2}
		(1,0) node[anchor=north] {3}
		(1.5,0) node[anchor=north] {4}
		(2,0) node[anchor=north] {5}
		(2.5,0) node[anchor=north] {6}
		(3,0) node[anchor=north] {7}
		(3.5,0) node[anchor=north] {8}
		(5,0) node[anchor=north] {11}		
		(5.5,0) node[anchor=north] {12}
		(6.5,0) node[anchor=north] {14}
		(7,0) node[anchor=north] {15};
\draw	(4.5,-1) node{{\scriptsize {\bf The critical function} $\mu_q(\Delta_{\alpha\beta\gamma}, |\cdot|, m)$}};

\draw[->] (-0.5,0) -- (-0.5,7) node[anchor=east] {$\mu_q(m)$};

\draw 	(-0.5,0.5) node[anchor=east] {1}
		(-0.5,1) node[anchor=east] {2}
		(-0.5,1.5) node[anchor=east] {3}
		(-0.5,2) node[anchor=east] {4}
		(-0.5,2.5) node[anchor=east] {5}
		(-0.5,3) node[anchor=east] {6}
		(-0.5,6.5) node[anchor=east] {13};
		
\draw[circle,fill=black] (-0.1,0);	

\filldraw [fill=black] (-0.1,0) circle (2pt);
\filldraw [fill=black] (0.4,6.5) circle (2pt);
\filldraw [fill=black] (0.9,3) circle (2pt);
\filldraw [fill=black] (1.4,2.5) circle (2pt);
\filldraw [fill=black] (1.9,2.5) circle (2pt);
\filldraw [fill=black] (2.4,2) circle (2pt);
\filldraw [fill=black] (2.9,2) circle (2pt);
\filldraw [fill=black] (3.4,1.5) circle (2pt);
\filldraw [fill=black] (4.9,1.5) circle (2pt);
\filldraw [fill=black] (5.4,1) circle (2pt);
\filldraw [fill=black] (6.4,1) circle (2pt);
\filldraw [fill=black] (6.9,0.5) circle (2pt);
    	
\draw[dotted] (-0.1,0) -- (0.4,6.5);
\draw[dotted] (0.4,6.5) -- (0.9,3);
\draw[dotted] (0.9,3) -- (1.4,2.5);
\draw[dotted] (1.4,2.5) -- (1.9,2.5);
\draw[dotted] (1.9,2.5) -- (2.4,2);
\draw[dotted] (2.4,2) -- (2.9,2);
\draw[dotted] (2.9,2) -- (3.4,1.5);
\draw[dotted] (3.4,1.5) -- (4.9,1.5);
\draw[dotted] (4.9,1.5) -- (5.4,1);
\draw[dotted] (5.4,1) -- (6.4,1);
\draw[dotted] (6.4,1) -- (6.9,0.5);
\draw[->][dotted] (6.9,0.5) -- (9,0.5);

\end{tikzpicture}

\noindent ${\bf (\alpha)}$ {\bf Coefficients for $u=2,3$.} 
We find
\begin{itemize}
\item[$\alpha$i)] $C_{q,s_{2}}^{\mathrm{left}}[(k,\underline{0})](r) \, = \, r^{2k}$.

\vskip5pt

\item[$\alpha$ii)] $C_{q,s_{3}}^{\mathrm{left}}[(k_1,k_2,\underline{0})](r) \, = \, \displaystyle \frac{3}{s(s-1)}M(k_1,k_2,\underline{0}) \displaystyle\frac{r^{k_1+k_2}}{1-r^2}$.

\vskip8pt

\item[$\alpha$iii)] If $1 \le k \le 3$, we obtain $$C_{q,s_{3}}^{\mathrm{left}}[(k, \underline{0})](r) \, = \, \left\{\begin{array}{ll}
0 &\mbox{if} \ k = 1, \\ [3pt] \frac32 r^4+\frac12 r^6 & \mbox{if} \ k = 2, \\  [3pt] 3 r^6 + \frac32 r^8 & \mbox{if} \ k = 3. \end{array}\right.$$

\item[$\alpha$iv)] If $k \ge 4$, we get 
\begin{eqnarray*}
\lefteqn{\hskip10pt C_{q,s_{3}}^{\mathrm{left}}[(k, \underline{0})](r) \ = \ \frac{3r^{2k}}{2(1-2r^2)}\Big(k(1-r^2)+\frac{1}{1-2r^2}\Big) \delta_{k < L/3}} \\ & + & \Big[ \big( 27r^{10} + 45r^{12} + 36r^{14} \big) \delta_{k=5} + \frac{9r^{2k} 2^{k/2}(3-4r^2)}{2(1-r^2)(1-2r^2)} \delta_{k \ge 6} \Big] \delta_{k \ge L/3}.
\end{eqnarray*}
\end{itemize} 
\vskip5pt

\noindent ${\bf (\beta)}$ {\bf Coefficients for $\Lambda$-regular terms.} 
We will use the $\Lambda$-estimates for regular terms only in the case $q=4$ $(r= 1/\sqrt{3})$ since $q(\Delta_{\alpha\beta\gamma},|\cdot|)=5$. In that situation we find the following coefficients for $\k \in \mathrm{Reg}_2(s)$   
\begin{itemize}
\item[$\beta$)] $\displaystyle C_{q,s_{4}}^{\mathrm{left},\mathbf{R}}[\k](r)=M(\k)\frac{36\sqrt{3}}{(2-\sqrt{3})(\sqrt{3}-1)}\Big(\frac{2}{3}\Big)^{|\k|}$.
\end{itemize}

\vskip5pt

\noindent ${\bf (\gamma)}$ {\bf Coefficients for $\Delta$-regular terms.} If $q > q(\Delta_{\alpha\beta\gamma}, |\cdot|) = 5$ and $\k \in \mathrm{Reg}(s)$, we use the $\Delta$-estimates. Given $s \ge 3$ and $0 \le r < \frac12$, we may rewrite estimates $\gamma$i), $\gamma$ii) and $\gamma$iii) for the free group replacing the sum $\mathcal{G}(\G,\psi,r) = 2nr/(1-(2n-1)r)$ there by our estimate $3r/(1-2r)$.

\vskip5pt

\noindent ${\bf (\delta)}$ {\bf Coefficients for pathological terms.} The coefficients $C_{q,s_u}^{\mathrm{left},\mathbf{P}}[\k](r)$ for pathological terms follow from the formulas established in Paragraph \ref{numericalF2} for $\F_n$ deduced from the $\Lambda$-estimates, with $N_k\leq 3\cdot 2^{k-1}$ in that case. Note that the inequality $$N_{\ell_{\xi_j^\star}} \le N_{\ell_{2j-1}}$$ (used in the proof, see Appendix C) does not necessarily hold since $N_k$ might not be increasing. However, our upper bounds $3\cdot 2^{k-1}$ are increasing, which is enough. 

\vskip5pt

\noindent ${\bf (\varepsilon)}$ {\bf Coefficients for super-pathological terms.} By implementing in a computer the formulas for pathological terms $(\delta)$ above, we identify the super-pathological terms which fail the required inequality. These are of the forms $(1,1,\underline{0})$, $(2,1,\underline{0})$, $(1,1,1,\underline{0})$ adding at the end as many zeros as needed to form $s$-tuples. Hence, we need more careful estimates for the sums $s_u[\l](r)$ with $\l$ belonging to the exceptional set $B_u$ where
$$B_4  =  \big\{ (1,1,1,1),(3,1,1,1) \big\}, \quad B_5  =  \big\{ (2,1,1,1,1) \big\}, \quad B_6  =  \big\{ (1,1,1,1,1,1) \big\}.$$
We have $B_u=\emptyset$ otherwise. Given $0\leq r \leq 1$, we find
$$\begin{array}{rclcl}
\varepsilon \mbox{i)}&&s_4[(1,1,1,1)](r)&\leq & 2r^4\alpha_1^2, \\ [2pt] 
\varepsilon \mbox{ii)}&&s_4[(3,1,1,1)](r)&\leq & 3r^6\alpha_1^2+r^6 \alpha_1\alpha_3, \\ [2pt] 
\varepsilon \mbox{iii)}&\widehat{f}(e)&s_5[(2,1,1,1,1)](r)&\leq &\frac{45}{2}r^6 \alpha_0\alpha_1\alpha_2 +\frac{15}{2}r^6 \alpha_1^3, \\ [2pt] 
\varepsilon \mbox{iv)}&&s_6[(1,1,1,1,1,1)](r)&\leq &5r^6 \alpha_1^3.
\end{array}$$

\begin{remark}\label{RkCleftReg-triangle}
\emph{Note that $|\k| \le 4$ for the $\alpha_{\k}$'s in the right hand side of the estimates above. Also, recall that $\k \in L_m(s)$ is regular when $|\k| \ge m + \mu_q(\Delta_{\alpha\beta\gamma}, |\cdot|, m)$. Since $2 \le m \le 3$ for the $\k$'s above and $\mu_q(\Delta_{\alpha\beta\gamma}, |\cdot|, m) \ge 6$ in that range, it is clear that the super-pathological sums do not affect the regular terms $\k \in \mathrm{Reg}(s)$. This means that for these terms we have $C_{q,s_{u}}^{\mathrm{left},\mathbf{S}}[\k](r)=0$.}
\end{remark}

\subsubsection{Proof of Theorem \emph{A2}}

As we did for the free group in Paragraph \ref{proofA1}, it remains to put together all the estimates established in $(\alpha)-(\varepsilon)$ above in order to prove for all $\k \in L_m(s)\setminus \{0\}$ $$C_q^{\mathrm{left}}[\k](r) = \sum_{u = 2 \vee (2m-1)}^{q \wedge (2m+1)} \binom{q}{u} C_{q,s_u}^{\mathrm{left}}[\k](r) \, \le \, C_q^{\mathrm{right}}[\k]=M(\k).$$ 

\noindent {\bf The regular case $\k \in \mathrm{Reg}_m(s)$.} By Remark \ref{RkCleftReg-triangle}, we will only use $(\alpha), (\beta)$ and $(\gamma)$ to fix the regular case. When $m=1$ and $\k =(k, \underline{0}) \in \mathrm{Reg}_1(s)=L_1(s)$ we need to prove the inequalities \eqref{Reg1-1-triangle}-\eqref{Reg1-5-triangle} in Appendix D, the case $k=1$ is trivial. It is crucial to observe that to make the inequalities \eqref{Reg1-4-triangle} and \eqref{Reg1-5-triangle} true for all $s\geq 2$, which corresponds to the case $k \geq L/3$, we need the condition $L\geq 15$. When $m=2$, we have $\k=(k_1,k_2,\underline{0}) \in \mathrm{Reg}_2(s)$ if $|\k| \geq 2+\mu_q(\Delta_{\alpha\beta\gamma},|\cdot|,2)=15$. If $s=2$ $(r= 1/\sqrt{3})$ we use $\alpha$ii) and $\beta)$ and the desired inequality follows by direct substitution. If $s \ge 3$ $(0 < r < \frac{1}{2})$ we use $\alpha$ii), $\gamma$i) and $\gamma$iii). The resulting inequality is \eqref{Reg2-triangle} in Appendix D. When $3\leq m \leq s$ and $\k \in \mathrm{Reg}_m(s)$, we reduce the desired inequality to \eqref{Regm-triangle} in Appendix D by means of $\gamma$i), $\gamma$ii) and $\gamma$iii). 

\noindent {\bf The pathological case $\k \in \mathrm{Pat}_m(s)$.} We proceed as for $\F_2$. We first fix by hand the terms of the form $\k=(1,1,\ldots,1,\underline{0}) \in \mathrm{Pat}_m(s)$ for $m\geq 15$, then we use computer assistance to treat the finitely many remaining cases. Since the parity condition still holds in this setting, we can similarly describe the sets $A'_{2m}[\k], A_{2m}[\k]$ and $A'_{2m+2}[\k]$. Using $(\delta)$, this yields for $0\leq r <1$
\begin{eqnarray*}
\lefteqn{C_q^{\mathrm{left}}[\k](r) \, = \, \binom{2s}{2m} 3^{m-1}\Big(r^{2m}+\sum_{\begin{subarray}{c} \ell=3 \\ \ell \ \mathrm{odd} \end{subarray}}^{2m-1} m r^{2m-1+\ell} \Big)} \\ \!\!\!\! & + & \!\!\!\! \binom{2s}{2m-1} \sum_{\begin{subarray}{c} \ell=2 \\ \ell \ \mathrm{even} \end{subarray}}^{2m-2} \frac{2m-1}{4} 3^{m-2}r^{2m-2+\ell} \, + \, \binom{2s}{2m+1} \sum_{\begin{subarray}{c} \ell=2 \\ \ell \ \mathrm{even} \end{subarray}}^{2m} \frac{2m+1}{4}3^{m}r^{2m+\ell} \\ [5pt] \!\!\!\! & \leq & \!\!\!\! \Big[ \binom{2s}{2m-1} \frac{2m-1}{36} + \binom{2s}{2m} \frac{m}{3} + \binom{2s}{2m+1} \frac{2m+1}{4}r^2 \Big] \frac{(3r^2)^{m}}{1-r^2}. 
\end{eqnarray*}
Since $C_q^{\mathrm{right}}[\k]=\binom{s}{m}$, we are reduced to prove inequality \eqref{Pat-triangle} in Appendix D. We conclude as in the free group case, by implementing $(\alpha),(\delta)$ and $(\varepsilon)$ in a computer to fix the terms in $\bigcup_{2\leq m \leq 14}\mathrm{Pat}_m(s)$. All these terms verify the test. \fin 

\subsection{Finite cyclic groups}

Although our combinatorial method could be used in the non-Markovian setting (see below), this is not the case of the extrapolation result in Appendix A. Therefore, it will be useful to know whether the word length for $\Z_n$ is conditionally negative, since we could not find it in the literature. Our argument in Appendix B could be of independent interest. It is easily checked that finite cyclic groups satisfy $$N_k \le 2, \quad \mathcal{G}(\Z_n, |\cdot|,r) \le \frac{2r}{1-r} \quad \mbox{and} \quad q(\Z_n, |\cdot|) = 2.$$
 
\subsubsection{Admissible lengths} \label{Adm-Cyclic}

Let us write $$\Z_n = \{e, g, g^2, \ldots, g^{n-1}\},$$ so that any element in $\Z_n$ has the form $g^{\pm \ell}$ with $0 \le \ell \le \lfloor \frac{n}{2} \rfloor$. Given $g_1, \ldots, g_u \in \Z_n$ with $g_j = g^{\varepsilon_j \ell_j}$ ($\varepsilon_j = \pm 1$ and $0 \le \ell_j \le \lfloor \frac{n}{2} \rfloor$), we have $g_1 g_2 \cdots g_u = e$ if and only if $n \, | \sum_j \varepsilon_j \ell_j$. Given $\l = (\ell_1, \ell_2, \ldots, \ell_u)$, we will write in what follows $n | \l$ whenever there exists a family of signs $\underline{\varepsilon} = (\varepsilon_1, \varepsilon_2, \ldots, \varepsilon_u)$ such that $n$ divides $\sum_j \varepsilon_j \ell_j$. This motivates the following definition of admissible lengths $$\mathrm{Adm}_u(n) \, = \, \Big\{ \l \in L(u) \; : \; \lfloor \mbox{$\frac{n}{2}$} \rfloor \ge \ell_1 \ge \ldots \ge \ell_u \ge 1, \ \ell_1 \le \sum_{j=2}^u \ell_j, \ n | \l \Big\}.$$ We will also implicitly use the following consequences of our definition 
\begin{itemize}
\item If $n$ is even and $\l \in \mathrm{Adm}_u(n)$, then $|\l|$ is even.

\item If $n$ is odd, $\l \in \mathrm{Adm}_u(n)$ and $|\l| < n$, then $|\l|$ is even.
\end{itemize}

\subsubsection{Estimates for $s_3(r)$}

Our estimates for the sum $s_3(r)$ usually require some information on the structure of the metric space $(\G,\psi)$. This is also the case for finite cyclic groups with the word length, which require a different approach compared to our estimates for free and triangular groups. By Markovianity, $f$ is assumed to have symmetric positive Fourier coefficients. In particular, since $N_k \le 2$ we set $$a_k = \widehat{f}(g^k) = \widehat{f}(g^{-k}) \ge 0 \quad \mbox{for} \quad 0 \le k \le \big\lfloor \frac{n}{2} \big\rfloor.$$ Therefore, the $\alpha_k$'s are described for $0 \le k \le \lfloor \frac{n}{2} \rfloor$ by the formula $$\alpha_k = \delta_{k,0} a_0^2 + 2 \Big( \sum_{j=1}^{\lfloor \frac{n}{2} \rfloor -1} \delta_{k,j} a_j^2 \Big) + 2^{\delta_{n \, \mathrm{odd}}} \delta_{k, \lfloor \frac{n}{2} \rfloor} a_{\lfloor \frac{n}{2} \rfloor}^2.$$

\begin{proposition} \label{s23cyclic}
We have 
\begin{eqnarray*}
\widehat{f}(e) s_3(r) & \le & \sum_{\begin{subarray}{c} \lfloor \frac{n}{2} \rfloor \ge \ell_1 \ge \ell_2 \ge \ell_3 \ge 1 \\ \ell_1 + \ell_2 + \ell_3 = n \end{subarray}} M(\l) \big( \mbox{$\frac12$} \alpha_0 \alpha_{\ell_1} + \mbox{$\frac14$} \alpha_{\ell_2} \alpha_{\ell_3} \big) r^n \\ & + & \sum_{\begin{subarray}{c} \lfloor \frac{n}{2} \rfloor \ge \ell_1 \ge \ell_2 \ge \ell_3 \ge 1 \\ \ell_1 = \ell_2 + \ell_3 \end{subarray}} M(\l) \big( \mbox{$\frac12$} \alpha_0 \alpha_{\ell_1} + \mbox{$\frac14$} \alpha_{\ell_2} \alpha_{\ell_3} \big) r^{2 \ell_1}.
\end{eqnarray*}
\end{proposition}

\dem If we consider the sets $$\mathrm{D}_3(n) = \Big\{ \d = (d_1,d_2,d_3) \in \Z^3 \, : \, \lfloor \mbox{$\frac{n}{2}$} \rfloor \ge |d_1| \ge |d_2| \ge |d_3| \ge 1, \ n \hskip1pt | \hskip1pt (d_1 + d_2 + d_3) \Big\},$$ then we may write the sums $s_3(r)$ as follows 
\begin{eqnarray*}
\widehat{f}(e) s_3(r) & = & \sum_{\begin{subarray}{c} g_1g_2g_3 = e \\ g_j \neq e \end{subarray}} \widehat{f}(e) \prod_{j=1}^3
\widehat{f}(g_j) r^{|g_j|} \\ & = & \sum_{\d \in \mathrm{D}_3(n)} M(\d) a_0 \prod_{j=1}^3 a_{|d_j|} r^{|d_j|} \\ [8pt] & \le & \frac12 \sum_{\d \in \mathrm{D}_3(n)} \underbrace{M(\d) \big( a_0^2 a_{|d_1|}^2 + a_{|d_2|}^2 a_{|d_3|}^2 \big) r^{|\d|}}_{\gamma(\d)}.
\end{eqnarray*}
Case 1. If $n$ is odd, we decompose $\mathrm{D}_3(n)$ into the disjoint union $$\bigcup_{k=0}^3 \mathrm{D}_3(n,k) \quad \mbox{with} \quad \mathrm{D}_3(n,k) = \Big\{ \d \in \mathrm{D}_3(n) \; : \; \big| \big\{ j \ \mbox{ s.t. } \, d_j < 0 \big\} \big| = k \Big\}.$$ Then, it is clear that $\mathrm{D}_3(n,k) = - \mathrm{D}_3(n,3-k)$ for all $0 \le k \le 3$. On the other hand, since the terms $\gamma(\d)$ are invariant with respect to change of sign $\d \mapsto -\d$, we deduce the following estimate $$\widehat{f}(e) s_3(r) \, \le \, \sum_{\d \in \mathrm{D}_3(n,0)} \gamma(\d) + \sum_{\d \in \mathrm{D}_3(n,1)} \gamma(\d).$$ If $\d \in \mathrm{D}_3(n,0)$, we have $0 < d_j \le \frac{n}{2}$ and $n | \sum_j d_j$, so that $d_1 + d_2 + d_3 = n$. On the other hand, if $\d \in \mathrm{D}_3(n,1)$ we have $0 < |d_j| \le \frac{n}{2}$ and only one of the $d_j$'s is negative. This implies $|\sum_j d_j| < n$, so that the divisibility condition $n \, | \sum_j d_j$ forces $\sum_j d_j = 0$. Now, since $|d_1| \ge |d_2 | \ge |d_3|$ we must have $|d_1| = |d_2| + |d_3|$ which implies $|\d| = 2 |d_1|$. Altogether, we get 
\begin{eqnarray*}
\widehat{f}(e) s_3(r) & \le & \sum_{\begin{subarray}{c} \lfloor \frac{n}{2} \rfloor \ge \ell_1 \ge \ell_2 \ge \ell_3 \ge 1 \\ \ell_1 + \ell_2 + \ell_3 = n \end{subarray}} M(\l) \big( a_0^2 a_{\ell_1}^2 + a_{\ell_2}^2 a_{\ell_3}^2 \big) r^n \\ & + & \sum_{\begin{subarray}{c} \lfloor \frac{n}{2} \rfloor \ge \ell_1 \ge \ell_2 \ge \ell_3 \ge 1 \\ \ell_1 = \ell_2 + \ell_3 \end{subarray}} M(\l) \big( a_0^2 a_{\ell_1}^2 + a_{\ell_2}^2 a_{\ell_3}^2 \big) r^{2 \ell_1}.
\end{eqnarray*}
To be rigorous, we should write $M(-\ell_1, \ell_2, \ell_3)$ instead of $M(\l) = M(\ell_1, \ell_2, \ell_3)$ in the second sum. Nevertheless, since $\ell_1 = \ell_2 + \ell_3$ and the latter are positive we obtain $\ell_1 > \max \{\ell_2, \ell_3\}$ and $M(-\ell_1, \ell_2, \ell_3) = M(\l)$. On the other hand, since $n$ is odd and $\ell_j \ge 1$ we have $2a_{\ell_j}^2 = \alpha_j$ for $j=1,2,3$ and the assertion follows. 

\noindent Case 2. If $n$ is even, we first distinguish those products $g_1 g_2 g_3 = e$ with $|g_j| = \frac{n}{2}$ for some $j$. Observe that this can only occur for one $j$. Since there is just one element in $\Z_n$ of length $\frac{n}{2}$, this means we should replace the set $\mathrm{D}_3(n)$ above by $\mathrm{D}_3'(n) \cup \mathrm{D}_3''(n)$ with 
\begin{eqnarray*}
\mathrm{D}_3'(n) & = & \Big\{ \hskip1pt ( \hskip1pt \mbox{$\frac{n}{2}$} \hskip1pt , d_2, d_3) \in \Z^3 \; : \; \mbox{$\frac{n}{2}$} > |d_2| \ge |d_3| \ge 1, \ n | (\mbox{$\frac{n}{2}$} + d_2 + d_3) \Big\}, \\
\mathrm{D}_3''(n) & = & \Big\{ (d_1, d_2, d_3) \in \Z^3 \; : \; \mbox{$\frac{n}{2}$} > |d_1| \ge |d_2| \ge |d_3| \ge 1, \ n | (d_1 + d_2 + d_3) \Big\}.
\end{eqnarray*}
Arguing as above, this gives 
\begin{eqnarray*}
\widehat{f}(e) s_3(r) & \le & \frac12 \sum_{\d \in \mathrm{D}_3'(n)} M(\d) \big( a_0^2 a_{\frac{n}{2}}^2 + a_{|d_2|}^2 a_{|d_3|}^2 \big) r^{|\d|} \\ & + & \frac12 \sum_{\d \in \mathrm{D}_3''(n)} M(\d) \big( a_0^2 a_{|d_1|}^2 + a_{|d_2|}^2 a_{|d_3|}^2 \big) r^{|\d|} \ =: \ \mathrm{S}' + \mathrm{S}''.
\end{eqnarray*}
To estimate the sum $\mathrm{S}'$ we split $\mathrm{D}_3'(n) = \mathrm{D}_3'(n,+) \cup \mathrm{D}_3'(n,-)$ which recollect those $(\frac{n}{2},d_2,d_3) \in \mathrm{D}_3'(n)$ with $d_2,d_3$ both positive/negative respectively. Note that it can not happen that $\mathrm{sgn}(d_2) \neq \mathrm{sgn}(d_3)$. Arguing again as above we get $$\mathrm{S}' \le \sum_{\begin{subarray}{c} \frac{n}{2} > \ell_2 \ge \ell_3 \ge 1 \\ \frac{n}{2} + \ell_2 + \ell_3 = n \end{subarray}} M(\l) \big( \mbox{$\frac12$}\alpha_0 \alpha_{\frac{n}{2}} + \mbox{$\frac18$} \alpha_{\ell_2} \alpha_{\ell_3} \big) r^n + \sum_{\begin{subarray}{c} \frac{n}{2} > \ell_2 \ge \ell_3 \ge 1 \\ \frac{n}{2} = \ell_2 + \ell_3 \end{subarray}} M(\l) \big( \mbox{$\frac12$}\alpha_0 \alpha_{\frac{n}{2}} + \mbox{$\frac18$} \alpha_{\ell_2} \alpha_{\ell_3} \big) r^n$$ where we have used crucially that $|d_j| < \frac{n}{2}$ ($j=1,2$) for elements in $\mathrm{D}_3'(n)$. On the other hand, since the same property holds for $j=1,2,3$ in $\mathrm{D}_3''(n)$, we still have $2a_{|d_j|}^2 = \alpha_{|d_j|}$ and the argument for the odd case yields
\begin{eqnarray*}
\mathrm{S}'' & \le & \sum_{\begin{subarray}{c} \lfloor \frac{n}{2} \rfloor > \ell_1 \ge \ell_2 \ge \ell_3 \ge 1 \\ \ell_1 + \ell_2 + \ell_3 = n \end{subarray}} M(\l) \big( \mbox{$\frac12$} \alpha_0 \alpha_{\ell_1} + \mbox{$\frac14$} \alpha_{\ell_2} \alpha_{\ell_3} \big) r^n \\ & + & \sum_{\begin{subarray}{c} \lfloor \frac{n}{2} \rfloor > \ell_1 \ge \ell_2 \ge \ell_3 \ge 1 \\ \ell_1 = \ell_2 + \ell_3 \end{subarray}} M(\l) \big( \mbox{$\frac12$} \alpha_0 \alpha_{\ell_1} + \mbox{$\frac14$} \alpha_{\ell_2} \alpha_{\ell_3} \big) r^{2 \ell_1}.
\end{eqnarray*}
Summing our estimates for the sums $\mathrm{S}'$ and $\mathrm{S}''$ we obtain the assertion. \fin 

\subsubsection{Numerical estimates for $\Z_n$}

If we keep the notation introduced in Paragraph \ref{numericalF2}, then \eqref{Reg-Pat} still holds true and we will see again that the super-pathological sums do not affect the regular terms. We collect below the formulas for the left coefficients. We refer to Appendix C for the details of the proofs. The critical function for $\Z_n$ with the word length is given by the graph

\begin{tikzpicture}

\draw[->] (1.5,0) -- (10,0) node[anchor=north] {$m$};
\draw	(1.5,0) node[anchor=north] {0}
		(2.5,0) node[anchor=north] {1}
		(3.5,0) node[anchor=north] {2}
		(4.5,0) node[anchor=north] {3}
		(5.5,0) node[anchor=north] {4}
		(6.5,0) node[anchor=north] {5}
		(7.5,0) node[anchor=north] {6}
		(8.5,0) node[anchor=north] {7};
		
\draw	(5.5,-1) node{{\scriptsize {\bf The critical function} $\mu_q(\Z_n, |\cdot|, m)$}};

\draw[->] (1.5,0) -- (1.5,4.5) node[anchor=east] {$\mu_q(m)$};

\draw 	(1.5,1) node[anchor=east] {1}
		(1.5,2) node[anchor=east] {2}
		(1.5,3) node[anchor=east] {3};
		
\draw[circle,fill=black] (-0.1,0);	

\filldraw [fill=black] (2.5,0) circle (2pt);
\filldraw [fill=black] (3.5,3) circle (2pt);
\filldraw [fill=black] (4.5,3) circle (2pt);
\filldraw [fill=black] (5.5,2) circle (2pt);
\filldraw [fill=black] (6.5,2) circle (2pt);
\filldraw [fill=black] (7.5,2) circle (2pt);
\filldraw [fill=black] (8.5,1) circle (2pt);
    	
\draw[dotted] (2.5,0) -- (3.5,3);
\draw[dotted] (3.5,3) -- (4.5,3);
\draw[dotted] (4.5,3) -- (5.5,2);
\draw[dotted] (5.5,2) -- (6.5,2);
\draw[dotted] (6.5,2) -- (7.5,2);
\draw[dotted] (7.5,2) -- (8.5,1);
\draw[->][dotted] (8.5,1) -- (10,1);

\end{tikzpicture}

\noindent ${\bf (\alpha)}$ {\bf Coefficients for $u=2,3$.} 
If $n \ge 6$, we find
\begin{itemize}
\item[$\alpha$i)] $C_{q,s_{2}}^{\mathrm{left}}[(k,\underline{0})](r) \, = \, r^{2k}$.

\vskip5pt

\item[$\alpha$ii)] $\displaystyle C_{q,s_{3}}^{\mathrm{left}}[(k_1,k_2, \underline{0})](r) \, = \, \left\{\begin{array}{ll}
\frac32 r^4 &\mbox{if} \ (k_1,k_2) = (1,1), \\ [3pt] 3 r^6 & \mbox{if} \ (k_1,k_2) = (2,1), \\  [3pt] 3 r^8 & \mbox{if} \ (k_1,k_2) = (3,1). \end{array}\right.$

\vskip5pt

\item[$\alpha$iii)] $\displaystyle C_{q,s_{3}}^{\mathrm{left}}[(k_1,k_2, \underline{0})](r) \, = \, \left\{\begin{array}{ll}
\frac14 r^6 &\mbox{if} \ (k_1,k_2) = (2,2) \mbox{ and } n=6, \\ [3pt] \frac34 r^7 & \mbox{if} \ (k_1,k_2) = (2,2) \mbox{ and } n = 7, \\  [3pt] \frac32 r^8 & \mbox{if} \ (k_1,k_2) = (2,2) \mbox{ and } n \ge 8. \end{array}\right.$

\vskip5pt

\item[$\alpha$iv)] When $k_1 + k_2 \ge 5$, we obtain $$C_{q,s_{3}}^{\mathrm{left}}[(k_1,k_2, \underline{0})](r) \, = \, 
\frac{3}{s(s-1)} M(k_1,k_2, \underline{0}) r^{k_1+k_2+2}.$$

\vskip5pt

\item[$\alpha$v)] When $1 \le k \le \lfloor \frac{n}{2} \rfloor$, we obtain $$C_{q,s_{3}}^{\mathrm{left}}[(k, \underline{0})](r) \, = \, \left\{\begin{array}{ll}
0 &\mbox{if} \ k = 1, \\ [3pt] \frac32 r^4 & \mbox{if} \ k = 2, \\  [3pt] 3(k+1) r^{2k} & \mbox{if} \ 3 \le k \le \lfloor \frac{n}{2} \rfloor. \end{array}\right.$$
\end{itemize} 
\vskip5pt

\noindent ${\bf (\beta)}$ {\bf Coefficients for $\Lambda$-regular terms.} 
No $\Lambda$-regular terms since $q(\Z_n,|\cdot|)=2$. 

\vskip5pt

\noindent ${\bf (\gamma)}$ {\bf Coefficients for $\Delta$-regular terms.} If $\k \in \mathrm{Reg}(s)$, we use $\Delta$-estimates. Given $s \ge 2$ and $0 \le r < 1$, we may rewrite estimates $\gamma$i), $\gamma$ii) and $\gamma$iii) for the free group replacing the sum $\mathcal{G}(\G,\psi,r) = 2nr/(1-(2n-1)r)$ there by our estimate $2r/(1-r)$.

\vskip5pt

\noindent ${\bf (\delta)}$ {\bf Coefficients for pathological terms.} The coefficients $C_{q,s_u}^{\mathrm{left},\mathbf{P}}[\k](r)$ follow from the formulas in Paragraph \ref{numericalF2} for $\F_n$ deduced from the $\Lambda$-estimates, with $N_k\leq 2$ in the present case. Again, the $N_k$'s are not increasing (a property which is applied in the proof, see Appendix C) but our upper bound $N_k \le 2$ is monotone and this suffices. 

\vskip5pt

\noindent ${\bf (\varepsilon)}$ {\bf Coefficients for super-pathological terms.} As usual, computer assistance allows us to identify the super-pathological terms. These are of the form $(1,1,\underline{0})$ adding at the end as many zeros as needed to form $s$-tuples. The exceptional sets $B_u$ are given by 
$$B_4  =  \big\{ (1,1,1,1),(3,1,1,1) \big\}.$$
We have $B_u=\emptyset$ otherwise. Given $0\leq r \leq 1$, we find
\begin{itemize}
\item[$\varepsilon$i)] $s_4[(1,1,1,1)](r) = \frac32 r^4\alpha_1^2$,

\item[$\varepsilon$ii)] $s_4[(3,1,1,1)](r) \le r^6 (\alpha_1^2 + \alpha_1 \alpha_3)$ for $n \ge 7$,

\item[$\varepsilon$iii)] $s_4[(3,1,1,1)](r) \le r^6 (\alpha_1^2 + 2\alpha_1 \alpha_3)$ for $n = 6$.
\end{itemize}

\begin{remark}\label{RkCleftReg-cyclic}
\emph{Note that $\k \in L_2(s)$ for the $\alpha_{\k}$'s in the right hand side of the estimates above. In particular, regularity means $|\k| \ge 2 + \mu_q(\Z_n,|\cdot|,2) = 5$. Since $|\k| \le 4$ for the $\alpha_{\k}$'s appearing above, the super-pathological sums do not affect the regular terms $\k \in \mathrm{Reg}(s)$. It implies that for regular terms we have $C_{q,s_{u}}^{\mathrm{left},\mathbf{S}}[\k](r)=0$.}
\end{remark}

\subsubsection{Proof of Theorem \emph{A3}}

Again, the goal is to show $$C_q^{\mathrm{left}}[\k](r) = \sum_{u = 2 \vee (2m-1)}^{q \wedge (2m+1)} \binom{q}{u} C_{q,s_u}^{\mathrm{left}}[\k](r) \, \le \, C_q^{\mathrm{right}}[\k]=M(\k).$$ 

\noindent {\bf The regular case $\k \in \mathrm{Reg}_m(s)$.} By Remark \ref{RkCleftReg-cyclic}, we will only need $(\alpha)$ and $(\gamma)$ to fix the regular case. When $m=1$ and $\k =(k, \underline{0}) \in \mathrm{Reg}_1(s)=L_1(s)$ we need to prove the inequalities \eqref{Reg1-1-cyclic}-\eqref{Reg1-2-cyclic} in Appendix D. As usual, the case $k=1$ is trivial. When $m=2$, we have $\k=(k_1,k_2,\underline{0}) \in \mathrm{Reg}_2(s)$ if $|\k| \geq 2+\mu_q(\Z_n,|\cdot|,2)=5$ and we use $\alpha$iv), $\gamma$i) and $\gamma$iii). The inequality is \eqref{Reg2-cyclic} in Appendix D. When $3\leq m \leq s$ and $\k \in \mathrm{Reg}_m(s)$, we reduce the desired inequality to \eqref{Regm-cyclic} in Appendix D by means of $\gamma$i), $\gamma$ii) and $\gamma$iii). 

\noindent {\bf The pathological case $\k \in \mathrm{Pat}_m(s)$.} We proceed as for $\F_2$. We first fix by hand the terms of the form $\k=(1,1,\ldots,1,\underline{0}) \in \mathrm{Pat}_m(s)$ for $m\geq 7$, then we use computer assistance to treat the finitely many remaining cases. According to our definition of admissible lengths $\mathrm{Adm}_u(n)$ in $\Z_n$, we may define
the corresponding sets $A'_{2m}(n)[\k], A_{2m}(n)[\k]$ and $A'_{2m+2}(n)[\k]$ for the $(\delta)$-estimates accordingly. The description of these sets depends on the parity of $n$. We include for clarity the concrete form of these sets below 
\begin{itemize}
\item If $\max\{7,\frac{s}{2}\} \le m \le s$
\begin{eqnarray*}
\hskip40pt A'_{2m}(n)[\k] \!\!\!\! & = & \!\!\!\! \big\{ (\ell,1,1,\ldots,1,0) \, : \, 2 \leq \ell \leq 2m-2 \big\} \subset \mathrm{Adm}'_{2m}(n), \\ A_{2m}(n)[\k] \!\!\!\! & = & \!\!\!\! \big\{ (\ell,1,1,\ldots,1,1) \, : \, 1 \leq \ell \leq 2m-1 \big\} \subset \mathrm{Adm}_{2m}(n), \\ A'_{2m+2}(n)[\k] \!\!\!\! & = & \!\!\!\! \big\{ (\ell,1,1,\ldots,1,0) \, : \, 2 \leq \ell \leq 2m \big\} \hskip4pt \subset \hskip4pt  \mathrm{Adm}'_{2m+2}(n).
\end{eqnarray*}

\item If $7 \le m < \frac{s}{2}$
\begin{eqnarray*}
\hskip40pt A'_{2m}(n)[\k] \!\!\!\! & = & \!\!\!\! \big\{ (\ell,1,1,\ldots,1,0) \, : \, 2 \leq \ell \leq 2m-2, \ \ell \ \mbox{even} \big\} \subset \mathrm{Adm}'_{2m}(n), \\ A_{2m}(n)[\k] \!\!\!\! & = & \!\!\!\! \big\{ (\ell,1,1,\ldots,1,1) \, : \, 1 \leq \ell \leq 2m-1, \ \ell \ \mbox{odd} \big\} \hskip1pt \subset \hskip2pt \mathrm{Adm}_{2m}(n), \\ A'_{2m+2}(n)[\k] \!\!\!\! & = & \!\!\!\! \big\{ (\ell,1,1,\ldots,1,0) \, : \, 2 \leq \ell \leq 2m, \ \ell \ \mbox{even} \big\} \hskip4pt \subset \hskip4pt  \mathrm{Adm}'_{2m+2}(n).
\end{eqnarray*}
\end{itemize}
Namely, if $\max\{7,\frac{s}{2}\} \le m \le s$ we just need to justify why we have excluded the case $\ell = 1$ from $A_{2m}'(n)[\k]$ and $A_{2m+2}'(n)[\k]$. Indeed, recall that $|\l|$ is odd in this case. If $n$ is even, this breaks the parity condition $|\l| \in 2 \Z$. If $n$ is odd, our assumption $q \le n$ implies crucially that $|\l| < n$. Since $|\l|$ is odd, $\l$ can not belong to $\mathrm{Adm}_{2m}'(n)$ as we noticed in Paragraph \ref{Adm-Cyclic}. Using $(\delta)$, this yields for $0 \le r < 1$
\begin{eqnarray*}
\lefteqn{C_q^{\mathrm{left}}[\k](r) \, = \, \binom{2s}{2m} \frac{m}{2} (2r^2)^m \sum_{\ell=1}^{2m-1} r^{\ell-1}} \\ \!\!\!\! & + & \!\!\!\! \binom{2s}{2m-1} \frac{2m-1}{16} (2r^2)^m \sum_{\ell=2}^{2m-2} r^{\ell-2} \, + \, \binom{2s}{2m+1} \frac{2m+1}{4} (2r^2)^m \sum_{\ell=2}^{2m} r^{\ell} \\ [5pt] \!\!\!\! & \leq & \!\!\!\! \Big[ \binom{2s}{2m-1} \frac{2m-1}{16} + \binom{2s}{2m} \frac{m}{2} + \binom{2s}{2m+1} \frac{2m+1}{4}r^2 \Big] \frac{(2r^2)^{m}}{1-r}. 
\end{eqnarray*}
Since $C_q^{\mathrm{right}}[\k]=\binom{s}{m}$, we are reduced to prove inequality \eqref{Pat1-cyclic} in Appendix D. In the second case $7 \le m < \frac{s}{2}$, the novelty comes from the parity condition $|\l| \in 2 \Z$ for the admissible lengths $\l$ appearing in $A_{2m}(n)[\k]$, $A_{2m}'(n)[\k]$ and $A_{2m+2}'(n)[\k]$ above. As explained in Paragraph \ref{Adm-Cyclic}, this is always the case when $n$ is even. If $n$ is odd we additionally need that $|\l| < n$. Note however that in this situation we have $$|\l| \le \ell + 2m \le 4m < 2s = q \le n.$$ Using $(\delta)$, this yields for $0 \le r < 1$
\begin{eqnarray*}
\lefteqn{C_q^{\mathrm{left}}[\k](r) \, = \, \binom{2s}{2m} \frac{m}{2} (2r^2)^m \sum_{\begin{subarray}{c} \ell=1 \\ \ell \ \mathrm{odd} \end{subarray}}^{2m-1} r^{\ell-1}} \\ \!\!\!\! & + & \!\!\!\! \binom{2s}{2m-1} \frac{2m-1}{16} (2r^2)^m \sum_{\begin{subarray}{c} \ell=2 \\ \ell \ \mathrm{even} \end{subarray}}^{2m-2} r^{\ell-2} \, + \, \binom{2s}{2m+1} \frac{2m+1}{4} (2r^2)^m \sum_{\begin{subarray}{c} \ell=2 \\ \ell \ \mathrm{even} \end{subarray}}^{2m} r^{\ell} \\ [5pt] \!\!\!\! & \leq & \!\!\!\! \Big[ \binom{2s}{2m-1} \frac{2m-1}{16} + \binom{2s}{2m} \frac{m}{2} + \binom{2s}{2m+1} \frac{2m+1}{4}r^2 \Big] \frac{(2r^2)^{m}}{1-r^2}. 
\end{eqnarray*}
Since $C_q^{\mathrm{right}}[\k]=\binom{s}{m}$, we are reduced to prove an inequality which is clearly weaker than \eqref{Pat1-cyclic} in Appendix D. We conclude as in the free group case, by implementing $(\alpha),(\delta)$ and $(\varepsilon)$ in a computer to fix the terms in $\bigcup_{2\leq m \leq 6} \mathrm{Pat}_m(s)$. All these terms verify the test. \fin 

\subsection{Comments} \label{RFnCritical}

We finish this section with some general comments:

\subsubsection{On the growth condition}

As explained before, the gap $4 \le q \le q(\G,\psi)$ between $4$ and the critical index usually requires different and sometimes ad hoc estimates to complete the argument. It is worth mentioning that $R(\G,\psi) = 1$ for pairs $(\G,\psi)$ of polynomial ---or even subexponential--- growth, so that $q(\G,\psi) = 2$ for this class of groups/lengths and we find no gap. 

\subsubsection{On the cancellation condition}

Any $(\G,\psi)$ satisfying our growth condition attains its spectral gap. That is, the set of elements $g \in \G$ with $\psi(g) = \sigma$ is nonempty. Given $m \ge 3$, an $m$-loop associated to the pair $(\G,\psi)$ is any relation of the form $g_1 g_2 \cdots g_m = e$ with $\psi(g_j) = \sigma$. The relations $g^2=e$ are not considered loops, so that $3$-loops are the smallest possible ones. Let us first show that the presence of $3$-loops in $(\G,\psi)$ yields nonstandard optimal hypercontractivity bounds. Namely, recall that we set $T(p,q,\sigma)$ for the \lq\lq expected" optimal time $\frac{1}{2\sigma} \log (\frac{q-1}{p-1})$ and $t_{p,q} = t_{p,q}(\G,\psi)$ for the optimal time. We always have $t_{p,q} \ge T(p,q,\sigma)$. Then the following result holds when $(\G,\psi)$ admits $3$-loops
\begin{equation} \label{Noloop3}
t_{p,q} = T(p,q,\sigma) \ \Rightarrow \ (p-1)^{\frac32}(q-2) \le (q-1)^{3/2} (p-2).
\end{equation}
This implies $t_{p,q} > T(p,q,\sigma)$ in the presence of $3$-loops for $p \le 2$ and $q > p$, it also yields nonstandard optimal times for $2 \le p \le 4$ and small values of $q>p$. The proof is not difficult. Indeed, take $g_1, g_2, g_3$ with $\psi(g_j) = \sigma$ such that $g_1g_2g_3=e$. Define $\Lambda = \sum_{j=1}^3 \lambda(g_j) + \lambda(g_j^{-1})$. Note that we do not assume $g_j \neq g_j^{-1}$ or $g_j \neq g_k$ ($j \neq k$) so that there might be repetitions in the sum defining $\Lambda$. Let us set $t = T(p,q,\sigma)$ and define $f = \mathbf{1} + \varepsilon \Lambda$ so that $\S_{\psi,t}f = \mathbf{1} + \varepsilon \mu \Lambda$ with $\mu = \exp(-t\sigma) = \sqrt{(p-1)/(q-1)}$. Since $\Lambda$ is self-adjoint and both $f$ and $\S_{\psi,t}f$ live in the unital $*$-algebra generated by $\Lambda$, we may assume in what follows that $\G$ is abelian since all our calculations will be made in this commutative algebra. In particular we may approximate $|\S_{\psi,t}f|^q$ by its Taylor expansion up to degree $2$ 
\begin{eqnarray*}
|\S_{\psi,t}f|^q & = & \big( \mathbf{1} + 2\varepsilon \mu \Lambda + \varepsilon^2 \mu^2 \Lambda^2 \big)^{\frac{q}{2}} \\ & \sim & \mathbf{1} + \varepsilon q \mu \Lambda + \varepsilon^2 \mu^2 \frac{q^2-q}{2} \Lambda^2 + \varepsilon^3 \mu^3 \frac{q^2-2q}{2} \Lambda^3 + \varepsilon^4 \mu^4 \frac{q^2-2q}{8} \Lambda^4. 
\end{eqnarray*}
Using the same formula, $\tau(\Lambda)=0$ and another Taylor expansion of order $1$, we get
\begin{eqnarray*}
\big( \tau |f|^p \big)^{\frac{q}{p}} & \sim & 1 + \varepsilon^2 \frac{q(p-1)}{2} \tau(\Lambda^2) + \varepsilon^3 \frac{q(p-2)}{2} \tau(\Lambda^3) + \varepsilon^4 \frac{q(p-2)}{8} \tau(\Lambda^4). 
\end{eqnarray*}
Assume now that $t_{p,q} = T(p,q,\sigma)$, then we must have 
$$\lim_{\varepsilon \to 0} \frac{\|f\|_p^q - \|\S_{\psi,t}f\|_q^q}{\varepsilon^3} = \frac{q}{2} \Big[ (p-2) - \Big( \frac{p-1}{q-1} \Big)^{\frac32} (q-2) \Big] \tau(\Lambda^3) \ge 0.$$ The presence of the loop $g_1g_2g_3=e$ easily implies $\tau(\Lambda^3) > 0$ and the result follows.

The simplest pair $(\G,\psi)$ admitting a $3$-loop is $\Z_3$ with the word length. Optimal logarithmic Sobolev inequalities and related estimates were obtained by Andersson and Diaconis/Saloff-Coste \cite{A1,DS}. Unfortunately, Gross equivalence between log Sobolev estimates and hypercontractivity does not give optimal hypercontractivity bounds for $\Z_3$, which are still open.  

On the other hand, in this paper we avoid small loops. As explained in the Introduction and illustrated for triangular groups above, only our estimate for the sum $s_3(r)$ requires this additional assumption. Observe that our estimate of $s_3(r)$ for $\Delta_{\alpha \beta \gamma}$ can be trivially extended to any other pair $(\G,\psi)$ not admitting small enough loops, we just need to use the decay properties of $\S_{\psi}$. On the other hand, as it was illustrated for finite cyclic groups, we may not include a parity condition for admissible lengths in the presence of odd loops. As in that case, this imposes a restriction on the set of $q \in 2 \Z_+$ which we can consider. Namely, using that we have $\mu_q(\G,\psi,m) \ge 1$ for all $m \ge 2$ and $q$ large, we conclude that $\k = (1,1,\ldots,1,\underline{0}) \in \mathrm{Pat}_m(s)$ for all $m$. In presence of an odd loop of size $2m+1$, we may consider $\l = (1,1,\ldots,1) \in \mathrm{Adm}_{2m+1}$. Completing squares, it is clear that $s_{2m+1}[\l](r)$ necessarily has a non-zero contribution to the term $\alpha_{\k}$ for $\k = (1,1,\ldots,1,\underline{0}) \in \mathrm{Pat}_m(s)$ with a coefficient $c(m) \cdot r^{2m+1}$. This means for $r = 1/ \sqrt{2s-1}$ that $$\frac{C_q^{\mathrm{left}}[\k](r)}{C_q^{\mathrm{right}}[\k]} \ge \frac{\binom{2s}{2m+1}}{\binom{s}{m}} \frac{c(m)}{(2s-1)^{\frac{2m+1}{2}}} = f(m,s)$$ and it is easily checked that $f(m,s)$ behaves as $\sqrt{s}$ for $m$ fixed and $s \to \infty$. 

\subsubsection{Beyond conditionally negative lengths}

According to Schoenberg's theorem $(\S_{\psi,t})_{t \ge 0}$ is Markovian iff the length $\psi$ is conditionally negative. This property is crucial in Appendix A to adapt Gross extrapolation method to our setting and extend our optimal time $L_2 \to L_q$ estimates for $q \in 2 \Z_+$ to other more general indices $1 < p \le q < \infty$. On the contrary, our combinatorial method for $q$ even seems unaffected if we remove this assumption. Indeed, we have only used it to reduce the problem to positive $f$ and subsequently to $f$ with symmetric positive Fourier coefficients. It is however relatively simple to adapt our arguments to a general $f \in L_2(\L(\G))$. Indeed, we clearly have $$\|\T_{\psi,r}f\|_q^q \le \sum_{g_1 g_2 \cdots g_q=e} \prod_{j=1}^s |\widehat{f}(g_{2j-1}^{-1})| \, |\widehat{f}(g_{2j})| r^{\psi(g_{2j-1}) + \psi(g_{2j})}.$$ This means that we should invert the $g_j$'s labelled by an odd $j$, which of course affects our expressions for the sums $s_u(r)$. Nevertheless, a careful reading of Section \ref{Method} will convince the reader that adapting our arguments to this more general framework only requires to introduce a little more complicated terminology to keep track of the powers of the $g_j$'s. Since we have not included an illustration of our method in the non-Markovian setting, we have decided to avoid this more general formulation for clarity in the exposition. 

\subsubsection{Hypercontractivity in $\F_\infty$}

Theorem A1 does not give any result for $\F_\infty$ since $q(n) \to \infty$ with $n$. However, we may provide optimal hypercontractivity bounds in $\F_\infty$ with respect to a weighted form of the word length. More precisely, let $c_1, c_2, \ldots$ denote the free generators of $\F_\infty$, consider a sequence $\mathrm{m}=(m_k)_{k\geq 1}$ of positive integers and define the length 
\begin{equation} \label{WeightedLength}
\psi_\mathrm{m}(\underbrace{c_{j_1}^{s_1}c_{j_2}^{s_2}\cdots c_{j_\ell}^{s_\ell}}_w) \, = \, \sum_{k=1}^\infty (2m_k + 1) \sum_{n: \, j_n = k} |s_n|
\end{equation}
for $w$ written in reduced form. It can be checked that $\psi_\mathrm{m}$ is conditionally negative for any choice of $m_k \in \R_+$. Indeed, it follows easily adapting Haagerup's original argument in \cite{Ha}. Note that we recover the word length when $m_k=0$ for all $k \ge 1$. Optimal hypercontractivity follows when $\mathrm{m}$ is strictly increasing $0 \leq m_1 < m_2 < \ldots$ Namely, consider the injective homomorphism $J_{\mathrm{m}}: c_j \in \F_\infty\mapsto  b^{-m_j}ab^{m_j} \in \F_2$ where $a,b$ denote the free generators of $\F_2$. This map clearly leads to an isometry $L_p(\mathcal{L}(\F_\infty)) \to L_p(\mathcal{L}(\F_2))$ for all $p\geq 1$, still denoted by $J_{\mathrm{m}}$. On the other hand, we have $| J_{\mathrm{m}}(w)| \leq \psi_{\mathrm{m}}(w)$ for any $w \in \F_\infty$. Hence, if $f \in L_2(\mathcal{L}(\F_\infty))$ and $\P_t$ denotes the free Poisson semigroup in $\L(\F_2)$, Theorem A1 i) implies 
\begin{eqnarray*}
\|\S_{\psi_{\mathrm{m}},t}f \|_{ L_2(\mathcal{L}(\F_\infty))} & = & \Big( \sum_{w \in \F_\infty} e^{-2t \psi_{\mathrm{m}}(w)} |\widehat{f}(w)|^2 \Big)^{\frac12} \\ & \le & \Big( \sum_{w \in \F_\infty} e^{-2t |J_{\mathrm{m}}(w)|} |\widehat{f}(w)|^2 \Big)^{\frac12} \\ & = & \|\P_t(J_{\mathrm{m}}(f))\|_{ L_2(\mathcal{L}(\F_2))} \ \le \|J_{\mathrm{m}}(f)\|_{ L_p(\mathcal{L}(\F_2))} = \|f\|_{L_p(\mathcal{L}(\F_\infty))}
\end{eqnarray*}
for those $1 < p \le 2$ whose conjugate index is in $2\Z_+$ and $t \ge -\frac{1}{2}\log(p-1)$. 
In conclusion, we get optimal hypercontractive $L_2 \to L_q$ estimates for the Poisson-like semigroup associated to $(\F_\infty, \psi_{\mathrm{m}})$ and $q$ an even integer. As usual, for other values of $(p,q)$ we find an additional $\log 3$ in our bounds. 

\subsubsection{Asymmetric square completions}

As mentioned in Paragraph \ref{LambdavsDelta}, the critical function we would obtain by using $\Lambda$-estimates is not uniformly bounded, which explains the necessity of $\Delta$-estimates to treat the regular terms. Let us justify this claim in $\F_n$, by showing that if $\tilde{\mu}_q(m)$ denotes the critical function for $(\F_n,|\cdot|)$ obtained from $\Lambda$-estimates, then we get $\tilde{\mu}_q(\frac{q}{2}) \to \infty$ as $q$ tends to infinity. Indeed, by using Proposition \ref{Lambda-est} i) and estimating as we do in the proof of $(\beta)$, we get for $\k \in L_s(s)$ and $0\leq r < \frac{1}{2n-1}$ $$ C_{q,s_{q}}^{\mathrm{left}}[\k](r) =\frac{(2s)!}{s!}M(\k)\Big(\frac{2nr}{1-(2n-1)r}\Big)^{s-1}\frac{r^{|\k|+1}}{1-r}.$$ Then, for $|\k| \ge s + \tilde{\mu}_q(\frac{q}{2})$ we should find $$C_{q,s_{q}}^{\mathrm{left}}[\k](r)  \leq C_{q}^{\mathrm{left}}[\k](r) \leq C_{q}^{\mathrm{right}}[\k]=M(\k)$$ 
for $0\leq r \leq \frac{1}{\sqrt{2s-1}}$. This means that $$\tilde{\mu}_q(s)\geq \frac{2}{\log(2s-1)}\log\left[\frac{(2s)!}{s!}\Big(\frac{2n}{\sqrt{2s-1}-2n+1}\Big)^{s-1}\frac{1}{\sqrt{2s-1}-1}\right]-s=:X(s).$$ Stirling's formula gives that $X(s) \to \infty$ as $s\to \infty$. This justifies our claim.

\section{{\bf Poisson-like lengths}} \label{SectThB}

In this section we prove the hypercontractivity inequality in Theorem B for Poisson-like lengths and give some bounds for the constant $\beta(\G,\psi)$. We also prove Theorem C and give examples of Poisson-like lengths for which Theorems B and C produce new estimates. Let us first adapt the notation used so far to the present context. Given $\G$ discrete equipped with a Poisson-like length $\psi: \G \to \R_+$, we know from the exponential growth condition that there must exists a sequence $0 = n_0 < n_1 < n_2 < \cdots$ such that $\psi(\G) = \{n_k \, : \, k \ge 0\}$. Also, the spectral gap must be attained, so that we have $\sigma = n_1$. Let us write $W_{n_k} = \{g \in \G \, : \, \psi(g) = n_k\}$ and $N_{n_k} = |W_{n_k}|$. Enumerate $W_{n_k}$ by $w_{n_k}(1), w_{n_k}(2), \ldots, w_{n_k}(N_{n_k})$. Then we set for any $f \in L_2(\L(\G))$ $$a_0 = |\widehat{f}(e)| \quad \mbox{and} \quad a_{n_k}(i)=|\widehat{f}(w_{n_k}(i))| \quad \mbox{for} \quad 1\le i \le N_{n_k}.$$
Given $k \ge 1$, we define $\alpha_0 = |\widehat{f}(e)|^2 = a_0^2$ and $\displaystyle \alpha_{n_k} = \sum_{g \in W_{n_k}} |\widehat{f}(g)|^2 = \sum_{i=1}^{N_{n_k}} a_{n_k}(i)^2$.

\subsection{Proof of Theorem B}

By interpolation and using Gross' argument recalled in Appendix A, it suffices to consider the case $(p,q) = (2,4)$. More precisely, if we set $r = e^{-t}$ and $\T_{\psi,r} = \S_{\psi,t}$ accordingly, we need to find some $0 < \mathcal{R}_{\G,\psi} \le 3^{-1/2\sigma}$ such that 
\begin{equation}\label{aim-ThmB}
\|\T_{\psi,r}f\|_4 \le \|f\|_2 \hskip5pt \mbox{for any} \hskip5pt f \in L_2(\mathcal{L}(\G)) \hskip5pt \mbox{and} \hskip5pt 0 \le r \le \mathcal{R}_{\G,\psi}.
\end{equation} 
In fact, we just need to prove it for $f \ge 0$ but we will prove this particular estimate for arbitrary $f$ without assuming conditional negativity (which is only required to apply Gross extrapolation argument), which will also serve to prove the last assertion of Theorem B. Once this estimate is known, we may deduce the required result with a factor $\log 3$. More precisely, if we have $$\mathcal{R}_{\G,\psi} = \Big(\frac{1}{\sqrt{3}}\Big)^{\frac{\gamma(\G,\psi)}{\sigma}} \ \mbox{for some $\gamma(\G,\psi) \ge 1$} \ \Rightarrow \ \beta(\G,\psi) = \log 3 \, \gamma(\G,\psi).$$ In order to prove \eqref{aim-ThmB}, fix $f \in L_2(\mathcal{L}(\G))$ and observe that $$\|f\|_2^4 \, = \, \Big( \sum_{g \in \G} |\widehat{f}(g)|^2 \Big)^{2} \, = \, \Big(\sum_{k \ge 0} \alpha_{n_k} \Big)^{2} \, = \, \sum_{k \ge 0}  \alpha_{{n_k}}^2 + 2 \sum_{j<k} \alpha_{n_j} \alpha_{n_k}.$$ On the other hand, the left hand side of \eqref{aim-ThmB} is easily expressed as 
\begin{eqnarray*}
\lefteqn{\big\| \T_{\psi,r}f \big\|_4^4} \\ [5pt] & = & \sum_{g_1^{-1}g_2g_3^{-1}g_4=e}^{\null}  \overline{\widehat{f}(g_1)} \widehat{f}(g_2) \overline{\widehat{f}(g_3)}\widehat{f}(g_4) r^{\psi(g_1)+\psi(g_2)+\psi(g_3)+\psi(g_4)} \\ & = & \sum_{h\in \G} \Big[ \underbrace{\sum_{g_1^{-1}g_2=h} \overline{\widehat{f}(g_1)}\widehat{f}(g_2) r^{\psi(g_1)+\psi(g_2)}}_{\phi_\G(h)} \Big] \Big[ \underbrace{\sum_{g_3^{-1}g_4=h^{-1}} \overline{\widehat{f}(g_3)}\widehat{f}(g_4) r^{\psi(g_3)+\psi(g_4)}}_{\phi_\G(h^{-1})} \Big].
\end{eqnarray*}
Since $\phi_\G(h^{-1}) = \overline{\phi_\G(h)}$, it suffices to prove 
\begin{equation} \label{Aim2-4}
 \sum_{h\in \G} |\phi_\G(h)|^2 \, \le \, \sum_{k \geq 0}  \alpha_{n_k}^2 + 2 \sum_{j<k}  \alpha_{n_j} \alpha_{n_k}. 
\end{equation}
We have $\phi_\G(e) = \sum_{k \ge 0} \alpha_{n_k} r^{2n_k}$, while for $h\neq e$ 
\begin{eqnarray*}
|\phi_\G(h)| & = & \Big| \sum_{g \in \G}^{\null}  \overline{\widehat{f}(g)} \widehat{f}(gh) r^{\psi(g) + \psi(gh)} \Big| \\ & \le & \sum_{\ell \ge 0} \sum_{g \in W_{n_\ell}} \big| \widehat{f}(g) \widehat{f}(gh) \big| r^{\psi(g) + \psi(gh)} \\ & \le & \underbrace{\sum_{\ell \ge 0} \sum_{\substack{\ell' \geq 0\\ n_{\ell'} \in X(\psi(h),\ell)}} \sum_{(i_1,i_2) \in Y_{\ell'}^\ell(h)}}_{\displaystyle \sum_{\delta \in \Delta(h)}} \underbrace{a_{n_\ell}(i_1) a_{n_{\ell'}}(i_2)}_{\gamma_\delta(h)} \underbrace{r^{n_\ell+n_{\ell'}}}_{\nu_\delta(h,r)} ,
\end{eqnarray*}
where $X(\psi(h),\ell) = \psi(\G) \cap [|\psi(h)-n_\ell|, \psi(h)+n_\ell]$ and 
\begin{eqnarray*}
\Delta(h) & = & \Big\{ (\ell,\ell',i_1,i_2) \, : \, \ell \ge 0, \, \ell'\geq 0, \, n_{\ell'} \in X(\psi(h),\ell), \, (i_1,i_2) \in Y_{\ell'}^\ell(h) \Big\}, \\ 
Y_{\ell'}^\ell(h) & = & \Big\{ (i_1,i_2) \in \{1,\ldots,N_{n_\ell}\} \times \{1,\ldots,N_{n_{\ell'}}\} \, : \, w_{n_\ell}(i_1)^{-1} w_{n_{\ell'}}(i_2) = h \Big\}.
\end{eqnarray*}
for $h \ne e$. Here we have used that $|\psi(h)-\psi(g)|\leq \psi(gh)\leq \psi(h)+\psi(g)$, which follows from the subadditivity of $\psi$. Observe also that the set $X(\psi(h),\ell)$ is finite by the exponential growth property. We now estimate the left hand side of \eqref{Aim2-4} as follows
\begin{eqnarray} \label{Est1ThB}
\lefteqn{\hskip-5pt |\phi_\G(e)|^2 +  \sum_{h \neq e} |\phi_\G(h)|^2} \\ \nonumber \hskip15pt 
& \le & \Big( \sum_{k \ge 0} \alpha_{n_k} r^{2n_k} \Big)^2 + \, \sum_{h\neq e} \sum_{\delta, \delta' \in \Delta(h)} \, \gamma_\delta(h) \gamma_{\delta'}(h) \nu_\delta(h,r) \nu_{\delta'}(h,r) \\ \nonumber \hskip15pt  & \le & \Big( \sum_{k \ge 0} \alpha_{n_k} r^{2n_k} \Big)^2 + \, \sum_{h \neq e} \sum_{\delta, \delta' \in \Delta(h)} \frac12\big( \gamma_\delta(h)^2 + \gamma_{\delta'}(h)^2 \big) \nu_\delta(h,r) \nu_{\delta'}(h,r) \\ \nonumber \hskip15pt & = & \Big( \sum_{k \ge 0} \alpha_{n_k} r^{2n_k} \Big)^2 + \, \sum_{h \neq e} \Big[ \sum_{\delta \in \Delta(h)} \gamma_\delta(h)^2 \nu_\delta(h,r) \Big] \Big[ \sum_{\delta \in \Delta(h)} \nu_{\delta}(h,r) \Big].    
\end{eqnarray}
Our argument in the sequel rests on the following claim for $0 \le r < \rho^{- \frac12}$
\begin{equation} \label{claim1ThB}
\sum_{\delta \in \Delta(h)}^{\null} \nu_{\delta}(h,r) \,\leq \, \Psi_\G(\psi(h),r) \, = \, 
C (\rho r)^{\psi(h)} \Big[ 1 + \rho \frac{2-\rho r^2}{1-\rho r^2} \Big]. 
\end{equation}
Let us complete the proof before justifying the claim. Given integers $0 \le j \le k$ with $(j,k) \neq (0,0)$, fix $(u_1,u_2)$ satisfying $1 \le u_1 \le N_{n_j}$ and $1 \le u_2 \le N_{n_k}$. Let $h_{j,k;u_1,u_2}=w_{n_j}(u_1)^{-1}w_{n_k}(u_2)$ and consider the following sets for $h\neq e$ $$\Delta[j,k;u_1,u_2](h) = \Big\{ \delta \in \Delta(h) \, : \, \big\{ w_{n_\ell}(i_1), w_{n_{\ell'}}(i_2) \big\} = \big\{ w_{n_j}(u_1), w_{n_k}(u_2)\big\} \Big\}.$$
Given $h \neq e$, we note that 
\begin{itemize}
\item $\gamma_\delta(h) = a_{n_j}(u_1) a_{n_k}(u_2)$ for $\delta \in \Delta[j,k;u_1,u_2](h)$,

\vskip3pt

\item $\Delta[j,k;u_1,u_2](h) = \emptyset$ for $h \notin \{h_{j,k;u_1,u_2},h_{j,k;u_1,u_2}^{-1}\}$,

\vskip3pt

\item $\big| \Delta[j,k;u_1,u_2](h_{j,k;u_1,u_2}) \big|=\big| \Delta[j,k;u_1,u_2](h_{j,k;u_1,u_2}^{-1}) \big| \le 2$ ,

\vskip3pt

\item $\Delta(h)$ decomposes as the disjoint union of $\Delta[j,k;u_1,u_2](h)$'s.
\end{itemize}
We may now continue using Fubini as follows
\begin{eqnarray*} 
\lefteqn{\hskip-10pt \sum_{h \neq e} \Psi_\G(\psi(h),r) \Big[ \sum_{\delta \in \Delta(h)} \gamma_\delta(h)^2 \nu_\delta(h,r) \Big]} \\  [6pt] \nonumber & = & \sum_{h \neq e} \Psi_\G(\psi(h),r) \Big[ \sum_{\begin{subarray}{c} j \le k \\ (j,k) \neq (0,0) \end{subarray}} \, \sum_{\begin{subarray}{c} 1 \le u_1 \le N_{n_j} \\ 1 \le u_2 \le N_{n_k} \\ \delta \in \Delta[j,k;u_1,u_2](h) \end{subarray}} a_{n_j}(u_1)^2 a_{n_k}(u_2)^2 \, r^{n_j+n_k} \Big] \\  \nonumber & \le & 4 \sum_{\begin{subarray}{c} j \le k \\ (j,k) \neq (0,0) \\ h_{j,k;u_1,u_2} \neq e \end{subarray}} \, \sum_{\begin{subarray}{c} 1 \le u_1 \le N_{n_j} \\ 1 \le u_2 \le N_{n_k} \end{subarray}} \Psi_\G \big( \psi(h_{j,k;u_1,u_2}),r \big) a_{n_j}(u_1)^2 a_{n_k}(u_2)^2 \, r^{n_j+n_k}  \\ [6pt] & \le & \nonumber 4 \sum_{\begin{subarray}{c} j \le k \\ (j,k) \neq (0,0) \\ h_{j,k;u_1,u_2} \neq e \end{subarray}}^{\null} \, \Big( \underbrace{\max_{\begin{subarray}{c} 1 \le u_1 \le N_{n_j} \\ 1 \le u_2 \le N_{n_k} \end{subarray}} \Psi_\G \big( \psi(h_{j,k;u_1,u_2}), r \big)}_{\Sigma(n_j,n_k,r)} \Big) \, \alpha_{n_j} \alpha_{n_k} \, r^{n_j+n_k}.  
\end{eqnarray*}
Recalling that $\max \{n_k-n_j,\sigma\} \le \psi(h_{j,k;u_1,u_2}) \le n_j + n_k$ and the function $\Psi_\G$ is decreasing in the first variable when $\rho r <1$ and increasing when $\rho r > 1$, we obtain the following estimate $$\Sigma(n_j,n_k,r) \, \le \, \begin{cases} \Psi_\G(n_j+n_k,r) & \mbox{if } \rho r > 1, \\ \min \big\{ \Psi_\G(\sigma,r), \Psi_\G(n_k-n_j,r) \big\} & \mbox{if } \rho r < 1. \end{cases}$$ Moreover, it follows from \eqref{Est1ThB} and claim \eqref{claim1ThB} that
\begin{eqnarray*} 
 \sum_{h \in \G} |\phi_\G(h)|^2  & \le & \alpha_0^2 + \sum_{k \ge 1} \big( \underbrace{r^{4n_k} + 4 \Sigma(n_k,n_k,r) r^{2n_k}}_{A_{n_k}(r)} \big) \alpha_{n_k}^2 \\ & + & 2 \sum_{0 \le j<k} \big( \underbrace{r^{2(n_j+n_k)} + 2 \Sigma(n_j,n_k,r) r^{n_j+n_k}}_{B_{n_j,n_k}(r)} \big) \alpha_{n_j} \alpha_{n_k}.
\end{eqnarray*}
Therefore, inequality \eqref{Aim2-4} reduces to \eqref{claim1ThB} and the upper bound 
\begin{equation} \label{claim2ThB}
\max \Big\{ \sup_{k \ge 1} \ A_{n_k}(r), \ \sup_{j < k} \ B_{n_j,n_k}(r) \Big\} \ \le \ 1 \quad \mbox{for} \quad 0 \le r \le \mathcal{R}_{\G,\psi} < \rho^{-\frac12}.
\end{equation}
Let us finally prove \eqref{claim1ThB} and \eqref{claim2ThB}. To prove \eqref{claim1ThB} we first observe that $$\sum_{\substack{\ell' \geq 0\\ n_{\ell'} \in X(\psi(h),\ell)}}  |Y_{\ell'}^\ell(h)|=|W_{n_\ell}|=N_{n_\ell} \Rightarrow \sum_{\substack{\ell' \geq 0\\ n_{\ell'} \in X(\psi(h),\ell)}}  \frac{|Y_{\ell'}^\ell(h)|}{N_{n_\ell}}=1.$$ Then we write
\begin{eqnarray*}
\lefteqn{\hskip-5pt \sum_{\delta \in \Delta(h)} \nu_{\delta}(h,r)} \\ & = & \sum_{\ell \ge 0} \sum_{\substack{\ell' \geq 0\\ n_{\ell'} \in X(\psi(h),\ell)}}  |Y_{\ell'}^\ell(h)|r^{n_\ell+n_{\ell'}} \\ & = & \sum_{\ell \ge 0} N_{n_\ell}r^{n_\ell} \sum_{\substack{\ell' \geq 0\\ n_{\ell'} \in X(\psi(h),\ell)}}  \frac{|Y_{\ell'}^\ell(h)|}{N_{n_\ell}}r^{n_{\ell'}} \\ & \le & \sum_{\ell \ge 0} N_{n_\ell}r^{n_\ell} \max \big\{ r^{n_{\ell'}} \, : \, n_{\ell'} \in X(\psi(h),\ell) \big\} \ = \ \sum_{\ell \ge 0} N_{n_\ell}r^{n_\ell+|\psi(h)-n_{\ell}|} \\ [5pt] & = & \underbrace{\sum_{\begin{subarray}{c} \ell \ge 0 \\ n_\ell \le \psi(h) \end{subarray}} N_{n_\ell}r^{\psi(h)}}_{(\mathrm{a})} + \underbrace{\sum_{\begin{subarray}{c} \ell \ge 0 \\ n_\ell > \psi(h) \\ n_\ell < \lfloor \psi(h)\rfloor +1 \end{subarray}} N_{n_\ell}r^{2n_\ell-\psi(h)}}_{(\mathrm{b})} + \underbrace{\sum_{\begin{subarray}{c} \ell \ge 0 \\ n_\ell \ge \lfloor \psi(h)\rfloor +1 \end{subarray}} N_{n_\ell}r^{2n_\ell-\psi(h)}}_{(\mathrm{c})}.
\end{eqnarray*}
To estimate these terms we use the exponential growth assumption
\begin{eqnarray*}
(\mathrm{a}) & = & r^{\psi(h)} \big| \big\{ g \in \G \; : \; \psi(g) \le \psi(h) \big\} \big| \ \le \ C (\rho r)^{\psi(h)}, \\ [18pt] (\mathrm{b}) & = & r^{-\psi(h)}\sum_{\substack{g\in \G\\ \psi(h) <\psi(g) < \lfloor \psi(h)\rfloor +1}}r^{2 \psi(g)} \ \le \ C\rho (\rho r)^{\psi(h)}, \\
(\mathrm{c}) & = & r^{-\psi(h)}\sum_{j\geq  \lfloor \psi(h)\rfloor +1}\sum_{\substack{k\geq 0\\ j\leq n_k<j+1}}\sum_{\substack{g\in \G \\ \psi(g)=n_k}}r^{2n_k}\\
& \le &  r^{-\psi(h)}\sum_{j\geq  \lfloor \psi(h)\rfloor +1}C\rho^{j+1}r^{2j} \ \le \ \frac{C \rho}{1-\rho r^2}(\rho r)^{\psi(h)},
\end{eqnarray*}
for $0 \le r < \rho^{-\frac12}$. This justifies our claim \eqref{claim1ThB}. To prove \eqref{claim2ThB} we will need to distinguish two cases according to the value of $\rho r$. Namely, when $\rho r > 1$ we know that $\Sigma(n_j,n_k,r) \le \Psi_\G(n_j+n_k,r)$ and $\Sigma(n_j,n_k,r) \le \min \{ \Psi_\G(\sigma,r), \Psi_\G(n_k-n_j,r) \}$   otherwise. This means that 
\begin{itemize}
\item If $\rho r > 1$, we find
\begin{eqnarray*}
A_{n_k}(r) & = & r^{4n_k} + 4 C \Big[ 1 + \rho \frac{2-\rho r^2}{1-\rho r^2} \Big] (\rho r^2)^{2n_k}, \\ 
B_{n_j,n_k}(r) & = & r^{2(n_j+n_k)} + 2 C \Big[ 1 + \rho \frac{2-\rho r^2}{1-\rho r^2} \Big] (\rho r^2)^{n_j+n_k}.
\end{eqnarray*}

\item If $\rho r < 1$, we find   
\begin{eqnarray*}
A_{n_k}(r) & \le & r^{4n_k} + 4 C \Big[ 1 + \rho \frac{2-\rho r^2}{1-\rho r^2} \Big] (\rho r)^{\sigma}  r^{2n_k}, \\ 
B_{n_j,n_k}(r) & \le & r^{2(n_j+n_k)} + 2 C \Big[ 1 + \rho \frac{2-\rho r^2}{1-\rho r^2} \Big] (\rho r)^{n_k-n_j} r^{n_j+n_k}.
\end{eqnarray*} 
\end{itemize}
Using $n_k \ge n_1 = \sigma$ and $\rho r^2 < 1$, we obtain in both cases $$\sup_{k \ge 1} A_{n_k}(r) \, \le \, r^{4\sigma} + 4 C \Big[ 1 + \rho \frac{2-\rho r^2}{1-\rho r^2} \Big] (\rho r^2)^{2\sigma} \, =: \, A_\sigma(r),$$ $$\sup_{j < k} B_{n_j,n_k}(r) \, \le \, r^{2 \sigma} + 2 C \Big[ 1 + \rho \frac{2-\rho r^2}{1-\rho r^2} \Big] (\rho r^2)^{\sigma} \ =: \ B_\sigma(r).$$ In conclusion, we always have the estimate $$\max \Big\{ \sup_{k \ge 1} \ A_{n_k}(r), \ \sup_{j < k} \ B_{n_j,n_k}(r) \Big\} \, \le \, \max \big\{ A_\sigma(r), B_\sigma(r) \big\} \, \to \, 0 \quad \mbox{as} \quad r \to 0.$$
This means we can always find $0 < \mathcal{R}_{\G,\psi} < \rho^{-1/2}$ such that inequality \eqref{claim2ThB} holds. Therefore, the assertion follows by taking $\beta(\G,\psi) = -2 \sigma \log (\mathcal{R}_{\G,\psi})$. \fin

\subsection{Behavior of the constant $\beta(\G,\psi)$}

As we explained in the Introduction, the expected optimal time is $\frac{1}{2 \sigma} \log (q-1/p-1)$ for many Poisson-like lengths. It is clear from our argument that the constant $\beta(\G,\psi)$ that we get might be far from the expected optimal value 1. It is however interesting to study how this constant depends on $\rho$, the exponential growth order of the pair $(\G,\psi)$. For instance, given a discrete group $\G$ finitely generated by $\mathrm{S} \subset \G$, we always have $N_k(\mathrm{S}) \le 2|\mathrm{S}| (2|\mathrm{S}|-1)^{k-1}$ for the number $N_k(\mathrm{S})$ of elements in $\G$ with word length $k$ with respect to the alphabet $\mathrm{S}$. This gives a rough bound 
\begin{eqnarray*}
\big| \big\{ g \in \G \, : \, \psi(g) \le R \big\} \big| & \le & 1 + \frac{2|\mathrm{S}|}{2|\mathrm{S}|-1} \sum_{k=1}^R (2|\mathrm{S}|-1)^k \\ & \le & \frac{|\mathrm{S}|}{|\mathrm{S}|-1} (2|\mathrm{S}|-1)^R \ =: \ C \rho^R 
\end{eqnarray*}
which is more and more accurate for groups with fewer relations. In particular, it is interesting to know the behavior of the constant $\beta(\G,\psi)$ for $C$ fixed and $\rho$ large. 

\begin{corollary} Assume that $\psi$ is Poisson-like with $$\big| \big\{ g \in \G \; : \; \psi(g) \le R \big\} \big| \, \le \, C \rho^R.$$ Then, the following estimates hold for $\rho$ large compared to $C \! :$ 
\begin{itemize}
\item[i)] $\beta(\G,\psi) \le \eta \log(\rho)$ for any $\eta > 1 + \sigma$.

\item[ii)] If $\psi: \G \to \Z_+$ and $N_k=|\{g \in \G \, : \, \psi(g) = k\}| \le C_0 \rho^k$, then $$\beta(\G,\psi) \le \eta \log(\rho) \quad \mbox{for any} \quad \eta > 1.$$
\end{itemize}
\end{corollary}

\dem The first assertion is equivalent to prove \eqref{claim2ThB} for $$r = \mathcal{R}_{\G,\psi} = \Big(\frac{1}{\sqrt{3}}\Big)^{\beta(\G,\psi)/\sigma \log 3} = \rho^{-\eta/2\sigma}$$ and $\rho$ large. In other words, we need to satisfy the inequalities
\begin{itemize}
\item $\rho r^2 = \rho^{1 - \eta/\sigma} < 1$,

\vskip4pt

\item $A_\sigma(r) = r^{4\sigma} + 4 C \Big[ 1 + \rho \frac{2-\rho r^2}{1-\rho r^2} \Big] (\rho r^2)^{2\sigma} \le 1$,

\item $B_\sigma(r) = r^{2 \sigma} + 2 C \Big[ 1 + \rho \frac{2-\rho r^2}{1-\rho r^2} \Big] (\rho r^2)^{\sigma} \le 1$.
\end{itemize}
If $\rho$ is large enough, it suffices to have $\eta > \max \{\sigma, \sigma + \frac12 , \sigma + 1\} = 1 + \sigma$. This completes the proof of i). To prove the second assertion we notice that in that case $\sigma = 1$. Moreover, in the proof of Theorem B we may replace \eqref{claim1ThB} by $$\sum_{\delta \in \Delta(h)}^{\null} \nu_{\delta}(h,r) \, \le \, (\mathrm{a}) + (\mathrm{b}) + (\mathrm{c}) \,\le \, \Phi_\G(\psi(h),r) \, := \, \Big[ C + \frac{C_0 \rho r^2}{1-\rho r^2} \Big] (\rho r)^{\psi(h)}$$ for $0 \le r < \rho^{-\frac12}$. Indeed, we still have $(\mathrm{a}) \le C (\rho r)^{\psi(h)}$, $(\mathrm{b})=0$ and $$(\mathrm{c}) \, = \, \sum_{\ell \geq \psi(h)+1} N_\ell r^{2\ell-\psi(h)} \, \le \, r^{-\psi(h)} \sum_{\ell \geq \psi(h)+1}  C_0 (\rho r^2)^\ell \, = \, \frac{C_0 \rho r^2 (\rho r)^{\psi(h)}}{1-\rho r^2}.$$ Using this new estimate with $\Phi_\G$ instead of $\Psi_\G$, we may proceed as in the proof of Theorem B to obtain a proof of 
\eqref{claim2ThB} which ultimately only requires the condition $\rho r^2 < 1$. If $\rho$ is large enough this follows from $\eta > \sigma = 1$ as desired. \fin

\begin{remark} \label{Polynomial}
\emph{Of course, Theorem B also applies for Poisson-like lengths with polynomial growth. However, our constants could be improved by giving a finer estimate of the sum $\sum_{\delta \in \Delta(h)} \nu_\delta(h,r)$. A quick look at our proof shows that the estimates for terms (a) and (b) can be easily adapted. Only the term (c) requires precise bounds for truncated power series with polynomial coefficients.} 
\end{remark}

\subsection{Examples of Poisson-like lengths}

Compared to the method used to prove Theorems A1-A3, our new argument for Theorem B is simpler and does not rely on computational estimates. Therefore, at the price of a worse constant (analyzed above) it applies at once to any Poisson-like length, even those which admit small loops. This also holds for Theorem C. Let us provide a few illustrations.

\begin{itemize}
\item[{\bf i)}] {\bf The word length.} Perhaps, the most significant application of Theorems B and C is that they provide a hypercontractivity/ultracontractivity bound for any finitely generated discrete group $\G$ equipped with the word length $|\cdot|_{\mathrm{S}}$ associated to any set $\mathrm{S}$ of generators. Indeed, any such length is Poisson-like. Namely, it is always subadditive, admits the spectral gap $\sigma = 1$ and has (at most) exponential growth $\rho = 2 |\mathrm{S}| - 1$. The conditional negativity might fail and lead to a non Markovian semigroup, see iv) below. For instance, Theorem B applies to the discrete Heisenberg group whose Cayley graph admits small loops. This could be used as in \cite{JMP} to obtain some related estimates on noncommutative tori, but we will recall in v) below a much simpler argument for these algebras.       

\vskip5pt

\item[{\bf ii)}] {\bf Finite-dimensional proper cocycles.} A natural question is determining the class of discrete groups which admit a Poisson-like length. We can provide a positive answer for groups admitting a finite-dimensional proper cocycle. A cocycle is a triple $(\mathcal{H},\alpha, b)$ given by a real Hilbert space $\mathcal{H}$, an orthogonal representation $\alpha:\G \to O(\mathcal{H})$ and a map $b:\G\to \mathcal{H}$ satisfying the cocycle law $b(gh)=\alpha_g(b(h))+b(g)$. It is called finite-dimensional when $\dim \mathcal{H} < \infty$ and proper when $\{g\in \G : \|b(g)\|_{\mathcal{H}}\leq R\}$ is finite for any $R > 0$. Let $\G$ be a discrete group admitting a finite-dimensional proper cocycle $(\mathcal{H},\alpha, b)$. We claim that $$\hskip20pt \psi(g) \, = \, \|b(g)\|_{\mathcal{H}} + \delta_{g\neq e}$$ defines a Poisson-like length. Let us first prove that it is conditionally negative. By Schoenberg's theorem, $g \mapsto \|b(g)\|^2_{\mathcal{H}}$ defines a conditionally negative length. Since conditional negativity is stable under square root and sum, we conclude that $\psi$ is conditionally negative. On the other hand, according to the cocycle law it is easy to see that $\psi$ is subadditive. Moreover, the spectral gap condition is ensured by the term $\delta_{g\neq e}$. Hence, it remains to show that $\psi$ has (at most) exponential growth. In fact, we will show that $\psi$ has polynomial growth. If $\dim \H = n$, it suffices to prove that $|b(\G)\cap \mathrm{B}_R(0)| \lesssim R^n$ where $\mathrm{B}_R(0) = \{ \xi \in \H \, | \, \|\xi\|_\H \le R \}$. Indeed, using the cocycle law we immediately deduce $$\hskip20pt \|b(g^{-1}h)\|_\H = \big\| b(g) - b(h) \big\|_{\mathcal{H}}.$$ In particular, since $b$ is proper we see that $$\hskip20pt \big| \big\{ h \in \G \, : \, b(h) = b(g) \big\} \big| \, = \, \big| \big\{ h \in \G \, : \, b(h) = 0 \big\} \big| \, \le \, K \, < \, \infty$$ for any $g \in \G$, which gives $|\{g \in \G : \|b(g)\| \le R \}| \le K \cdot |b(\G)\cap \mathrm{B}_R(0)|$. On the other hand,  $\mathrm{B}_R(0)$ can be covered by approximately $R^n$ translates of $\mathrm{B}_1(0)$, let us call them $\mathrm{B}_1(\xi_j)$ for $1 \le j \le R^n$. Therefore, it is enough to show that $|b(\G) \cap \mathrm{B}_1(\xi_j)|$ is uniformly bounded in $j$. However, any two points $b(g), b(h)$ in this set satisfy $\|b(g^{-1}h)\|_\H = \|b(g) - b(h)\|_\H \le 2$, which immediately gives $$\hskip20pt \sup_{\xi \in \H} \big| b(\G) \cap \mathrm{B}_1(\xi) \big| \, \le \, 1 + \big| b(\G) \cap \mathrm{B}_2(0) \big| \, < \, \infty.$$ This proves our assertion. Recall that $\G$ is called a-T-menable or is said to satisfy the Haagerup property when it admits a proper cocycle (not necessarily finite-dimensional), see e.g. \cite{CCJJV} for more on this property. At the time of this writing, we do not know if we may find a Poisson-like length for any group satisfying the Haagerup property. Let us note from \cite{CTV} that a-T-menable groups admit conditionally negative lengths with arbitrary slow growth, which lead to arbitrary large growth of the group. 

\vskip5pt

\item[{\bf iii)}] {\bf Kazhdan's property (T).} The Haagerup property is a strong negation of the so-called Kazhdan's property (T), which refers to those groups whose conditionally negative length functions are all bounded, see \cite{HV} for more on this notion.  It is clear that infinite groups satisfying Kazhdan's property (T) do not admit Poisson-like lengths, since the exponential growth can not hold. Therefore, Theorem B does not give hypercontractivity estimates for this class of groups. In fact, it can be shown that no hypercontractivity bound is possible for any pair $(\G,\psi)$ with $\G$ an infinite group satisfying Kazhdan's property (T) and $\psi$ a conditionally negative length. Assume that $\psi(g) \le K$ for all $g\in \G$ and that the pair $(\G,\psi)$ satisfies a hypercontractive estimate. Then, arguing via log-Sobolev inequalities as in Appendix A, we may find $t(p) \ge  0$ for $1 < p < 2$ satisfying $$\hskip20pt e^{-t(p)K} \|f\|_2 \le \|\S_{\psi,t(p)}f\|_2 \le \|f\|_p \le \|f\|_2$$ for any $f \in L_p^+(\L(\G))$. This implies that the $L_p$-norm is equivalent to the $L_2$-norm for all $1 < p \neq 2 < \infty$, which contradicts the fact that $\L(\G)$ is an infinite-dimensional algebra.    

\vskip5pt

\item[{\bf iv)}] {\bf Non-Markovian semigroups.} If $\psi$ is not conditionally negative, the semigroup $\S_{\psi,t}$ may not be positive preserving. In that case, Theorem B may yield interesting $L_2 \to L_4$ estimates. Note that we can interpolate with the (trivial) contraction $\S_{\psi,t}: L_2 \to L_2$. However, Markovianity seems to be inherent to Gross extrapolation argument, so we may not produce more general hypercontractivity bounds using that method. On the contrary, our ultracontractivity estimates in Theorem C hold even in the non-Markovian case. Moreover, since the condition $\psi(g) = \psi(g^{-1})$ implies self-adjointness of the semigroup $(\S_{\psi,t})_{t \ge 0}$, we deduce $\|\S_{\psi,2t}\|_{1 \to \infty}  \le \|\S_{\psi,t}\|_{2 \to \infty}^2$. Therefore Theorem C gives apparently new estimates for non-Markovian semigroups $$\|\S_{\psi,t}f\|_\infty \, \le \, \Big[ 1 + 2 e^{-t/4} \Big( \sum_{g \neq e} e^{-t \psi(g)/2} \Big)^{\frac12} \Big]^2 \|f\|_1.$$

\vskip5pt

\item[{\bf v)}] {\bf Classical and noncommutative tori.} As recalled in the Introduction, Weissler \cite{W} obtained the optimal time hypercontractivity inequalities for the one-dimensional Poisson semigroup in $L_\infty(\mathbb{T})=\L(\Z)$. One may think of two possible extensions of this result to the $n$-dimensional case. On the one hand, regarding $\mathbb{T}$ as the one-dimensional sphere $\mathbb{S}^1$, we may consider the Poisson semigroup on $\mathbb{S}^n$. Beckner obtained optimal hypercontractivity bounds in this context \cite{Be2}. The other direction is to consider $n$-dimensional tori $L_\infty(\mathbb{T}^n)=\L(\Z^n)$. By using Bonami's induction trick for cartesian products, Weissler's theorem extends to $\G = \Z^n$ equipped with the word length $|\k|_1=\sum_j |k_j|$. We might also consider the lengths $$\hskip20pt |\k|_u=\Big(\sum_{j=1}^n |k_j|^u\Big)^{\frac{1}{u}}$$ for any $1 \le u \le 2$, and the associated Poisson-like semigroups $$\hskip20pt \S_{|\cdot|_u,t}f = \sum_{\k \in \Z^n} e^{-t|\k|_u} \widehat{f}(\k) e^{2\pi i \langle \k, \cdot\rangle}.$$ Comparing lengths, it follows directly from Weissler's result for $u=1$ that $$\hskip20pt \S_{|\cdot|_u,t}: L_p(\L(\Z^n)) \to L_q(\L(\Z^n)) \quad \mbox{for} \quad  t \ge \frac{n^{1-1/u}}{2} \log\Big(\frac{q-1}{p-1} \Big).$$ If $u=2$ we obtain the usual Poisson semigroup and hypercontractivity estimates with a constant $\sqrt{n}$. It is reasonable to expect optimal estimates here, but (to the best of our knowledge) these are still open. Adapting Theorem B for polynomial growth lengths could improve this constant, see Remark \ref{Polynomial}. On the other hand, if $q \in 2 \Z_+$ the hypercontractivity $L_2 \to L_q$ estimates that we get for $\Z^n$ are easily transferred to noncommutative tori $\mathcal{A}_\Theta$ since noncommutativity in the algebraic expressions for the $q$-norms only adds an additional phase to the Fourier coefficients, which disappears by unconditionality in $L_2$. Moreover, Gross extrapolation technique also applies in this case, so that we obtain hypercontractive inequalities for arbitrary indices $1 < p \le q < \infty$. We refer e.g. to \cite{JMP} for the formal definition of these algebras.  
\end{itemize}


\subsection{Ultracontractivity}

Given $q > 2$, the hypercontractivity problem consists in finding the optimal time $t_{2,q}$ beyond which we obtain $L_2 \to L_q$ contractions. A related problem is to find an absolute constant $C_{\G,\psi}(t)$ so that $$\|\S_{\psi,t}f\|_\infty \ = \, \sup \Big\{ \|\mathcal{S}_{\psi,t}f\|_q \; : \; q >2 \Big\} \, \le \, C_{\G,\psi}(t) \|f\|_2.$$ If such a constant exists for some $t_0 > 0$, then the semigroup property and the contractivity of $\S_{\psi,t}: \L(\G) \to \L(\G)$ imply that it also exists for any larger $t$ and $C_{\G,\psi}(t)$ is a decreasing function for $t \ge t_0$. According to our lower bound for the optimal hypercontractivity time $\frac{1}{2\sigma} \log(q-1)$, it is also easy to show that we must have $C_{\G,\psi}(t) > 1$ for all $t > 0$ and $C_{\G,\psi}(t) \to \infty$ as $t \to 0^+$. This behavior of the semigroup was studied for the first time by Edward B. Davies and Barry Simon who called it ultracontractivity, see \cite{CM,DaSi}. Shortly after and motivated by the work of Moser \cite{M1,M2}, Varopoulos investigated in \cite{V} the relation between this behavior of the semigroup and Sobolev type inequalities associated to infinitesimal generators of Markovian semigroups. We also refer to the work of Jolissaint \cite{J1,J2} in the context of group von Neumann algebras. In our context, ultracontractivity for the free Poisson semigroup is a known consequence of Haagerup's inequality for homogeneous polynomials \cite{Ha} $$\Big\| \sum_{|g|=k} \widehat{f}(g) \lambda(g) \Big\|_{L_\infty(\L(\F_\infty))} \, \le \, (k+1) \, \Big( \sum_{|g|=k} |\widehat{f}(g)|^2 \Big)^{\frac12}.$$ Indeed, together with Cauchy-Schwarz inequality we obtain 
\begin{eqnarray} \label{HaagerupIneq}
\|\P_t f \|_\infty & \le & \sum_{k \geq 0}  e^{-tk} (k+1) \Big( \sum_{|g|=k} |\widehat{f}(g)|^2 \Big)^{\frac12} \\ \nonumber & \le & \Big(\sum_{k \geq 0}  e^{-2tk}(k+1)^2\Big)^{\frac12} \|f\|_2 \ = \  \Big(\frac{1+e^{-2t}}{(1-e^{-2t})^3}\Big)^{\frac12} \|f\|_2.
\end{eqnarray}
The following result establishes a similar behavior for arbitrary groups/lengths.

\begin{corollary} \label{CorollaryC}
We have $$\|\S_{\psi,t}f\|_\infty \, \le \, \Big[ 1 + 2 e^{-t/2} \Big( \sum_{g \neq e} e^{-t \psi(g)} \Big)^{\frac12} \Big] \|f\|_2,$$ for any discrete group $\G$ and any length $\psi$ satisfying $\psi(e)=0$ and $\psi(g) = \psi(g^{-1})$. 
\end{corollary}

\dem Assume first that $\psi$ is an integer-valued Poisson-like length. It suffices to show that $\|\S_{\psi,t}f\|_q \le C_{\G,\psi(t)} \|f\|_2$ for $q \in 2\Z_+$ arbitrarily large. Decompose $f$ as follows $$f = f_1 + f_2 = \big( f - \widehat{f}(e) \mathbf{1} \big) + \big( \widehat{f}(e) \mathbf{1} \big).$$ Clearly $\|\mathcal{S}_{\psi,t} f_2 \|_q = |\widehat{f}(e)| \le \|f\|_2$ for any $q$. We claim that 
\begin{equation}\label{f1}
\|\mathcal{S}_{\psi,t} f_1 \|_{L_\infty(\mathcal{L}(\G))} \, \le \, 2 e^{-t/2} \Big( \sum_{g \neq e} e^{-t\psi(g)} \Big)^{\frac12} \|f_1\|_{L_2(\mathcal{L}(\G))}. 
\end{equation}
Since $\|f_1\|_2\leq \|f\|_2$ this will prove the assertion for Poisson-like lenghts. The key point is that $f_1$ is mean zero, which makes the computations much easier. Namely, in the terminology of the combinatorial method in Section \ref{Method} (which is still valid for non-integer-valued lengths) we deduce $$\widehat{f_1}(e)=0 \ \Rightarrow \ \|\S_{\psi,t}f_1\|_q^q = \widehat{f}_1(e)^q +\sum_{u=1}^q \binom{q}{u} \widehat{f}_1(e)^{q-u} s_u(r) = s_q(r).$$ Following the proof of Lemmas \ref{EvenLemma} and \ref{method2 better than method1}, we get for $r=e^{-t}$ $$\|\S_{\psi,t} f_1\|_\infty = \lim_{q \to \infty} s_q(r)^{\frac{1}{q}} \le \lim_{s \to \infty} \Big( \sum_{\k \in L_s(s)} \sum_{\delta \in \Delta_{2s}[\k]} \gamma_{\delta} \nu_\delta(r) \Big)^{\frac{1}{2s}}$$ by declaring $B_q = \emptyset$. On the other hand, Proposition \ref{method1} i) for $m=s$ gives 
\begin{eqnarray*}
\|\S_{\psi,t} f_1\|_\infty \!\!\!\! & \le & \!\!\!\! \lim_{s \to \infty} \Big( \sum_{\k \in L_s(s)} \Gamma(s,r) r^{|\k|+1} M(\k) \alpha_{\k} \Big)^{\frac{1}{2s}} \\ \!\!\!\! & \le & \!\!\!\! \lim_{s \to \infty} \Big( \, \Gamma(s,r) \, r^{s+1} \sum_{\k \in L_s(s)} M(\k) \, \alpha_{\k} \Big)^{\frac{1}{2s}} = \lim_{s \to \infty} \Big( \, \Gamma(s,r) \, r^{s+1} \Big)^{\frac{1}{2s}} \|f_1\|_2 
\end{eqnarray*}
where $$\Gamma(s,r) \, = \, \frac{(2s)!}{(s!)^2} \mathcal{G}(\G,\psi,r)^{s-1} \, = \, \frac{(2s)!}{(s!)^2} \Big( \sum_{g \neq e} r^{\psi(g)} \Big)^{s-1}.$$ According to Stirling's formula we get $$\|\S_{\psi,t} f_1\|_\infty \, \le \, 2 e^{-t/2} \Big( \sum_{g \neq e} e^{-t\psi(g)} \Big)^{\frac12} \|f_1\|_2.$$ This completes the proof for Poisson-like lengths. However, our assumptions only impose $\psi(e) = 0$ and $\psi(g) = \psi(g^{-1})$ for all $g \in \G$. In this general setting, we need to eliminate our prior assumptions which include conditional negativity, subadditivity and the existence of a spectral gap with $\sigma=1$ since we also imposed $\psi$ to be integer valued. First, subadditivity was only relevant to define the set of admissible lengths which did not play any role in the argument above. Second, we only used $\sigma = 1$ when applying Proposition \ref{method1}, but the latter can be easily adapted to the case $\sigma \neq 1$ (even $\sigma=0$) just replacing $r^{|\k|+1}$ by $r^{|\k|}$. Finally, to remove conditional negativity we observe that the mean-zero condition of $f_1$ still gives
\begin{eqnarray*}
\|\S_{\psi,t}f_1\|_q^q & = & \tau \Big[ \Big| \sum_{g,h \in \G} \overline{\widehat{f}(g)} \widehat{f}(h) \lambda(g^{-1}h) \Big|^s \Big] \\ & \le & \sum_{\begin{subarray}{c} g_1 g_2 \cdots g_q = e \\ g_j \neq e \end{subarray}} \prod_{j=1}^s |\widehat{f}(g_{2j-1}^{-1})| e^{-t \psi(g_{2j-1})} |\widehat{f}(g_{2j})| e^{-t \psi(g_{2j})}.
\end{eqnarray*}
Therefore, we get an asymmetric form of $s_q(r)$ which still can be bounded by following the ideas of Lemmas \ref{EvenLemma} and \ref{method2 better than method1} and Proposition \ref{method1}. We leave to the reader the (easy) task of adapting the arguments there to this case. \fin

Our estimate in Corollary \ref{CorollaryC} is very close to the trivial estimate which arises by applying Cauchy-Schwarz inequality. Namely, if we set $N_R$ for the cardinality of points in the $\psi$-sphere of radius $R$ then we may trivially obtain the estimate $\|f\|_\infty \le \sqrt{N_R} \|f\|_2$ for any $f$ supported by that $\psi$-sphere. Then, arguing as in \eqref{HaagerupIneq} we deduce the following estimate
\begin{eqnarray*}
\|\S_{\psi,t}f\|_\infty & \le & \Big( \sum_{R \in \psi(\G)} e^{-2tR} N_R \Big)^{\frac12} \|f\|_2 \\ & = & \Big( \sum_{g \in \G} e^{-2t \psi(g)}  \Big)^{\frac12} \|f\|_2 \ \le \ \Big[ N_0^\frac12 + e^{-t \sigma/2} \Big( \sum_{g \notin N_0} e^{-t \psi(g)} \Big)^{\frac12} \Big] \|f\|_2,
\end{eqnarray*}
where $\sigma = \inf_{\psi(g) \neq 0} \psi(g)$ is the spectral gap after removing zero-length elements. A quick comparison between both estimates shows that Corollary \ref{CorollaryC} only gives substantially new information when there exits a large concentration around $0$. This means that our estimates are better when $N_0 > 1$ or $\sigma < 1$. 

\section*{{\bf Appendix}} 

\renewcommand{\theequation}{A.1}
\addtocounter{equation}{-1}

\subsection*{A. Logarithmic Sobolev inequalities}

In this appendix we shall adapt Gross extrapolation technique \cite{G1,G2} to obtain general hypercontractivity inequalities in group von Neumann algebras from $L_2 \to L_4$ ones. Given a conditionally negative length $\psi: \G \to \R_+$, consider the associated Markov semigroup $\S_{\psi,t}$ and assume that we know $\S_{\psi,t}: L_2(\L(\G)) \to L_4(\L(\G))$ for any $t$ greater than or equal to the expected optimal time $\frac{1}{2\sigma} \log 3$. The starting point is to use Stein's interpolation method \cite{Stein}. Namely, define $$\S_{\psi,z} f \, = \, \sum_{g \in \G} e^{-[(1-z)t_{2,2} + zt_{2,4}] \psi(g)} \widehat{f}(g) \lambda(g) \quad \mbox{with} \quad (t_{2,2}, t_{2,4}) = (0,\mbox{$\frac{1}{2\sigma} \log 3$}).$$ A straightforward application of Stein's interpolation gives
\begin{itemize}
\item $\displaystyle \|\S_{\psi,t}f\|_{q(t)} \le \|f\|_2$ for $q(t) = \frac{2 \log 3}{\log 3 - \sigma t}$ and $0 \le t \le \frac{1}{2\sigma} \log 3$,

\vskip3pt

\item $\displaystyle \|\S_{\psi,t}f\|_2 \le \|f\|_{p(t)}$ for $p(t) = \frac{2 \log 3}{\log 3 + \sigma t}$ and $0 \le t \le \frac{1}{2\sigma} \log 3$.
\end{itemize}
Moreover, according to the last inequality above, we deduce that $$\frac{d\Phi}{dt}(0) \ge 0 \quad \mbox{for} \quad \Phi(t) = \|f\|_{p(t)}^2 - \|\S_{\psi,t}f\|_2^2.$$ Indeed, it is a positive smooth function vanishing at $0$. Let us write $A_\psi$ to denote the infinitesimal generator of $\S_{\psi,t}$. Then, differentiating $\Phi$ at time $0$ produces the following inequality, known as logarithmic Sobolev inequality  
\begin{equation} \label{LogSobolev}
\tau \big( |f|^2 \log |f|^2 \big) - \|f\|_2^2 \log \|f\|_2^2 \, \le \, 2 \log 3 \, \big\langle f, A_{\psi}f \big\rangle.
\end{equation}
The next result that we need is the analog of Gross inequality for the generator $A_\psi$. This follows from the $L_p$-regularity of the associated Dirichlet form, which in turn was proved by Olkiewicz and Zegarlinski in the tracial case in \cite[Theorem 5.5]{OZ}. Namely, given $f \ge 0$ and $1 < p < \infty$, it follows that  
\renewcommand{\theequation}{A.2}
\addtocounter{equation}{-1}
\begin{equation} \label{Lp-regularity}
\big\langle f^{p/2}, A_{\psi} f^{p/2} \big\rangle \, \le \, \frac{p^2}{4(p-1)} \, \big\langle f, A_{\psi}f^{p-1} \big\rangle.
\end{equation}
Replacing $f$ by $f^{p/2}$ in \eqref{LogSobolev} combined with \eqref{Lp-regularity} gives
\renewcommand{\theequation}{A.3}
\addtocounter{equation}{-1}
\begin{equation} \label{Lp-Log-Sobolev}
\tau \big( f^p \log f^p \big) - \|f\|_p^p \log \|f\|_p^p \, \le \, \log 3 \, \frac{p^2}{2(p-1)} \, \big\langle f, A_{\psi}f \big\rangle
\end{equation}
for $f \ge 0$ and $1 < p < \infty$, which is nothing but an $L_p$-analog of the logarithmic Sobolev inequality. Once we have these estimates at hand, we may prove general hypercontractivity bounds as follows. Let us redefine $q(t)$ to be the expected optimal index $q$ for which we have $L_p \to L_q$ hypercontractivity at time $t$ up to a constant $\log 3$. In other words, set $$q(t) \, = \, 1 + (p-1) \exp \Big( \frac{2 \sigma t}{\log 3} \Big).$$ The goal is to show that
\renewcommand{\theequation}{A.4}
\addtocounter{equation}{-1}
\begin{equation} 
\|\S_{\psi,t}f\|_{q(t)} \, \le \, \|f\|_p \quad \mbox{for all} \quad t \ge 0.
\end{equation}
If we set $\Psi(t,p) = \|\S_{\psi,t}f\|_{q(t)}$, we clearly have that $\Psi(0,q(0)) = \Psi(0,p) = \|f\|_p$. In particular, it suffices to show that $\Psi(t,q(t))$ is a decreasing function of $t$. Moreover, since $\S_{\psi,t}$ has positive maximizers, we may assume that $f \ge 0$. Differentiating at time $t$, the result follows by applying \eqref{Lp-Log-Sobolev} for $(f,p) = (\S_{\psi,t}f, q(t))$.  

\renewcommand{\theequation}{\arabic{equation}}
\numberwithin{equation}{section}

\begin{remark}
\emph{We have adapted Gross extrapolation method to the algebras $\L(\G)$ assuming optimal time estimates for $(p,q) = (2,4)$, which is the case in Theorems A1-A3. In particular, the logarithmic Sobolev inequalities \eqref{LogSobolev} and \eqref{Lp-Log-Sobolev} hold for the families of free, triangular and cyclic groups considered in this paper with the group length. On the other hand, if we start with a non-optimal time estimate for $(p,q) = (2,4)$ ---as we do in Theorem B--- it is not difficult to keep tract of constants to show that only a $\log 3$ is lost in the extrapolation process. Note that Markovianity of the semigroup is implicitly used in the argument above, so that Theorem B can not be sharpened beyond $(p,q)=(2,4)$ for non conditionally negative lengths, at least using these ideas.}
\end{remark}

\renewcommand{\theequation}{B.1}
\addtocounter{equation}{-1}

\subsection*{B. The word length in $\Z_n$}

Let $\Z_n \sim \{ \exp (2 \pi i k/n) \, | \, 0 \le k < n \} \subset \mathbb{T}$ denote the finite cyclic group with $n$ elements. If we identify $\exp(2 \pi i k/n)$ with $k$ $(\mathrm{mod} \, n)$ as usual, our goal here is to show that the word length $\psi_n(k) = \min(k,n-k)$ is conditionally negative. We start with a trigonometric inequality. If $n \in 2 \Z_+$ and $j \in \Z$, we find
\begin{equation} \label{TrigIneq}
\Phi(n,j) \, := \, \frac{n}{2} + 2 \sum_{k=1}^{\frac{n}{2}-1} \big( \frac{n}{2} - k \big) \cos \big( \frac{2\pi kj}{n} \big) \, \ge \, 0.
\end{equation}
Indeed, let $\mathrm{M} = \frac{n}{2}-1$ and $z = \exp \big( \frac{2 \pi i j}{n} \big)$, so that
\begin{eqnarray*}
\Phi(n,j) & = & \mathrm{M}+1 + 2 \, \mathrm{Re} \Big( \sum_{k=1}^{\mathrm{M}} \big( \mathrm{M}+1 - k \big) z^k \Big) \\ & = & \mathrm{M}+1 + 2 \, \mathrm{Re} \Big( (\mathrm{M}+1) \frac{z-z^{\mathrm{M}+1}}{1-z} - \sum_{k=1}^{\mathrm{M}} k z^k \Big) \\ & = & \mathrm{Re} \Big( \mathrm{M}+1 + 2 (\mathrm{M}+1) \, \frac{z-z^{\mathrm{M}+1}}{1-z} - 2 z \frac{1 - (\mathrm{M}+1) z^{\mathrm{M}} + \mathrm{M}z^{\mathrm{M}+1}}{(1-z)^2} \Big) \\ [5pt] & = & \mathrm{Re} \Big( (\mathrm{M}+1) \frac{1+z}{1-z} + 2 \frac{z^{\mathrm{M}+2} - z}{(1-z)^2} \Big) \ = \ \mathrm{Re} \big( \mathrm{A} + \mathrm{B} \big) \ = \ \mathrm{Re} (\mathrm{B}).
\end{eqnarray*}
The last identity follows from $|z| =1$. On the other hand, it is easily checked that $z^{\mathrm{M}+2} = (-1)^j z$ so that $\mathrm{B} = 0$ for $j$ even. When $j$ is odd, we get $\mathrm{B} = -4z / (1-z)^2$ and $\mathrm{sgn} (\mathrm{Re} \mathrm{B}) \! = \! \mathrm{sgn} ( \mathrm{Re}(-z(1-\overline{z})^2)) \! = \! \mathrm{sgn}(1 - \cos(2\pi j/n)) \! = \! +1$. This proves \eqref{TrigIneq}.

Let us now prove that $\psi_n: \Z_n \to \Z_+$ is conditionally negative. We may assume without loss of generality that $n \in 4 \Z$, since for other values of $n$ we can always embed $j: k \in \Z_n \mapsto 4k \in \Z_{4n}$ and use that $$4 \psi_n \, = \, {\psi_{4n}}_{\mid_{j(\Z_n)}}.$$ Consider now a sequence $a_1, a_2, \ldots, a_n$ of complex numbers with $\sum_j a_j = 0$. In what follows we shall identify $a_j$ with $a_{j + n\Z}$ for obvious reasons. The goal is to show the non-positivity of the double sum $\sum \overline{a}_j a_k \psi_n(k-j)$. Let us begin by rewriting this sum in a more convenient way, which will allow us to reduce the problem to the case of $\R$-valued coefficients. We have    
\renewcommand{\theequation}{B.2}
\addtocounter{equation}{-1}
\begin{eqnarray} \label{RealCoeffcs}  
\lefteqn{\hskip-5pt \sum_{j,k=1}^n \overline{a}_j a_k \psi_n(k-j)} \\ \nonumber & = & \frac{n}{2} \sum_{j=1}^{\frac{n}{2}} |a_j + a_{j+ \frac{n}{2}}|^2 + \sum_{k=1}^{\frac{n}{2}-1} k \Big( \sum_{j=1}^n |a_j + a_{j+k}|^2 \Big) - \frac{n^2}{4} \sum_{j=1}^n |a_j|^2.
\end{eqnarray}
Namely, if we set $A_m := \sum_{j=1}^n \overline{a}_j a_{j+m}$, then we find $$\sum_{j,k=1}^n \overline{a}_j a_k \psi_n(k-j) \, = \, \sum_{m=1}^{n-1} \psi_n(m) A_m.$$ On the other hand, if $1 \le k \le \frac{n}{2}-1$ we obtain the following identities
\begin{itemize}
\item $\displaystyle A_{\frac{n}{2}} = 2 \sum_{j=1}^{\frac{n}{2}} \mathrm{Re} (\overline{a}_j a_{j+\frac{n}{2}}) = \sum_{j=1}^{\frac{n}{2}} |a_j + a_{j + \frac{n}{2}}|^2 - \sum_{j=1}^n |a_j|^2$,

\item $\displaystyle A_k + A_{n-k} = 2 \sum_{j=1}^n \mathrm{Re} (\overline{a}_j a_{j+k}) = \sum_{j=1}^n |a_j + a_{j+k}|^2 - 2 \sum_{j=1}^n |a_j|^2$.
\end{itemize}
Combining the identities above and rearranging terms, we deduce \eqref{RealCoeffcs}. According to it, we may assume in what follows that $a_j \in \R$ for all $1 \le j \le n$ since otherwise we may split the problem into real and imaginary parts by using $\sum_j \mathrm{Re}(a_j) = \sum_j \mathrm{Im}(a_j) = 0$ and the identity $|z|^2 = \mathrm{Re}(z)^2 + \mathrm{Im}(z)^2$. Let us consider the function $f: \R^n \to \R$ given by $$f(x) \, = \, \frac{n^2}{4} \sum_{j=1}^n x_j^2 - \frac{n}{2} \sum_{j=1}^{\frac{n}{2}} (x_j + x_{j+\frac{n}{2}})^2 - \sum_{k=1}^{\frac{n}{2}-1} k \Big( \sum_{j=1}^n (x_j + x_{j+k})^2 \Big),$$ where we still use the identifications $x_j = x_{j + n \Z}$. It suffices to show that $f \ge 0$ when restricted to the hyperplane $\Pi = \{\sum_j x_j = 0\}$. Moreover, since $f(\lambda x) = \lambda^2 f(x)$ the sign of $f$ in $\Pi$ is completely determined by the sign of $f$ in the compact set $\Omega = \mathbb{S}^{n-1} \cap \Pi$. Thus, our assertion will follow by showing that the absolute minimum of $f_{\mid_\Pi}$ (which exists since it coincides with the absolute minimum of $f_{\mid_\Omega}$) is not negative. Define $g(x) = \sum_j x_j$, according to the Lagrange multiplier method we have to solve the system $g=0$ and $\nabla f = \gamma \nabla g$. We claim that the solutions $x$ of this system satisfy
\begin{itemize}
\item[i)] $f(x) = 0$,

\item[ii)] $\nabla [f + \frac{n}{2} g^2] (x) = 0$,

\item[iii)] $\mathrm{Hess} [f + \frac{n}{2} g^2] (x) \ge 0$.
\end{itemize}  
In particular, all the solutions of the system are local minimums of $f_{\mid_\Pi} = (f + \frac{n}{2}g^2)_{\mid_\Pi}$. Let us prove the claim. It is clear that we have $\nabla g = (1,1,\ldots,1)$, while the partial derivarives of $f$ are given by $$-\frac12 \partial_jf(x) \, = \, \sum_{k=1}^{\frac{n}{2}} k (x_{j+k} + x_{j-k}) - \frac{n}{2} x_{j+\frac{n}{2}}.$$ In particular, setting $\mu = -\frac12 \gamma$ our system can be written as follows
\begin{itemize}
\item[(R)] $\displaystyle \sum_{j=1}^n x_j \, = \, 0$,

\item[($j$-th)] $\displaystyle \sum_{k=1}^{\frac{n}{2}} k (x_{j-k} + x_{j+k}) - \frac{n}{2} x_{j + \frac{n}{2}} \, = \, \mu$ \hskip5pt for \hskip5pt $1 \le j \le n$.
\end{itemize}
Now, operating with these equations we obtain the following identities
\begin{itemize}
\item $\displaystyle \sum_{j=1}^n (\mbox{$j$-th}) - \frac{n^2}{4} (\mbox{R}) \Leftrightarrow n \mu = \Big( 2 \big( \sum_{k=1}^{\frac{n}{2}} k \big) - \frac{n}{2} - \frac{n^2}{4} \Big) \Big( \sum_{j=1}^n x_j \Big) = 0$.

\item $\displaystyle \frac12 \Big[ \mbox{(R)} + ((\mbox{$j +\frac{n}{2}+1$)-th}) - ((\mbox{$j+\frac{n}{2}$)-th}) \Big] \Leftrightarrow \sum_{s=j+1}^{j+\frac{n}{2}} x_s = 0$ \ for \ $1 \le j \le n$.
\end{itemize}
The first identity gives $\mu = \gamma = 0$, while the second one is equivalent to $x_j = x_{j+\frac{n}{2}}$ for all $1 \le j \le n$ combined with (R). Now we are ready to prove the first assertion i) of our claim. Namely, since $x_j = x_{j+\frac{n}{2}}$ and $\sum_j x_j = 0$, we only have $\frac{n}{2}-1$ variables $x_1, x_2, \ldots, x_{\frac{n}{2}-1}$ and we may evaluate $f$ at the points satisfying these restrictions as follows
\begin{eqnarray*}
f(x) & = & \big( \frac{n^2}{2} - 2n \big) \sum_{j=1}^{\frac{n}{2}} x_j^2 - 2 \sum_{k=1}^{\frac{n}{2}-1} k \sum_{j=1}^{\frac{n}{2}} (x_j + x_{j+k})^2 \\ & = & \big( \frac{n^2}{2} - 2n - 4 \sum_{k=1}^{\frac{n}{2}-1} k \big) \sum_{j=1}^{\frac{n}{2}} x_j^2 - 4 \sum_{k=1}^{\frac{n}{2}-1} k \sum_{j=1}^{\frac{n}{2}} x_j x_{j+k} \\ & = & n \Big( \sum_{j=1}^{\frac{n}{2}-1} x_j \Big)^2 -n \sum_{j=1}^{\frac{n}{2}-1} x_j^2 - 4 \sum_{\begin{subarray}{c} j,k=1 \\ j+k \neq \frac{n}{2} \end{subarray}}^{\frac{n}{2}-1} k x_j x_{j+k} \ \equiv \ 0.
\end{eqnarray*}
This proves the first assertion, while $\nabla [f + \frac{n}{2} g^2] = (\gamma + n \sum_j x_j ) \nabla g = 0$ by (R) and identity $\gamma = 0$ proved above. It remains to justify assertion iii), for which we start by noticing that $\frac12 \partial_{jk} [f + \frac{n}{2} g^2] (x) = \frac{n}{2} - \psi_n(j-k)$. Using the unitary matrices $\lambda(s) = \sum_{j \in \Z_n} e_{j,j+s}$ in $\mathcal{B}(\ell_2(\Z_2))$ and the $*$-homomorphism $$\lambda(s) \in \mathcal{L}(\Z_n) \mapsto \exp(2 \pi i \cdot/n) \in L_\infty(\Z_n),$$ we may write the Hessian matrix of $f + \frac{n}{2} g^2$ as follows 
\begin{eqnarray*}
\frac12 \mathrm{Hess} [f + \frac{n}{2} g^2] (x) & = & \frac{n}{2} \lambda(0) + \sum_{k=1}^{\frac{n}{2}-1} \big( \frac{n}{2} - k \big) \big( \lambda(k) + \lambda(n-k) \big) \\ & \sim & \frac{n}{2} + 2 \sum_{k=1}^{\frac{n}{2}-1} \big( \frac{n}{2} -k \big) \cos \big( \frac{2 \pi k \, \cdot}{n} \big) \ = \ \Phi(n, \cdot).
\end{eqnarray*} 
By the given $*$-homomorphism, the positivity of the Hessian matrix of $f + \frac{n}{2} g^2$ at every $x$ reduces to the trigonometric inequality which we proved in \eqref{TrigIneq}. \fin

\subsection*{C. Numerical analysis}

In this appendix we justify all the numerical estimates for free groups, triangular groups and cyclic groups which we used in the proofs of Theorems A1-A3. As we did above, we label these estimates by $\alpha, \beta, \gamma, \delta, \varepsilon$. Let us start with the free group estimates. 

\noindent {\bf C1. Estimates for free groups.} We will give in this section detailed proofs of all the numerical estimates for free groups, which will serve us to omit some details in the case of triangular groups below.

\noindent {\bf Estimates $(\alpha)$ for $\F_n$.} Assertion $\alpha$i) follows from the identity $s_2(r) = \sum_{k \ge 1} r^{2k} \alpha_k$. According to Proposition \ref{s23}, the second assertion $\alpha$ii) follows from the inequalities $\frac12 M(k_1+k_2-2m, k_1,k_2) \le \frac32 M(k_1,k_2)= \frac{3}{s(s-1)}M(k_1,k_2,\underline{0})$. 
To prove $\alpha$iii) and $\alpha$iv) we use Proposition \ref{s23} again, so that we have to estimate 
\begin{eqnarray*}
A(k,r) & = & \sum_{\begin{subarray}{c} \ell_2 \ge \ell_3 \ge 1  \\ \ell_2 + \ell_3 = k \end{subarray}} \frac12 M(k, \ell_2, \ell_3) \, r^{2k} \\ & + & \sum_{\ell_2 \ge \ell_3 \ge 1} \sum_{\begin{subarray}{c} m=1 \\ \ell_2 + \ell_3 -2m = k \end{subarray}}^{\lfloor \ell_3/2 \rfloor} \frac12 (2n-2) (2n-1)^{m-1} M(k,\ell_2,\ell_3) \, r^{\ell_2+\ell_3+k}.
\end{eqnarray*}
Let us consider the sets $$\Omega(k) \, = \, \Big\{ (\ell_2, \ell_3, m) \, : \ \ell_2 + \ell_3 - 2m = k, \  0 \le m \le \frac{\ell_3}{2} \ \mbox{and} \ \ell_2 \ge \ell_3 \ge 1 \Big\}.$$ It is easy to check that we have 
\begin{eqnarray*}
\Omega(1) & = & \emptyset, \\ \Omega(2) & = & \big\{ (1,1,0), (2,2,1) \big\}, \\ \Omega(3) & = & \big\{ (2,1,0), (3,2,1) \big\}, \\ \Omega(4) & = & \big\{ (3,1,0), (2,2,0), (4,2,1), (3,3,1), (4,4,2) \big\}. 
\end{eqnarray*}
By simple computations we may show that 
\begin{eqnarray*}
\big( A(1,r), A(2,r), A(3,r) \big) & = & \Big( 0, \frac32 r^4+(n-1)r^6, 3 r^6 + 3 (n-1)r^8 \Big), \\ A(4,r) & = & \frac92 r^8 + 6 (n-1)r^{10} + (n-1)(2n-1) r^{12}.
\end{eqnarray*}
When $k \ge 5 \delta_{n=2} + 3 \delta_{n \ge 3}$ and $(2n-1)r^2 < 1$, we use $M(k,\ell_2,\ell_3) \le 6$ to obtain 
\begin{eqnarray*}
A(k,r) & \le & \frac32 k r^{2k} \ + \ \frac{6(n-1)}{2n-1} \, \Big( \frac{r}{\sqrt{2n-1}} \Big)^{k} \!\!\! \sum_{\begin{subarray}{c} k \ge \ell_2 \ge \ell_3 \ge 1 \\ k < \ell_2 + \ell_3 \\ \ell_2 + \ell_3 \equiv k \ \mathrm{mod} \ 2 \end{subarray}} (\sqrt{2n-1}r)^{\ell_2+\ell_3} \\ & \le & \frac32 k r^{2k} \ + \ \frac{6(n-1)}{2n-1} \, \Big( \frac{r}{\sqrt{2n-1}} \Big)^{k} \sum_{\ell_3=1}^{k} \sum_{\begin{subarray}{c} \ell_2 = \ell_3 \vee (k - \ell_3 + 1) \\ \ell_2 \equiv k - \ell_3 \ \mathrm{mod} \ 2 \end{subarray}}^{k} (\sqrt{2n-1}r)^{\ell_2+\ell_3} \\ & \le & \frac32 k r^{2k} \ + \ \frac{6(n-1)}{2n-1} \, r^{2k} \sum_{\ell_3=1}^{k/2} \sum_{\begin{subarray}{c} j=1 \\ j \ \mathrm{even} \end{subarray}}^{\ell_3} (\sqrt{2n-1}r)^{j} \\ & + & \frac{6(n-1)}{2n-1} \, r^{2k} \sum_{\ell_3 = k/2}^{k} \, \sum_{\begin{subarray}{c} j=2\ell_3-k \\ j \ \mathrm{even} \end{subarray}}^{\ell_3} (\sqrt{2n-1}r)^{j} \\ [5pt] & \le & \frac32 k r^{2k} \ + \ \frac{6(n-1)}{2n-1} \, r^{2k} \hskip5pt \sum_{\ell_3=1}^{k/2} \frac{(2n-
 1)r^2}{1-(2n-1)r^2} \\ & + & \frac{6(n-1)}{2n-1} \, r^{2k} \sum_{\ell_3 = k/2}^{\infty} \frac{(\sqrt{2n-1}r)^{2\ell_3-k}}{1-(2n-1)r^2} \\ [5pt] & = & \frac32 k r^{2k} \ + \ \frac{3(n-1) kr^{2k+2}}{1-(2n-1)r^2} \ + \ \frac{6(n-1)}{2n-1}\frac{ r^{2k}}{(1-(2n-1)r^2)^2}.
\end{eqnarray*}
If $r=\frac{1}{\sqrt{2n-1}}$, we obtain
\begin{eqnarray*}
A(k,r) \!\!\! & \le & \!\!\! \frac32 k r^{2k} \\ & + & \frac{6(n-1)}{2n-1} \, r^{2k} \Big( \sum_{\ell_3=1}^{k/2} \sum_{\begin{subarray}{c} j=1 \\ j \ \mathrm{even} \end{subarray}}^{\ell_3} (\sqrt{2n-1}r)^{j} + \sum_{\ell_3 = k/2}^{k} \, \sum_{\begin{subarray}{c} j=2\ell_3-k \\ j \ \mathrm{even} \end{subarray}}^{\ell_3} (\sqrt{2n-1}r)^{j} \Big) \\ \!\!\! & \le & \!\!\! \frac32 k (2n-1)^{-k} \hskip1pt + \ 6(n-1)(2n-1)^{-k-1} \Big( \sum_{\ell_3=1}^{k/2} \frac{\ell_3}{2} + \sum_{\ell_3=k/2}^{k} \frac{k - \ell_3 + 2}{2} \Big) \\ [8pt] \!\!\! & \le & \!\!\! \frac32 k (2n-1)^{-k} \hskip1pt + \ 6(n-1)(2n-1)^{-k-1} \Big( \sum_{\ell_3=1}^{k/2} \ell_3 + \frac{k}{2} + 1 \Big) \\ [5pt]  \!\!\! & \le & \!\!\!  \  \frac14 (2n-1)^{-k-1} \Big(3(n-1)k^2+(30n-24)k+24(n-1)\Big) . \hskip60pt \square
\end{eqnarray*}

\noindent {\bf Estimates $(\beta)$ for $\F_2$.} We claim that
\begin{enumerate}
\item[a)] If $\k \in L_m(s)$ and $\underline{\ell} \in \mathrm{Adm}_{2m}$ with $\ell_{\xi_j} = k_j$ for some $\x \in C_{2m}$
$$\hskip20pt M(\l) \, \le \, \frac{(2m)!}{m!} M(k_1, k_2, \ldots, k_m) \, = \, \frac{(2m)!(s-m)!}{s!} M(\k).$$

\item[b)] If $\k \in L_m(s)$ and $\x \in C_{2m}$ we have 
\begin{eqnarray*}
\lefteqn{\hskip-20pt \sum_{\begin{subarray}{c} \l \in \mathrm{Adm}_{2m} \setminus B_{2m} \\ \ell_{\xi_j} = k_j \end{subarray}} \Big( \prod_{j=2}^{m} N_{\ell_{\xi_j^\star}} r^{\ell_{\xi_j^\star}} \Big) r^{\ell_{\xi_1^\star}} \ =: \ \Theta_m(\k,\x,r)} \\ \hskip40pt & \le & 
\left\{\begin{array}{ll}
2^{m-2}
& \mbox{ if }  r =\frac{1}{3}, \\ [6pt]
\displaystyle\frac{(3r)^{|\k|}}{1-r}\Big(\displaystyle\frac{4r}{3r-1}\Big)^{m-1}
& \mbox{ if } \frac{1}{3}< r < \frac{1}{\sqrt{3}}, \\ [9pt]
\displaystyle\frac{4|\k|}{3-\sqrt{3}}& \mbox{ if }  r =\frac{1}{\sqrt{3}} \mbox{ and } s=m=2.
\end{array}\right.
\end{eqnarray*}

\item[c)] We have $|C_{2m}| = 2^m$.
\end{enumerate}
A moment of thought shows that $\beta$i) follows from Proposition \ref{Lambda-est} and the claim above. Set $j(k_1, k_2, \ldots, k_m) = (j_1, j_2, \ldots, j_m)$ and $j(\l) = (j_1', j_2', \ldots, j_{2m}')$. To prove a), we observe that $\sum_{i \ge i_0} j_i$ gives the number of blocks in $(k_1, k_2, \ldots, k_m)$ of size greater than or equal to $i_0$. Similarly, the sum $\sum_{i \ge i_0} j_i'$ gives the number of blocks in $(\ell_1, \ell_2, \ldots, \ell_{2m})$ of size greater than or equal to $i_0$. Since $\ell_{\xi_j} = k_j$ we find that $\sum_{i \ge i_0} j_i \, \le \, \sum_{i \ge i_0} j_i'.$ This immediately gives $$\prod_{i=1}^m (i!)^{j_i} = \prod_{i=1}^m i^{j_i + \ldots + j_m} \le \prod_{i=1}^{2m} i^{j_i' + \ldots + j_{2m}'} = \prod_{i=1}^{2m} (i!)^{j_i'},$$ which implies the first inequality in a). The identity is clear by the combinatorial interpretation of $M(\k)$. 
We now turn to the proof of b). If $r=\frac{1}{3}$, then $N_j r^j = \frac43 (3r)^j = 4/3$ for all $j \ge 1$ from which we deduce  
\begin{equation*}
\Theta_m(\k,\x,r) \, \le \, \Big(\frac{4}{3}\Big)^{m-1} \sum_{\ell_{\xi_1^\star} \ge \cdots \ge \ell_{\xi_m^\star} \ge 1} r^{\ell_{\xi_1^\star}} \, = \, \Big(\frac{4}{3}\Big)^{m-1} \frac{r}{(1-r)^m} \, = \, 2^{m-2}. 
\end{equation*}
Indeed, by induction it is easy to see that for $r<1$
$$\sum_{\ell_1 \ge \cdots \ge \ell_{m} \ge 1} r^{\ell_{1}} = \sum_{\ell_2 \ge \cdots \ell_{m} \ge 1} \Big(\sum_{\ell_1 \ge \ell_2} r^{\ell_{1}}\Big) = \frac{1}{1-r}\sum_{\ell_2 \ge \cdots \ell_{m} \ge 1} r^{\ell_{2}}= \frac{r}{(1-r)^m}.$$
If $\frac{1}{3}< r < \frac{1}{\sqrt{3}}$, we use $$\ell_{\xi_j} = k_j \mbox{ for some } \x \in C_{2m} \ \Rightarrow \ \max \big\{ k_{j+1},1 \big\} \le \ell_{\xi_j^\star} \le k_{j-1} \mbox{ for } 2 \le j \le m.$$ Thus, we obtain the desired estimate as follows 
\begin{eqnarray*}
\Theta_m(\k,\x,r) & \le &  \Big(\frac{4}{3}\Big)^{m-1} \Big(\sum_{\ell_{\xi_1^\star} \geq 0} r^{\ell_{\xi_1^\star}}\Big) 
\prod_{j=2}^{m} \Big(\sum_{\ell_{\xi_j^\star}=\max\{k_{j+1},1\}}^{k_{j-1}} (3r)^{\ell_{\xi_j^\star}}\Big) \\
&\le &  \Big(\frac{4}{3}\Big)^{m-1} \frac{1}{1-r} \frac{(3r)^{\sum_{j=1}^{m-1}k_j +m-1}}{(3r-1)^{m-1}} \ \le \ \frac{(3r)^{|\k|}}{1-r}\Big(\frac{4r}{3r-1}\Big)^{m-1}.
\end{eqnarray*}
If  $r =\frac{1}{\sqrt{3}}$ and $s=m=2$, the estimate becomes
$$\Theta_4(\k,\x,r) \, \le \, \frac{4}{3}\sum_{\ell_{\xi_2^\star}=1}^{k_1} (\sqrt{3})^{\ell_{\xi_2^\star}}\sum_{\ell_{\xi_1^\star}\geq \ell_{\xi_2^\star}} \Big(\frac{1}{\sqrt{3}}\Big)^{\ell_{\xi_1^\star}} \, = \, \frac{4}{3}\frac{k_1}{1-\frac{1}{\sqrt{3}}}\le \displaystyle\frac{4|\k|}{3-\sqrt{3}}.$$ Since c) is clear, this completes the proof of the even case $\beta$i). The estimates $\beta$ii) and $\beta$iii) are proved similarly, by estimating the right hand side of the inequality ii) of Proposition \ref{Lambda-est} when $\k \in L_m(s)$ and $\k \in L_{m-1}(s)$ respectively. The claim to be proved becomes in the odd case 
\begin{enumerate}
\item[a$'$)] If $\k \in L_m(s)\cup L_{m-1}(s)$ and $\l \in \mathrm{Adm}'_{2m}$ with $\ell_{\xi_j} = k_j$ for some $\x \in C_{2m}$ 
$$\hskip20pt M(\l) \, \le \, \left\{\begin{array}{ll}
\displaystyle\frac{(2m)!(s-m)!}{s!} M(\k)& \mbox{ if } \k \in L_{m}(s),\\ [5pt] \displaystyle\frac{(2m)!(s-m+1)!}{s!} M(\k) & \mbox{ if } \k \in L_{m-1}(s). \end{array}\right.$$

\item[b$'$)] If $\k \in L_m(s)\cup L_{m-1}(s)$ and $\x \in C_{2m}$ we have
\begin{eqnarray*}
\hskip20pt \lefteqn{\hskip-20pt \sum_{\begin{subarray}{c} \l \in \mathrm{Adm}'_{2m} \setminus B'_{2m} \\ \ell_{\xi_j} = k_j \end{subarray}} \Big( \prod_{j=2}^{m} N_{\ell_{\xi_j^\star}} r^{\ell_{\xi_j^\star}} \Big) r^{\ell_{\xi_1^\star}}
\ =: \ \Theta_m'(\k,\x,r)} \\ \hskip50pt & \le & \left\{\begin{array}{ll} 2^{m-3} & \mbox{ if }  r =\frac{1}{3} \mbox{ and } \k \in L_{m}(s), \\ [8pt] 2^{m-2} & \mbox{ if }  r =\frac{1}{3} \mbox{ and } \k \in L_{m-1}(s), \\ [3pt] \displaystyle\frac{(3r)^{|\k|}}{1-r}\Big(\frac{4r}{3r-1}\Big)^{m-2} & \mbox{ if } \frac{1}{3}< r < \frac{1}{\sqrt{3}} \mbox{ and } \k \in L_{m}(s), \\ [3pt] \displaystyle\frac{(3r)^{|\k|}}{1-r}\Big(\displaystyle\frac{4r}{3r-1}\Big)^{m-1} & \mbox{ if } \frac{1}{3}< r < \frac{1}{\sqrt{3}} \mbox{ and } \k \in L_{m-1}(s). 
\end{array}\right.
\end{eqnarray*}

\item[c$'$)] We have $|C'_{2m}| = 2^{m-1}$.
\end{enumerate}
To see a$'$), it suffices to note that if $\k \in L_{m-1}(s)$, then $$M(\k)  =  \binom{s}{m-1} M(k_1, k_2, \ldots, k_{m-1}) = \frac{1}{m} \binom{s}{m-1} M(k_1, k_2, \ldots, k_{m})$$ and argue as for a). For the assertions b$'$) and c$'$) we use the fact that $\xi_m$ is completely determined for $\k \in L_m(s) \cup L_{m-1}(s)$ and $\l \in \mathrm{Adm}'_{2m}$ satisfying $\ell_{\xi_j} = k_j$ for some $\x \in C_{2m}$. More precisely, if $\k \in L_m(s)$ then $\xi_m=2m-1$, hence $\xi^\star_m=2m$ and $\ell_{\xi^\star_m}=0$. If $\k \in L_{m-1}(s)$, then $\xi_m=2m$. This proves $\beta$ii) and $\beta$iii). \fin

\noindent {\bf Estimates $(\gamma)$ for $\F_n$.} All the assertions follow directly from Proposition \ref{method1}. \fin

\noindent {\bf Estimates $(\delta)$ for $\F_n$.} Since $\x \in C_{2m} \Rightarrow N_{\ell_{\xi_j^\star}} \le N_{\ell_{2j-1}}$, Proposition \ref{Lambda-est} yields
\begin{eqnarray*}
\alpha_0^{s-m} \sum_{\lambda \in \Lambda_{2m}[\k]} \gamma_\lambda \nu_\lambda(r) &\le & \frac{1}{2^m} \Big[ \sum_{\x \in C_{2m}} \sum_{\begin{subarray}{c} \l \in \mathrm{Adm}_{2m} \setminus B_{2m} \\ \ell_{\xi_j} = k_j \end{subarray}}  M(\l)\Big( \prod_{j=2}^{m} N_{\ell_{\xi_j^\star}}\Big)r^{|\l|} \Big] \alpha_{\k} \\ & \le & \frac{1}{2^m} \Big[ \sum_{\l \in A_{2m}[\k]} \sum_{\begin{subarray}{c} \x \in C_{2m} \\ \ell_{\xi_j} = k_j \end{subarray}}  M(\l)r^{|\l|} \Big] \Big( \prod_{j=2}^{m} N_{\ell_{2j-1}}\Big)  \alpha_{\k}. 
\end{eqnarray*}
Then $\delta$i) follows from the fact that $$\big| \big\{ \x \in C_{2m} \, : \, \ell_{\xi_j} = k_j \ \mbox{for all} \ j=1,2,\ldots,m \big\} \big|=2^{P(\l)}$$ for $\k \in L_m(s)$ and $\l \in A_{2m}[\k]$. The estimates $\delta$ii) and $\delta$iii) are proved similarly. We just need to observe that if $\k \in L_m(s)$ and $\l \in \mathrm{Adm}'_{2m}$ satisfy $\ell_{\xi_j} = k_j$ for some $\x \in C_{2m}$, then $\xi^\star_m=2m$ and $\ell_{\xi^\star_m}=0$. 
Hence $N_{\ell_{\xi^\star_m}}=1$ and the product in the right hand side of the inequality ii) of Proposition \ref{Lambda-est} runs over $2\leq j \leq m-1$. \fin

\renewcommand{\theequation}{C1.1}
\addtocounter{equation}{-1}

\noindent {\bf Estimates $(\varepsilon)$ for $\F_2$.} To estimate the sums $\varepsilon$4), $\varepsilon$5), $\varepsilon$6), $\varepsilon$7), $\varepsilon$8), $\varepsilon$10) and $\varepsilon$11) we can be slightly less precise than for the other terms (see below) and use a general formula.  The estimates of the $s_4$-sums are $\varepsilon$4), $\varepsilon$5), $\varepsilon$6), $\varepsilon$7), $\varepsilon$8) and follow from the formula 
$$s_4[\l](r) \le \frac{1}{2} M(\l) \big( N_{\ell_3}\alpha_{\ell_1}\alpha_{\ell_4}+N_{\ell_4}\alpha_{\ell_2}\alpha_{\ell_3} \big) r^{|\ell|}.$$
To prove it, we write as in the proof of Lemma \ref{EvenLemma} $$s_4[\l](r) \ \le \ \sum_{\d \sim \l}^{\null} \sum_{\substack{\ibar \\ w_{d_1}(i_1) \cdots w_{d_{4}}(i_{4}) = e}} \Big( \prod_{j=1}^{4} a_{\ell_j}(i_{\sigma_{\d}^{-1}(j)}) \Big) r^{|\l|}.$$ Now, we estimate the sum in the right hand side as follows
\begin{eqnarray*}
\sum_{\d, \ibar} \prod_{j=1}^4 a_{\ell_j}(i_{\sigma_{\d}^{-1}(j)}) & \le & \frac{1}{2} \sum_{\d, \ibar} a_{\ell_1}(i_{\sigma_{\d}^{-1}(1)})^2 a_{\ell_4}(i_{\sigma_{\d}^{-1}(4)})^2 + a_{\ell_2}(i_{\sigma_{\d}^{-1}(2)})^2 a_{\ell_3}(i_{\sigma_{\d}^{-1}(3)})^2  \\ & \le & \frac{1}{2} \sum_{\d \sim \l} \sum_{i_{\sigma_{\d}^{-1}(1)},i_{\sigma_{\d}^{-1}(3)}, i_{\sigma_{\d}^{-1}(4)}} a_{\ell_1}(i_{\sigma_{\d}^{-1}(1)})^2 a_{\ell_4}(i_{\sigma_{\d}^{-1}(4)})^2 \\ & + & \frac12 \sum_{\d \sim \l} \sum_{i_{\sigma_{\d}^{-1}(2)},i_{\sigma_{\d}^{-1}(3)}, i_{\sigma_{\d}^{-1}(4)}} a_{\ell_2}(i_{\sigma_{\d}^{-1}(2)})^2 a_{\ell_3}(i_{\sigma_{\d}^{-1}(3)})^2 \\
& = & \frac{1}{2} M(\l) \big( N_{\ell_3} \alpha_{\ell_1} \alpha_{\ell_4} + N_{\ell_4} \alpha_{\ell_2} \alpha_{\ell_3} \big).
\end{eqnarray*}  
For the $s_5$-sums $\varepsilon$10) and $\varepsilon$11), we use respectively the formulas 
\begin{eqnarray} \label{s5Est}
\widehat{f}(e)s_5[\l](r) & \le & \frac{1}{2} M(\l) \big( N_{\ell_4}N_{\ell_5}\alpha_0 \alpha_{\ell_1}\alpha_{\ell_2}+N_{\ell_2}\alpha_{\ell_3}\alpha_{\ell_4}\alpha_{\ell_5} \big) r^{|\ell|}, 
\end{eqnarray}
\vskip-10pt
\renewcommand{\theequation}{C1.2}
\addtocounter{equation}{-1}
\begin{eqnarray}
\widehat{f}(e)s_5[\l](r) & \le & \frac{1}{2} M(\l) \big( N_{\ell_3}N_{\ell_4}\alpha_0 \alpha_{\ell_1}\alpha_{\ell_5}+N_{\ell_5}\alpha_{\ell_2}\alpha_{\ell_3}\alpha_{\ell_4}\big) r^{|\ell|}.
\end{eqnarray}
The proofs are similar to our estimate for $s_4[\l](r)$, we omit the details. Let us now prove the estimates for the remaining sums. Since we want to estimate sums of small size, we can use computer assistance to compute exactly these sums, and then complete squares by hand. Namely, it is quite simple to teach a computer to multiply words in a free group and determine whether a product of them equals $e$. This allows us to identify the summation index of the sums $s_u[\l](r)$ for small parameters $(u,\l)$ as it is the case for the estimates $(\varepsilon)$ we are about to prove. Once the summation index is given by the computer, we obtain an identity in terms of $a_k(i)$'s which can be estimated by $\alpha_k$'s completing squares. In other words, the first identity in our estimates below for the sums $\varepsilon$1), $\varepsilon$2) and $\varepsilon9)$ ---the ones where we need to be more careful, since they contribute to the most pathological term $\alpha_0^{s-2}\alpha_1^2$--- is given by the computer. This will serve as an illustration of how to proceed for the other terms, for which we will omit some details 
\begin{eqnarray*}
\varepsilon1) & s_4[(1,1,1,1)](r) & = \ \Big[ 6 \big( a_1(1)^4 +a_1(2)^4 \big) +16a_1(1)^2a_1(2)^2 \Big] r^4 \\  
& & = \ \frac{3}{2}r^4\alpha_1^2+4r^4 a_1(1)^2a_1(2)^2 \ \leq \ 2r^4\alpha_1^2.\\ [5pt]
\varepsilon2) & s_4[(3,1,1,1)](r) & = \ 8 \Big[a_1(1)^3a_3(i_1)+a_1(2)^3a_3(i_2) \Big] r^6 \\ & & + \ 8 \Big[a_1(1)^2a_1(2)\sum_{j=3}^{10} a_3(i_j)+a_1(1)a_1(2)^2\sum_{j=11}^{18} a_3(i_j)\Big]r^6 \\ & & \le \ 36 r^6 \Big[ a_1(1)^4+a_1(2)^4 \Big] +r^6 \alpha_1\alpha_3 \ \leq \ 9r^6\alpha_1^2+r^6 \alpha_1\alpha_3.
\end{eqnarray*}
The estimate for $\varepsilon$9) is a bit more involved, since the goal now is to avoid the most pathological term $\alpha_0^{s-2} \alpha_1^2$ 
\begin{eqnarray*}
\lefteqn{\varepsilon9) \quad \widehat{f}(e) s_5[(2,1,1,1,1)](r)} \\ & = & a_0 \Big[ 50 \big( a_1(1)^3a_1(2)+a_1(1)a_1(2)^3 \big) \Big(\sum_{j=3}^6 a_2(i_j) \Big) \\ & + & 60a_1(1)^2a_1(2)^2 \big( a_2(i_1)+a_2(i_2) \big) + 40 \big( a_1(1)^4a_2(i_1)+a_1(2)^4a_2(i_2) \big) \Big] r^6 \\ [6pt] & \le & \Big[20a_1(1)^2(a_1(1)^4+a_0^2a_2(i_1)^2) +20a_1(2)^2(a_1(2)^4+a_0^2a_2(i_2)^2) \\
[3pt] & + & 25 \Big( \big(a_1(1)^2+a_1(2)^2 \big) \Big[ 4a_1(1)^2a_1(2)^2+a_0^2 \Big(\sum_{j=3}^6 a_2(i_j)^2\Big) \Big] \Big) \\
&+ & 30a_1(2)^2(a_1(1)^4+a_0^2a_2(i_1)^2) +30a_1(1)^2(a_1(2)^4+a_0^2a_2(i_2)^2)
\Big]r^6 \\ & \leq & \frac{15}{2}r^6\alpha_0\alpha_1\alpha_2 +\frac{65}{12}r^6\alpha_1^3.
\end{eqnarray*}
For the sums $\varepsilon$3), $\varepsilon$12), $\varepsilon$13), $\varepsilon$14) we just list the computer outcome. Showing that these expressions are bounded above by the given terms in the corresponding $(\varepsilon)$-estimates follows by completing squares   
\begin{align*}
\varepsilon3)\quad & s_4[(2,2,1,1)](r) \ = \ 8 \Big[ a_1(1)^2a_2(i_1)^2 + a_1(2)^2a_2(i_2)^2 \Big] r^6 \\ 
& \hskip28pt + \ 8 \Big[ 2 a_1(1)a_1(2) \big( a_2(i_1)+a_2(i_2) \big) \Big(\sum_{j=3}^6a_2(i_j)\Big) \Big] r^6 \\ & \hskip28pt + \  8 \Big[ \big( a_1(1)^2+a_1(2)^2 \big) \Big(2 \sum_{1 \le i \le 6} a_2(i)^2+\sum_{3\le j<j'\le 6} a_2(i_j)a_2(i_{j'})\Big) \Big] r^6. \end{align*}
\begin{align*}
\varepsilon12) \quad & s_6[(1,1,1,1,1,1)](r) \\ & \hskip55pt = \ \Big[ 20 \big( a_1(1)^6 +a_1(2)^6 \big) +96 \big( a_1(1)^4a_1(2)^2+a_1(1)^2a_1(2)^4 \big) \Big]r^6.
\end{align*}
\begin{align*}
\varepsilon13) \quad& \widehat{f}(e) s_7[(2,1,1,1,1,1,1)](r) \\
= & \ a_0 \Big[ a_2(i_1) \Big( 420 a_1(1)^2a_1(2)^4+728a_1(1)^4a_1(2)^2+210 a_1(1)^6 \Big) \Big] r^8 \\
+ & \ a_0 \Big[ a_2(i_2) \Big(420 a_1(1)^4a_1(2)^2 + 728a_1(1)^2a_1(2)^4 + 210 a_1(2)^6 \Big) \Big] r^8 \\
+ & \ a_0 \Big[ \Big( \sum_{3 \le j \le 6} a_2(i_j) \Big) \Big( 742a_1(1)^3a_1(2)^3 + 308 \big( a_1(1)^5a_1(2) + a_1(1)a_1(2)^5 \big) \Big)\Big] r^8.
\end{align*}
\begin{eqnarray*}
\lefteqn{\varepsilon14) \quad s_8[(1,1,1,1,1,1,1,1)](r)} \\ 
& \hskip-10pt = \ \Big[ 70 \big( a_1(1)^8 + a_1(2)^8 \big) + 512 \big( a_1(1)^6a_1(2)^2 + a_1(1)^2a_1(2)^6 \big) + 928 a_1(1)^4a_1(2)^4 \Big] r^8.
\end{eqnarray*}
\fin

\noindent {\bf Estimates $(\varepsilon')$ for $\F_n$.} For $n\geq 3$ and $\underline{j}, \underline{j}' \in \N^n$ we introduce the notation
$$\omega_n(\underline{j},\underline{j}')
=(\underbrace{c_1,\cdots,c_1}_{j_1},\underbrace{c_1^{-1},\cdots, c_1^{-1}}_{j'_1}, \cdots, 
\underbrace{c_n,\cdots,c_n}_{j_n},\underbrace{c_n^{-1},\cdots, c_n^{-1}}_{j'_n}) \in \F_n^{|\underline{j}|+|\underline{j}'|}.$$
We start with the proof of $\varepsilon'1)$. Let $2\leq m \leq s$. 
To estimate $s_{2m}[(1,\ldots,1)](r)$, we need to look at the words $g_1,\ldots,g_{2m}$ of length $1$ satisfying $g_1\cdots g_{2m}=e$. 
Since this clearly implies that $|\{j \, : \, g_j=c_i\}|=|\{j \, : \, g_j=c_i^{-1}\}|$ for all $1\leq i \leq n$, by the multinomial theorem we obtain
\begin{eqnarray*}
s_{2m}[(1,\ldots,1)](r) 
&=& \!\!\! \sum_{\begin{subarray}{c} \underline{j} \in \N^n \\ |\underline{j}|=m \ \null \end{subarray}} A_{n,m}(\underline{j}) a_1(1)^{2j_1}\cdots a_1(n)^{2j_n} r^{2m}\\
& \leq & \!\!\! A_{n,m} \Big[\sum_{\begin{subarray}{c}\underline{j} \in \N^n \\ |\underline{j}|=m \ \null \end{subarray}} \binom{m}{j_1,\cdots,j_n}\prod_{i=1}^n a_1(i)^{2j_i}\Big] r^{2m} 
 =  A_{n,m} \Big(\frac{\alpha_1}{2}\Big)^m r^{2m},
\end{eqnarray*}
where 
\begin{eqnarray*}
 A_{n,m}(\underline{j}) &=& \Big|\Big\{(g_1,\cdots,g_{2m}) \in \F_n^{2m}\, : \, g_1\ldots g_{2m}=e, (g_i)_{1\leq i \leq 2m} \sim \omega_n(\underline{j},\underline{j})\Big\}\Big|, \\
 A_{n,m} &=& \max \left\{ \frac{A_{n,m}(\underline{j})}{\binom{m}{j_1,\cdots,j_n}} \; : \ \underline{j} \in \N^n, \ |\underline{j}|=m\right\}.
\end{eqnarray*}
We claim that 
\renewcommand{\theequation}{C1.3}
\addtocounter{equation}{-1}
\begin{equation}\label{Anm}
 A_{n,m} \leq 2^m \frac{(2m)!}{m! (m+1)!} \quad ( = \mbox{ if } n\ge m).
\end{equation}
This completes the proof of $\varepsilon'1)$. 
For $\varepsilon'2)$, we consider $2\leq m \leq s-1$ and words $g_1,\ldots,g_{2m+1}$ with $(|g_i|)_{1\leq i \leq 2m+1}\sim (2,1,\cdots,1)$ satisfying that $g_1 g_2\cdots g_{2m+1}=e$. 
Similarly, we write
\begin{eqnarray*}
\lefteqn{\hskip-15pt \widehat{f}(e) s_{2m+1}[(2,1,\cdots,1)](r)} \\ [5pt]
&=& 2(2m+1)r^{2m+2}\Big[\sum_{k=1}^na_0\widehat{f}(c_k^2)\sum_{\begin{subarray}{c}\underline{j} \in \N^n \\ |\underline{j}|=m-1 \ \null  \end{subarray}} 
B_{n,m}(\underline{j},k,k) \prod_{i=1}^n a_1(i)^{2j_i} a_1(k)^2 \\
&+& \sum_{\begin{subarray}{c} 1\leq k <\ell \leq n \\ \varepsilon_1,\varepsilon_2 = \pm 1\end{subarray}} a_0\widehat{f}(c_k^{\varepsilon_1}c_{\ell}^{\varepsilon_2})
\sum_{\begin{subarray}{c}\underline{j} \in \N^n \\ |\underline{j}|=m-1 \ \null  \end{subarray}} 
B_{n,m}(\underline{j},k,\ell) \prod_{i=1}^n a_1(i)^{2j_i} a_1(k)a_1(\ell) \Big],\\
\end{eqnarray*}
where 
\begin{eqnarray*}
B_{n,m}(\underline{j},k,\ell) &=& \Big|\Big\{(g_1,\ldots,g_{2m+1}) \in \F_n^{2m+1} \; : \; g_1 \cdots g_{2m+1}=e, \, g_1=c_kc_{\ell}, \\
&& \quad (g_i)_{2\leq i \leq 2m+1} \sim \omega_n(\underline{j},\underline{j}') \mbox{ with } j'_{i}=j_i+\delta_{i,k}+\delta_{i,\ell} \ \forall \; i\Big\}\Big|.
\end{eqnarray*}
By completing squares as follows 
$$a_0\widehat{f}(c_k^{\varepsilon_1}c_{\ell}^{\varepsilon_2}) a_1(k)a_1(\ell) \leq \frac12 \Big(a_0^2\widehat{f}(c_k^{\varepsilon_1}c_{\ell}^{\varepsilon_2})^2+ a_1(k)^2a_1(\ell)^2\Big) $$
for $1\leq k \leq \ell \leq n$, $\varepsilon_1,\varepsilon_2 = \pm 1$, 
and rearranging we find 
\begin{eqnarray*}
\lefteqn{\hskip-15pt \widehat{f}(e) s_{2m+1}[(2,1,\cdots,1)](r)} \\ [5pt]
&\leq & \frac{2m+1}{2^m} B_{n,m} r^{2m+2}\alpha_0\alpha_1^{m-1}\alpha_2
+\frac{2m+1}{2^{m+1}}n(2n-1)  \widetilde{B}_{n,m}r^{2m+2}\alpha_1^{m+1},
\end{eqnarray*}
where 
\begin{eqnarray*}
B_{n,m}&=&\max \left\{ \frac{B_{n,m}(\underline{j},k,\ell)}{\binom{m-1}{j_1,\ldots,j_n}} \; : \; \underline{j} \in \N^n, \ |\underline{j}|=m-1, \ 1\leq k \leq \ell \leq n\right\}, \\
\widetilde{B}_{n,m}&=&\max 
\left\{\begin{array}{l} 
\displaystyle\frac{B_{n,m}(\underline{j}',k,\ell)}{\binom{m+1}{j_1,\ldots,j_n}} \; : \; \underline{j} \in \N^n, \  |\underline{j}|=m+1, \ 1\leq k \leq \ell \leq n, \\[7pt] 
 j'_i=j_i-\delta_{i,k}-\delta_{i,\ell} \ \forall \; i, \, j_k, j_\ell \geq 1 \mbox{ if } k<\ell, \  j_k \geq 2 \mbox{ if } k=\ell  
\end{array}\right\}.
\end{eqnarray*}
We claim that 
\renewcommand{\theequation}{C1.4}
\addtocounter{equation}{-1}
\begin{equation}\label{Bnm}
 B_{n,m}, \widetilde{B}_{n,m} \leq 3\cdot 2^{m-1} \frac{(2m)!}{(m-1)! (m+2)!},
\end{equation} 
which ends the proof of $\varepsilon'2)$. 
It remains to prove the claims \eqref{Anm} and \eqref{Bnm}. 
Observe that for $\underline{j} \in \N^n$ with $|\underline{j}|=m$ we have 
$$  A_{n,m}(\underline{j}) \leq \mathcal{C}(m) \binom{m}{j_1,\ldots,j_n} 2^m, $$
where $\mathcal{C}(m)=\binom{2m}{m}\frac{1}{m+1}$ denotes the Catalan number, which counts the number of non-crossing pair partitions of the set $\{1,\cdots,2m\}$. 
Indeed, to order the $2m$ letters $\omega_n(\underline{j},\underline{j})$ so that the product gives $e$, 
we may choose a non-crossing partition of $\{1,2,\ldots,2m\}$, then associate to the $m$ pairs obtained one of the $j_i$ pairs $(c_i,c_i^{-1})$ for each $i$, and choose for each the order of the two letters. 
This proves the inequality  in the claim \eqref{Anm}. 
The equality comes from the choice $\underline{j}\sim (1,\ldots,1,\underline{0}) \in \N^n$ if $n\geq m$, since in this case the process described above counts exactly once each possibility.     
We use a similar argument to show \eqref{Bnm}. To compute $B_{n,m}(\underline{j},k,\ell)$ for $\underline{j} \in \N^n$ with $|\underline{j}|=m-1$, we need to order the $2m+2$ letters $c_k,c_\ell, \omega_n(\underline{j},\underline{j}')$ with $j'_{i}=j_i+\delta_{i,k}+\delta_{i,\ell}$ for all $1\leq i \leq n$ so that the product gives $e$, knowing that the two first positions are occupied by $c_k$ and $c_\ell$. Hence, we need to fix a non-crossing pair partition of  $\{1,2,\ldots,2m+2\}$ such that $(1,2)$ is not a pair. 
We denote by $\widetilde{\mathcal{C}}(m+1)$ the number of such partitions. 
Since the two pairs containing $1$ and $2$ are necessarily associated to the letters  $(c_k,c_k^{-1})$ and $(c_\ell,c_\ell^{-1})$ respectively, 
it remains to associate the other $m-1$ pairs to the letters $\omega_n(\underline{j},\underline{j}')\setminus \{c_k^{-1},c_\ell^{-1}\}=\omega_n(\underline{j},\underline{j})$. 
We obtain 
$$  B_{n,m}(\underline{j},k,\ell) \leq \widetilde{\mathcal{C}}(m+1) \binom{m-1}{j_1,\ldots,j_n} 2^{m-1}.$$
Since $\mathcal{C}(m)$ also counts the number of non-crossing pair partitions of $\{1,\cdots,2m+2\}$ such that $(1,2)$ is a pair, we compute 
$$\widetilde{\mathcal{C}}(m+1)=\mathcal{C}(m+1) -\mathcal{C}(m)=\frac{3\cdot (2m)!}{(m-1)!(m+2)!},$$
which yields the first part of \eqref{Bnm}. 
If $n \geq m-1$, the choice $\underline{j}\sim (1,\cdots,1, \underline{0})\in \N^n$ gives an equality. 
For the second part, fix $1\leq k \leq \ell \leq n$. Applying the preceding result to $B_{n,m}(\underline{j}',k,\ell)$ for $|\underline{j}|=m+1$ and $j'_i=j_i-\delta_{i,k}-\delta_{i,\ell}$ for all $1\leq i \leq n$, it suffices to observe that 
$$ \binom{m-1}{j'_1,\ldots,j'_n}\binom{m+1}{j_1,\ldots,j_n}^{-1}=
 \left\{\begin{array}{ll}
\displaystyle\frac{j_k(j_k-1)}{m(m+1)} & \mbox{ if } k=\ell \\
\displaystyle\frac{j_k j_\ell}{m(m+1)} & \mbox{ if } k<\ell 
\end{array}\right.
\leq 1,$$
recalling that if $k=\ell$ then $j_k \geq 2$, and if $k< \ell$ then $j_k, j_\ell \geq 1$. \fin

\vskip5pt

\noindent {\bf C2. Estimates for triangular groups.} We now turn to the triangular group estimates. Since the computations are similar to the ones detailed above in the free group case, we will only prove the $(\varepsilon)$-estimates and one $(\alpha)$-estimate.

\noindent {\bf Estimates $(\alpha)$ for $\Delta_{\alpha\beta\gamma}$.} The only $(\alpha)$-estimate which differs in nature from those for the free group is $\alpha$iv) for $k \ge L/3 \ge 5$. In the terminology of Proposition \ref{s23-triangle}, note that $k = \ell_2 + \ell_3 - 2m \ge L/3$ iff $m \le K(\ell_2,\ell_3,L)$. In particular, we have to estimate $$\sum_{\ell_2 \ge \ell_3 \ge 1} \sum_{\begin{subarray}{c} m=0 \\ \ell_2 + \ell_3 -2m = k \end{subarray}}^{\lfloor \ell_3/2 \rfloor} \frac32 2^{\ell_3-1} M(k,\ell_2,\ell_3) \, r^{\ell_2+\ell_3+k}.$$
When $k=5$, our estimate follows from $$\Omega(5) = \big\{ (3,2,0), (4,1,0), (5,2,1), (4,3,1), (5,4,2) \big\},$$ where we recall that $$\Omega(k) = \Big\{ (\ell_2, \ell_3,m) \, : \, \ell_2 + \ell_3 - 2m = k, \ 0 \le m \le \frac{\ell_3}{2} \ \mbox{and} \ \ell_2 \ge \ell_3 \ge 1 \Big\}.$$ For $k \ge 6$, using $M(k,\ell_2,\ell_3) \le 6$, the goal is to show that
$$\frac92 r^{2k} \sum_{\begin{subarray}{c} 1 \le i \le j \le k \\ i+j-k \ \mathrm{even} \end{subarray}} 2^i r^{i+j-k} \ \le \ \frac92 r^{2k} \frac{2^{k/2}(3-4r^2)}{(1-r^2)(1-2r^2)}.$$
The sum in the left hand side can be written as follows
\begin{eqnarray*}
\hskip5pt \sum_{\begin{subarray}{c} 1 \le i \le j \le k \\ i+j-k \ \mathrm{even} \end{subarray}} 2^i r^{i+j-k} & = & \sum_{i=1}^k 2^i \sum_{\begin{subarray}{c} j=i \vee(k-i) \\ i+j-k \ \mathrm{even} \end{subarray}}^k r^{i+j-k} \\ & \le & \sum_{i=1}^{k/2} 2^i \sum_{s=0}^{\infty} r^{2s} + \sum_{i=k/2}^{\infty} 2^i r^{2i-k} \sum_{s=0}^{\infty} r^{2s} \\ [7pt] & \le & \frac{1}{1-r^2} \Big( 2^{\frac{k}{2}+1} + r^{-k} \frac{(2r^2)^{k/2}}{1-2r^2} \Big) \ = \ \frac{2^{k/2}(3-4r^2)}{(1-r^2)(1-2r^2)}. \hskip10pt  \square
\end{eqnarray*}

\noindent {\bf Estimates $(\varepsilon)$ for $\Delta_{\alpha\beta\gamma}$.} The estimate $\varepsilon \mbox{iii})$ follows from \eqref{s5Est}. For the sums $\varepsilon \mbox{i})$, $\varepsilon \mbox{ii})$ and $\varepsilon \mbox{iv})$ we need to be more careful. 
Observe that for $\l \in B_u$ we have $|\l|\leq 6$. Hence if the smallest loop $L$ is greater or equal than $7$, then we are in the free situation and the super-pathological sums $s_u[\l](r)$ associated to the triangular group $\Delta_{\alpha\beta\gamma}$ coincide with the corresponding one for the free group $\Z_2\ast \Z_2\ast \Z_2$ when $\l \in B_u$. As we did for $\F_2$, we can easily modelize the group $\Z_2\ast \Z_2\ast \Z_2$ with a computer and estimate precisely the required sums $s_u[\l](r)$ for $\l \in B_u$. We detail below the outcome of the computer and the way we used to complete squares. In the sequel we have 
$$ \alpha_1= \sum_{i=1}^3 a_1(i)^2, \quad \alpha_2= 2\sum_{i=1}^3 a_2(i)^2 \quad \mbox{and} \quad \alpha_3=2\sum_{j=1}^3 a_3(i_j)^2+\sum_{1\leq i\neq j \leq 3} a_3(k_{ij})^2,$$ where $\{k_{i \ell} \,: \, 1 \leq i \leq 3, \, i\neq \ell\} \cap \{k_{\ell j} \, : \,1\leq j \leq 3, \, j\neq \ell \} =\emptyset$ for $1 \le \ell \le 3$ fixed. We get 

\begin{itemize}
\item[$\varepsilon$\mbox{i})] $\displaystyle s_4[(1,1,1,1)](r) \, = \, \Big[ \sum_{i=1}^3 a_1(i)^4 + 4\sum_{1\leq i<j\leq 3}a_1(i)^2a_1(j)^2 \Big] r^4 \, \le \, 2r^4\alpha_1^2$. 
\end{itemize}

\begin{eqnarray*}
\lefteqn{\hskip2pt \varepsilon \mbox{ii}) \hskip5pt s_4[(3,1,1,1)](r)} \\ & = & \Big[ 8a_1(1)a_1(2)a_1(3) \sum_{j=1}^3 a_3(i_j) +4 \sum_{1\leq i\neq j \leq 3}a_1(i)^2a_1(j)a_3(k_{i,j})\Big] r^6 \\ & \leq &  \Big[ 4 \sum_{1\leq i<j\leq 3} a_1(i)^2a_1(j)^2+\frac{4}{3}\Big(\sum_{i=1}^3 a_1(i)^2\Big)\Big(\sum_{j=1}^3 a_3(i_j)^2\Big) \Big] r^6 \\ [5pt] & + & \Big[ \sum_{1\leq i\neq j \leq 3} \big(a_1(i)^4 + a_1(j)^2a_3(k_{i,j})^2 + a_1(i)^2a_1(j)^2+a_1(i)^2a_3(k_{i,j})^2\big) \Big] r^6 \\ [5pt] & \leq &  3r^6\alpha_1^2+r^6\alpha_1\alpha_3.
\end{eqnarray*}

\begin{eqnarray*}
\lefteqn{\hskip9pt \varepsilon \mbox{iv}) \hskip5pt s_6[(1,1,1,1,1,1)](r)} \\
& = & \Big[ \sum_{i=1}^3a_1(i)^6 + 9\sum_{1\leq i \neq j\leq 3} a_1(i)^4a_1(j)^2+30a_1(1)^2a_1(2)^2a_1(3)^2 \Big] r^6 
\ \le \ 5r^6\alpha_1^3. 
\end{eqnarray*}
\fin

\noindent {\bf C3. Estimates for finite cyclic groups.} Only $(\alpha)$ and $(\varepsilon)$ estimates are new.

\noindent {\bf Estimates $(\alpha)$ for $\Z_n$.} The coefficient in $\alpha$i) follows once again from the identity we gave for $s_2(r)$ in Section \ref{Method}. To obtain the left-coefficients for $s_3(r)$ we use Proposition \ref{s23cyclic}. When $\k = (k_1,k_2,\underline{0})$, we have to estimate 
$$\frac14 \Big( \sum_{\begin{subarray}{c} \lfloor \frac{n}{2} \rfloor \ge \ell \ge k_1 \ge k_2 \ge 1 \\ \ell + k_1 + k_2 = n \end{subarray}} M(\ell,k_1,k_2) r^n + \sum_{\begin{subarray}{c} \lfloor \frac{n}{2} \rfloor \ge \ell \ge k_1 \ge k_2 \ge 1 \\ \ell = k_1 + k_2 \end{subarray}} M(\ell,k_1,k_2) r^{2 k} \Big),$$ which in turn can be written as the sum of two single terms as follows $$\frac14 \Big( \delta_{k_1 \le n-k_1-k_2 \le \lfloor \frac{n}{2} \rfloor} M(n-k_1-k_2,k_1,k_2) r^n + \delta_{k_1+k_2 \le \lfloor \frac{n}{2} \rfloor} M(k_1+k_2,k_1,k_2) r^{2(k_1+k_2)} \Big).$$ Our estimates $\alpha$ii) and $\alpha$iii) then follow by direct substitution and elementary inequalities taking into account that we are assuming $q \le n \, (\Leftrightarrow r \ge 1/\sqrt{n-1})$ when $n$ is odd and $n \ge 6$. For $\alpha$iv) it suffices to observe that $M(\cdot,k_1,k_2) \le 3 M(k_1,k_2)$ and $\delta_{k_1 \le n-k_1-k_2} r^n + \delta_{k_1+k_2 \le \lfloor \frac{n}{2} \rfloor} r^{2(k_1+k_2)} \, \le \, 2 r^{k_1+k_2+2}$. It remains to prove $\alpha$v). By Proposition \ref{s23cyclic}, we need to estimate $$\frac12 \sum_{\begin{subarray}{c} \lfloor \frac{n}{2} \rfloor \ge k \ge \ell_2 \ge \ell_3 \ge 1 \\ k = \ell_2 + \ell_3 \end{subarray}} M(k,\ell_2,\ell_3) r^{2k} + \frac12 \sum_{\begin{subarray}{c} \lfloor \frac{n}{2} \rfloor \ge k \ge \ell_2 \ge \ell_3 \ge 1 \\ k + \ell_2 + \ell_3 = n \end{subarray}} M(k,\ell_2,\ell_3) r^n \, = \, \mathrm{A} + \mathrm{B}.$$ Using $M(k,\ell_2,\ell_3) \le 3 \delta_{k=2} + 6 \delta_{k > 2}$ we get 
\begin{eqnarray*}
\mathrm{A} & \le & \mbox{$\frac32$} r^4 \delta_{k=2} + \mbox{$\frac32$} k r^{2k} \delta_{k > 2}, \\ \mathrm{B} & \le & 3 \Big| \Big\{ (i,j) \; : \; 1 \le i \le j \le k \le \lfloor \mbox{$\frac{n}{2}$} \rfloor, \ i+j+k=n \Big\} \Big| \, r^n.
\end{eqnarray*}
Any $(i,j)$ belonging to the set above satisfies $j=n-k-i$ and $2i+k \le n \le i+2k$. The later condition implies $n-2k \le i \le \frac12 (n-k)$, which means that the cardinal of that set is bounded by $\frac12 (n-k-2n+4k) +1 \le \frac12k+1$ since $k \le \frac{n}{2}$. On the other hand, $r^n \le r^{2k}$ for the same reason and $\mathrm{B} \le (\frac32 k + 3) r^{2k} \delta_{k \ge 3}$. \fin


\noindent {\bf Estimates $(\varepsilon)$ for $\Z_n$.} If $|g|=1$, we have by commutativity 
\begin{eqnarray*}
\hskip10pt s_4[(1,1,1,1)](r) & = & \binom{4}{2} \widehat{f}(g)^4 r^4 \ = \ \frac32 \alpha_1^2 r^4, \\ [3pt] s_4[(3,1,1,1)](r) & = & \binom{4}{1} \Big( \delta_{n=6} \widehat{f}(g^3) 2 \widehat{f}(g)^3 + \delta_{n>6} 2 \widehat{f}(g^3) \widehat{f}(g)^3 \Big) r^6 \\ [5pt] & \le & 4 \big( a_3^2 a_1^2 + a_1^4 \big) r^6 \ = \  \big( \delta_{n=6} 2 \alpha_1 \alpha_3 + \delta_{n>6} \alpha_1 \alpha_3 + \alpha_1^2 \big) r^6. \hskip10pt \ \square
\end{eqnarray*}

\subsection*{D. Technical inequalities}

In this appendix we list all the technical inequalities which appear in the proof of Theorems A1-A3. Given $q \in 2 \Z$ with $q=2s$, we need to verify these inequalities for $0 \le r \le 1/\sqrt{2s-1}$. Nevertheless, it suffices to assume for simplicity that $r = 1/\sqrt{q-1} = \exp(-t_{2,q})$ in what follows, since for smaller values of $r = \exp(-t) < \exp(-t_{2,q})$ we may use the semigroup property for $t = t_{2,q} + \delta$ so that $$\|\T_{\psi,r}\|_{2 \to q} = \|\S_{\psi,t}\|_{2 \to q} \le \|\S_{\psi,\delta}\|_{2 \to 2} \|\S_{\psi,t_{2,q}}\|_{2 \to q} \le 1.$$

\renewcommand{\theequation}{D0.1}
\addtocounter{equation}{-1}

\noindent \textbf{D0. Positivity test for polynomials} We begin describing an algorithmic way to find the best positive integer $s_0$ associated to a given polynomial with positive dominant coefficient $P$ such that $P(s)\geq 0$ for all integers $s\geq s_0$. 
This will be used to prove most of the technical inequalities below. We will use the two following well-known facts. Let $P(s)=c_0+c_1s+\cdots +c_d s^d \in \R[X]$ such that $c_d>0$
\begin{enumerate}
\item[i)] Cauchy's bound. If $t$ is a root of $P$ then 
$$|t|\leq u=\max_{1\leq i\leq d-1}\Big( \frac{|c_i|}{|c_d|}\Big)+1.$$

\vskip3pt

\item[ii)] If $P^{(i)}(v)\geq 0$ for all $0\leq i\leq d-2$, then $P(s)\geq 0$ for all $s\geq v$.
\end{enumerate} 
In order to find the best positive integer $s_0$ such that 
\begin{equation}\label{s0}
P(s)\geq 0 \quad \mbox{for all integers} \quad s\geq s_0,
\end{equation} 
we proceed as follows. We first compute the Cauchy's bound $u$, which satisfies \eqref{s0}. By making the computations for all integers $0 \le v \le u$, we may find an integer $v$ such that $P^{(i)}(v)\geq 0$ for all $0\leq i\leq d-1$. Let $w$ be the smallest such $v$ if it exists, otherwise we set $w=u$. Then $w$ still satisfies \eqref{s0}. Hence the best $s_0$ possible is $\leq w$, and to find that integer it remains to decide whether $P(s)\geq 0$ for all integers $s=w-1\cdots 0$. Finally we set $$s_0 \, = \, \min \Big\{ k \; : \; 0 \leq k\leq w,\; P(s)\geq 0 \ \forall \ k \le s \leq w \Big\}.$$

\noindent \textbf{D1. Technical inequalities for free groups.} Define $$q(n) = 4 \delta_{n=2} + \big( (22n)^{44n} + 2 \big) \delta_{n \ge 3}.$$ Let $n \ge 2$, $s \geq \frac{q(n)}{2}$ and $r = \frac{1}{\sqrt{2s-1}}$, then the following inequalities hold:
\renewcommand{\theequation}{D1.1}
\addtocounter{equation}{-1}

\vskip-5pt

\begin{equation} \label{Reg1-1}
\binom{2s}{2}r^{4}+\binom{2s}{3} \big( \mbox{$\frac32$} r^4 + (n-1)r^6 \big)  \leq s.
\end{equation}

\vskip-5pt

\renewcommand{\theequation}{D1.2}
\addtocounter{equation}{-1}

\begin{equation}\label{Reg1-2}
\binom{2s}{2}r^{6}+\binom{2s}{3} \big( 3 r^6 + 3 r^8 \big)  \leq s.
\end{equation}

\vskip-5pt

\renewcommand{\theequation}{D1.3}
\addtocounter{equation}{-1}

\begin{equation}\label{Reg1-3}
\binom{2s}{2}r^{8}+\binom{2s}{3} \big( \mbox{$\frac92$} r^8 + 6 r^{10} + 3 r^{12} \big)  \leq s.
\end{equation}

\vskip-5pt

\renewcommand{\theequation}{D1.4}
\addtocounter{equation}{-1}

\begin{eqnarray}\label{Reg1-4}
\lefteqn{\hskip18pt \mbox{We have}} \\ 
\nonumber && \binom{2s}{2}r^{2k}+\binom{2s}{3} \Big( \frac{3 k (1-r^2)r^{2k}}{2(1-(2n-1)r^2)} + \frac{6(n-1)r^{2k}}{(2n-1)(1-(2n-1)r^2)^2}\Big)  \leq s. \\ \nonumber \lefteqn{\hskip18pt \mbox{when $n=2$ with $k \geq 5$ and $s\geq 3$ or when $n \ge 3$ with $k \ge 3$ and $s \ge \frac{q(n)}{2}$.}}
\end{eqnarray}

\vskip-15pt

\renewcommand{\theequation}{D1.5}
\addtocounter{equation}{-1}

\begin{eqnarray}\label{Reg2}
\lefteqn{\hskip27pt \mbox{We have}} \\ \nonumber && \frac{r^k}{s(s-1)} \Big[ \binom{2s}{3} \frac{3}{1-r^2}
+\binom{2s}{4} \frac{24nr^{2}}{1-(2n-1)r} + \binom{2s}{5} \frac{80 n^2 r^{3}}{(1-(2n-1)r)^2} \Big] \leq 1. \\ \nonumber \lefteqn{\hskip27pt \mbox{when $n=2$ with $k \geq 16$ and $s\geq 6$ or when $n \ge 3$ with $k \ge 4$ and $s \ge \frac{q(n)}{2}$.}}
\end{eqnarray}

\renewcommand{\theequation}{D1.6}
\addtocounter{equation}{-1}

\vskip-15pt

\begin{eqnarray} \label{Regm}
\lefteqn{\hskip15pt \mbox{If $3\leq m \leq s$, $k\geq m+\mu_q(\F_n,|\cdot|,m)$ and $s\geq 6 \delta_{n=2} + \frac{q(n)}{2} \delta_{n \ge 3}$}} \\ \nonumber && \binom{2s}{2m-1} A_m(s,r,k) + \binom{2s}{2m} B_m(s,r,k) + \binom{2s}{2m+1} C_m(s,r,k) \leq 1, \\ [5pt] \nonumber \lefteqn{\hskip15pt \mbox{where}} \\ [5pt] \nonumber && A_m(s,r,k) = \frac{(2m-1)!(s-m)!}{(m-1)! s!}\Big(\frac{2nr}{1-(2n-1)r}\Big)^{m-2} r^{k+1}, \\ \nonumber \hskip10pt && B_m(s,r,k) = \frac{(2m)!(s-m)!}{m! s!}\Big(\frac{2nr}{1-(2n-1)r}\Big)^{m-1} r^{k+1}, \\ \nonumber \hskip10pt && C_m(s,r,k) = \frac{(2m+1)!(s-m)!}{(m+1)! s!}\Big(\frac{2nr}{1-(2n-1)r}\Big)^{m}r^{k+1}.
\end{eqnarray}

\renewcommand{\theequation}{D1.7}
\addtocounter{equation}{-1}

\begin{eqnarray} \label{Pat} 
\lefteqn{\hskip23pt \mbox{If $22\leq m \leq s$}} \\ \nonumber && \frac{4^m r^{2m}}{1-r^2} \Big[ \binom{2s}{2m-1} \frac{2m-1}{4^3} + \binom{2s}{2m} \frac{m}{4} + \binom{2s}{2m+1} \frac{2m+1}{4} r^2 \Big] \ \leq \ \binom{s}{m}.
\end{eqnarray} 

\vskip-10pt

\renewcommand{\theequation}{D1.8}
\addtocounter{equation}{-1}

\begin{eqnarray} \label{Pat-1-Fn} 
\lefteqn{\hskip12pt \mbox{If $n \ge 3$ and $3\leq m \leq s$}} \\ \nonumber && \binom{2s}{2m-1} A_m(r,s) + \binom{2s}{2m} B_m(r,s) + \binom{2s}{2m+1} C_m(r,s) \leq \binom{s}{m}, \\ [5pt] \nonumber \lefteqn{\hskip12pt \mbox{where}} \\ [5pt] \nonumber 
&& A_m(r,s) = \Big(\frac{3n(2n-1)(2m-1)!r^{2m}}{4(m-2)!(m+1)!}+ (2m-1)\frac{(2n)^{m-2}}{4}\frac{r^{2m+2}}{1-r^2}\Big), \\ [5pt] \nonumber && B_m(r,s) = \Big(\frac{(2m)!}{m!(m+1)!}r^{2m}+ m (2n)^{m-1}\frac{r^{2m+2}}{1-r^2} \Big), 
\\ [5pt] \nonumber && C_m(r,s) = (2m+1)\frac{(2n)^{m}}{4}\frac{r^{2m+4}}{1-r^2}.
\end{eqnarray} 

\vskip-15pt

\renewcommand{\theequation}{D1.9}
\addtocounter{equation}{-1}

\begin{eqnarray} \label{Pat-2-Fn} 
\lefteqn{\hskip-8pt \mbox{If $n \ge 3$}} \\ \nonumber && \binom{2s}{3}\frac32 r^4 
\, + \,\binom{2s}{4} \big(2r^4+ 4n r^{6} \big) 
\, + \, \binom{2s}{5}5n^2 r^{8}
\leq  \ \frac{s(s-1)}{2}.
\end{eqnarray} 

\vskip-10pt

\renewcommand{\theequation}{D1.10}
\addtocounter{equation}{-1}

\begin{eqnarray} \label{Pat-3-Fn} 
\lefteqn{\hskip15pt \mbox{If $n \ge 3$ and $3\leq m \leq s$}} \\ \nonumber && \binom{2s}{2m-1} A_m(r,s) + \binom{2s}{2m} B_m(r,s) + \binom{2s}{2m+1} C_m(r,s) \leq \binom{s}{m}, \\ [5pt] \nonumber \lefteqn{\hskip15pt \mbox{where}} \\   \nonumber
&& A_m(r,s) = \frac{(2n)^{m-2}}{4m}(2m-1)(2m-2) \frac{r^{2m+2}}{1-r^2} 
\Big[1+ \frac{(2n-1)}{4}(2m-3)\Big] \\[5pt] \nonumber && B_m(r,s) = 
(2n)^{m-1}(2m-1)\frac{r^{2m+2}}{1-r^2} 
\Big[ 1 +\frac{2n-1}{4}(2m-2)r^2\Big] \\[5pt] \nonumber
&& C_m(r,s) =  
 \frac{\frac32(2m+1)!}{m!(m+2)!}r^{2m+2}  \\ [5pt] \nonumber && \hskip38pt + \, 
 \frac{(2n)^{m}}{2}(2m+1)\frac{r^{2m+4}}{1-r^2} 
\Big[1 +  \frac{2n-1}{4}(2m-1)\Big].
\end{eqnarray} 

\vskip-10pt

\renewcommand{\theequation}{D1.11}
\addtocounter{equation}{-1}

\begin{eqnarray} \label{Pat-4-Fn} 
\lefteqn{\mbox{If $n \ge 3$}} \\ \nonumber && \hskip-12pt \binom{2s}{3}3 r^6 
\, + \,\binom{2s}{4}\frac{r^{6}}{1-r^2} (12n+6n(2n-1)r^2) 
\\ \nonumber && \hskip30pt + \hskip4pt \binom{2s}{5}r^6\Big( \frac{15}{2}+\frac{r^2}{1-r^2}5n^2(6n+1)\Big)
\leq  \ s(s-1).
\end{eqnarray} 

\noindent \textbf{D2. Technical inequalities for triangular groups.} 

\noindent Let $s\geq 2$ and $r = \frac{1}{\sqrt{2s-1}}$, then the following inequalities hold:

\vskip-5pt

\renewcommand{\theequation}{D2.1}
\addtocounter{equation}{-1}

\begin{equation}\label{Reg1-1-triangle}
\binom{2s}{2}r^{4}+\binom{2s}{3} \big( \mbox{$\frac32$} r^4 + \mbox{$\frac12$} r^6 \big)  \leq s.
\end{equation}

\vskip-5pt

\renewcommand{\theequation}{D2.2}
\addtocounter{equation}{-1}

\begin{equation}\label{Reg1-2-triangle}
\binom{2s}{2}r^{6}+\binom{2s}{3} \big( 3 r^6 + \mbox{$\frac32$} r^8 \big)  \leq s.
\end{equation}

\vskip-5pt

\renewcommand{\theequation}{D2.3}
\addtocounter{equation}{-1}

\begin{equation} \label{Reg1-3-triangle}
\binom{2s}{2}r^{8}+\binom{2s}{3}\frac{3r^{8}}{2(1-2r^2)} \Big( 4(1-r^2) + \frac{1}{1-2r^2} \Big) \leq s.  
\end{equation}

\vskip-5pt

\renewcommand{\theequation}{D2.4}
\addtocounter{equation}{-1}

\begin{equation} \label{Reg1-4-triangle}
\binom{2s}{2}r^{10}+\binom{2s}{3}\big( 27r^{10} + 45r^{12} + 36r^{14} \big) \leq s.  
\end{equation}

\vskip-10pt

\renewcommand{\theequation}{D2.5}
\addtocounter{equation}{-1}

\begin{eqnarray} \label{Reg1-5-triangle}
\lefteqn{\hskip-48pt \mbox{If $ k \ge 6$}} \\ 
\nonumber && \binom{2s}{2}r^{2k}+\binom{2s}{3} \frac{9r^{2k} 2^{k/2}(3-4r^2)}{2(1-r^2)(1-2r^2)}  \leq s. 
\end{eqnarray}

\vskip-10pt

\renewcommand{\theequation}{D2.6}
\addtocounter{equation}{-1}

\begin{eqnarray}\label{Reg2-triangle}
\lefteqn{\hskip32pt \mbox{If $ k \geq15$ and $s\geq 3$}} \\ \nonumber && \binom{2s}{3}\frac{3r^{k}}{s(s-1)(1-r^2)}
+\binom{2s}{4}\frac{36}{s(s-1)}\frac{r^{k+2}}{1-2r}
+\binom{2s}{5}\frac{180}{s(s-1)}\frac{r^{k+3}}{(1-2r)^2}\leq 1. 
\end{eqnarray}

\vskip-10pt

\renewcommand{\theequation}{D2.7}
\addtocounter{equation}{-1}

\begin{eqnarray} \label{Regm-triangle}
\lefteqn{\hskip33pt \mbox{If $3\leq m \leq s$ and $k\geq m+\mu_q(\Delta_{\alpha\beta\gamma},|\cdot|,m)$}} \\ \nonumber && \Big(\frac{3r}{1-2r}\Big)^{m-2}r^{k+1} \left[ \binom{2s}{2m-1} A_m(s) + \binom{2s}{2m} B_m(s) + \binom{2s}{2m+1} C_m(s) \right] \leq 1, 
\end{eqnarray} 
\noindent where $$\hskip20pt \underbrace{\frac{(2m-1)!(s-m)!}{(m-1)! s!}}_{A_m(s)}, \quad \underbrace{\frac{(2m)!(s-m)!}{m! s!}\frac{3r}{1-2r}}_{B_m(s)}, \quad \underbrace{\frac{(2m+1)!(s-m)!}{(m+1)! s!}\Big(\frac{3r}{1-2r}\Big)^2}_{C_m(s)}.$$

\vskip-10pt

\renewcommand{\theequation}{D2.8}
\addtocounter{equation}{-1}

\begin{eqnarray} \label{Pat-triangle}
\lefteqn{\hskip20pt \mbox{If $15\leq m \leq s$}} \\ \nonumber && \frac{(3r^2)^{m}}{1-r^2} \Big[ \binom{2s}{2m-1} \frac{2m-1}{36} + \binom{2s}{2m} \frac{m}{3} + \binom{2s}{2m+1} \frac{2m+1}{4} r^2 \Big] \leq \binom{s}{m}.
\end{eqnarray}

\noindent \textbf{D3. Technical inequalities for finite cyclic groups.}

\noindent Let $s\geq 2$ and $r = \frac{1}{\sqrt{2s-1}}$, then the following inequalities hold:

\vskip-5pt

\renewcommand{\theequation}{D3.1}
\addtocounter{equation}{-1}

\begin{equation}\label{Reg1-1-cyclic}
\binom{2s}{2}r^{4}+\binom{2s}{3} \mbox{$\frac32$} r^4  \leq s.
\end{equation}

\vskip-10pt

\renewcommand{\theequation}{D3.2}
\addtocounter{equation}{-1}

\begin{eqnarray}\label{Reg1-2-cyclic}
\lefteqn{\hskip-65pt \mbox{If $ k \ge 3$}} 
\\ \nonumber && \binom{2s}{2}r^{2k} + \binom{2s}{3} 3 (k+1) r^{2k}  \leq s.
\end{eqnarray}

\vskip-15pt

\renewcommand{\theequation}{D3.3}
\addtocounter{equation}{-1}

\begin{eqnarray}\label{Reg2-cyclic}
\lefteqn{\hskip-16pt \mbox{If $ k \geq 5$}} \\ \nonumber && \frac{r^k}{s(s-1)} \Big[ \binom{2s}{3} 3r^{2}
+ \binom{2s}{4} \frac{24 r^{2}}{1-r} + \binom{2s}{5} \frac{80 r^{3}}{(1-r)^2} \Big] \leq 1. 
\end{eqnarray}

\vskip-10pt

\renewcommand{\theequation}{D3.4}
\addtocounter{equation}{-1}

\begin{eqnarray} \label{Regm-cyclic}
\lefteqn{\hskip33pt \mbox{If $3\leq m \leq s$ and $k\geq m+\mu_q(\Z_n,|\cdot|,m)$}} \\ \nonumber && \Big(\frac{2r}{1-r}\Big)^{m-2}r^{k+1} \left[ \binom{2s}{2m-1} A_m(s) + \binom{2s}{2m} B_m(s) + \binom{2s}{2m+1} C_m(s) \right] \leq 1,
\end{eqnarray} 
\noindent where $$\hskip20pt \underbrace{\frac{(2m-1)!(s-m)!}{(m-1)! s!}}_{A_m(s)}, \quad \underbrace{\frac{(2m)!(s-m)!}{m! s!}\frac{2r}{1-r}}_{B_m(s)}, \quad \underbrace{\frac{(2m+1)!(s-m)!}{(m+1)! s!}\Big(\frac{2r}{1-r}\Big)^2}_{C_m(s)}.$$

\vskip-10pt

\renewcommand{\theequation}{D3.5}
\addtocounter{equation}{-1}

\begin{eqnarray} \label{Pat1-cyclic}
\lefteqn{\hskip20pt \mbox{If $7 \le m \le s$}} \\ \nonumber && \frac{(2r^2)^{m}}{1-r} \Big[ \binom{2s}{2m-1} \frac{2m-1}{16} + \binom{2s}{2m} \frac{m}{2} + \binom{2s}{2m+1} \frac{2m+1}{4} r^2 \Big] \leq \binom{s}{m}.
\end{eqnarray}

\vskip5pt

\noindent \textbf{D4. Proofs.} By clear similarities in the arguments, we shall only prove the technical inequalities for free groups. 
First note that the left hand sides of inequalities \eqref{Reg1-4}, \eqref{Reg2} and \eqref{Regm} are decreasing in $k$, hence it suffices to prove them for the smallest value of $k$. Moreover, inequalities \eqref{Reg1-1}-\eqref{Reg1-3} and \eqref{Reg1-4} for $n=2, k=5$ can be rewritten as polynomial inequalities in the variable $s$, since $r^2 = \frac{1}{2s-1}$ and $r$ is raised to even powers. In particular, these inequalities can be justified with computer assistance by means of the positivity test for polynomials in D0 above. On the other hand, the inequality \eqref{Reg2} for $n=2$ and $k=16$ is equivalent to some polynomial inequality in the variables $s$ and $\sqrt{2s-1}$. We may adapt the positivity test for polynomials presented in D0 to decide whether such polynomial in those variables is positive for all integer $s \ge s_0$. Thus \eqref{Reg2} for $n=2$ and $k=16$ can also be proved by using computer assistance. In the case $n \ge 3$, inequalities \eqref{Reg1-1}, \eqref{Reg1-4}, \eqref{Reg2}, \eqref{Pat-2-Fn} and \eqref{Pat-4-Fn} follow easily from the large order of growth of $q(n)$ we impose. We will only detail here the proof of the crucial inequality (D1.6) for free groups, which yields the bounded critical functions $\mu_q(\F_2,|\cdot|,m)$ and $\mu_q(\F_n,|\cdot|,m)$ given in Paragraph \ref{numericalF2}. In fact the argument we will use to prove (D1.6) in the $\F_2$-case can be adapted to show the similar technical inequalities (D2.6) and (D3.4) in the case of triangular groups and cyclic groups respectively. 
In the proof of (D1.6) for $\F_n$, it will appear that we need a large order of growth for $q(n) \sim n^n$ to get a critical function which is uniformly bounded in $n$. Then, the remaining inequalities \eqref{Pat}, \eqref{Pat-1-Fn} and \eqref{Pat-3-Fn} can be proved by using the same ideas. 

The key point in (D1.6) for $\F_2$, (D2.6) and (D3.4) is that these inequalities hold true with $\mu_q(\G,\psi,m)=1$ when $M_0\leq m \leq s$ for some $M_0\geq 3$. 
For $N$ fixed ($N=4$ for $\F_2$, $N=3$ for triangular groups and $N=2$ for cyclic groups) we set 
\begin{eqnarray*} 
&& A_m^N(s,r,k) = \frac{(2m-1)!(s-m)!}{(m-1)! s!}\Big(\frac{Nr}{1-(N-1)r}\Big)^{m-2} r^{k+1}, \\ 
 \hskip10pt && B_m^N(s,r,k) = \frac{(2m)!(s-m)!}{m! s!}\Big(\frac{Nr}{1-(N-1)r}\Big)^{m-1} r^{k+1}, \\ 
 \hskip10pt && C_m^N(s,r,k) = \frac{(2m+1)!(s-m)!}{(m+1)! s!}\Big(\frac{Nr}{1-(N-1)r}\Big)^{m}r^{k+1},
\end{eqnarray*}
for $3\leq s \leq m$, $0\leq r <\frac{1}{N-1}$ and $k \in \N$. 
Let us prove that 
$$\binom{2s}{2m-1} A_m^N(s,r,m+1) + \binom{2s}{2m} B_m^N(s,r,m+1) + \binom{2s}{2m+1} C_m^N(s,r,m+1) \leq 1$$
for $M_0\leq m \leq s$ and $r=1/\sqrt{2s-1}$, and we will give the numerical proof for $N=4$ (in the $\F_2$-case). 
The (finitely many) remaining inequalities for (D1.6) in the $\F_2$-case, namely for $3 \leq m \leq M_0-1$, can be justified by using our adapted positivity test for polynomials. We write
\begin{eqnarray*}
\binom{2s}{2m-1} A_m^N(s,r,m+1) & = & \frac{4(2N)^{m-2}\prod_{j=1}^{m-2}(2s-2j-1)}{(m-1)! (2s-1)^{m/2}(\sqrt{2s-1}-(N-1))^{m-2}},\\[5pt]
\binom{2s}{2m} B_m^N(s,r,m+1) & = & \frac{2(2N)^{m-1}\prod_{j=1}^{m-1}(2s-2j-1)}{m! (2s-1)^{m/2}(\sqrt{2s-1}-(N-1))^{m-1}},\\[5pt]
\binom{2s}{2m+1} C_m^N(s,r,m+1) & = & \frac{(2N)^{m} (2s-2m) \prod_{j=1}^{m-1}(2s-2j-1)}{(m+1)! (2s-1)^{m/2}(\sqrt{2s-1}-(N-1))^{m}}.
\end{eqnarray*}
Observe that for $m \geq m_0$ we have
\begin{itemize}
\item $\prod_{j=1}^{m-2}(2s-2j-1) \leq (2s-1)^{(m-2)/2}(2s-m_0)^{(m-2)/2}$, 
\vspace{0.2cm}
\item $\prod_{j=1}^{m-1}(2s-2j-1) \leq (2s-1)^{(m-1)/2}(2s-m_0)^{(m-1)/2}$,
\vspace{0.2cm}
\item $\prod_{j=1}^{m-1}(2s-2j-1)(2s-2m) \leq (2s-1)^{m/2}(2s-m_0)^{m/2}$. 
\end{itemize}
This implies the following estimates for $m\geq m_0$ 
\begin{eqnarray*}
\binom{2s}{2m-1} A_m^N(s,r,m+1) & \le & \frac{4(2N)^{m-2}}{(m-1)! (2s-1)} [k_{m_0}(s)]^{\frac{m-2}{2}},\\[5pt]
\binom{2s}{2m} B_m^N(s,r,m+1) & \le & \frac{2(2N)^{m-1}}{m! \sqrt{2s-1}}[k_{m_0}(s)]^{\frac{m-1}{2}},\\[5pt]
\binom{2s}{2m+1} C_m^N(s,r,m+1) & \le & \frac{(2N)^{m}}{(m+1)!}[k_{m_0}(s)]^{\frac{m}{2}},\\
\end{eqnarray*}
where 
$$ k_{m_0}(s)=\frac{2s-m_0}{2s+N^2-2N-2(N-1)\sqrt{2s-1}}.$$
We may find $s(m_0)$ such that $k_{m_0}(s)$ decreases for $s\geq s(m_0)$. 
Hence by setting $\lambda(m_0)=2N\sqrt{k_{m_0}(s(m_0))}$, for $m_0\leq m \leq s$ and $s\geq s(m_0)$ we get
\begin{eqnarray*}
\lefteqn{\hskip-35pt \binom{2s}{2m-1} A_m^N(s,r,m+1) + \binom{2s}{2m} B_m^N(s,r,m+1) + \binom{2s}{2m+1} C_m^N(s,r,m+1)} \\
& \leq & \frac{4\lambda(m_0)^{m-2}}{(m-1)!(2m-1)}+ \frac{2\lambda(m_0)^{m-1}}{m!\sqrt{2m-1}} + \frac{\lambda(m_0)^{m}}{(m+1)!}\ =: \ H_{m_0}(m).
\end{eqnarray*}
Since the sequence $u_m=\frac{\lambda^{m-2}}{(m-1)!}$ decreases for $\lambda \leq m$, we deduce that $H_{m_0}(m)$ decreases for $m \geq \lambda(m_0)$.  
Thus it suffices to find $M_0\geq \max\{m_0,s(m_0), \lambda(m_0)\}$ satisfying $H_{m_0}(M_0)\leq 1$. 
This implies the desired inequality for any $M_0\leq m \leq s$. 
In the $\F_2$-case $N=4$ and we take $m_0=22$. Then $s(m_0)=25$, $\lambda(m_0)=2\sqrt{28}\simeq 10.6$ and $M_0=25$ give the required inequality $H_{m_0}(M_0)\leq 1$. 

We now turn to the proof of (D1.6) in the general case $\F_n$, for any $n\geq 3$. We need to prove 
$$\binom{2s}{2m-1} A_m^{2n}(s,r,m+2) + \binom{2s}{2m} B_m^{2n}(s,r,m+2) + \binom{2s}{2m+1} C_m^{2n}(s,r,m+2) \leq 1$$
for all $n\geq 3$, $3\leq m \leq s$ and $s \geq \frac{q(n)}{2}$. 
In that case we will be more brutal than before and use 
\begin{itemize}
\item $2s-j\leq 2s-1$ for any $j\geq 1$,
\item $\displaystyle\frac{1}{\sqrt{2s-1}-2n+1}\leq \displaystyle\frac{2}{\sqrt{2s-1}}$ for $s\geq 8n^2-8n+\frac{5}{2}$,
\item Stirling's formula: $m! \geq \sqrt{2\pi m}\Big(\displaystyle\frac{m}{e}\Big)^m$.
\end{itemize}
Since $8e \leq 22$, we obtain 
\begin{eqnarray*}
\binom{2s}{2m-1} A_m^{2n}(s,r,m+2) & \le & \frac{4e}{\sqrt{2\pi}(m-1)^{3/2}(2s-1)^{3/2}}\Big(\frac{22n}{m-1}\Big)^{m-2},\\[5pt]
\binom{2s}{2m} B_m^{2n}(s,r,m+2) & \le &  \frac{2e}{\sqrt{2\pi}m^{3/2}(2s-1)}\Big(\frac{22n}{m}\Big)^{m-1},\\[5pt]
\binom{2s}{2m+1} C_m^{2n}(s,r,m+2) & \le &  \frac{e}{\sqrt{2\pi}(m+1)^{3/2}\sqrt{2s-1}}\Big(\frac{22n}{m+1}\Big)^{m}.\\
\end{eqnarray*}
We claim that each term in the right hand side above is less than or equal to $\frac{1}{3}$ for $n\geq 3$, $3\leq m \leq s$ and $s \geq \frac{q(n)}{2}$, which will complete the proof. Let us prove it for the third one, the proof of the two other estimates being similar. 
We can show that 
$$c:= \frac{e}{\sqrt{2\pi}(m+1)^{3/2}\sqrt{2s-1}}\Big(\frac{22n}{m+1}\Big)^{m} \leq \max \Big\{\frac{e}{8\sqrt{2\pi}66^{66}}, \frac{e}{512\sqrt{2\pi} \cdot 66^2}\Big\} \leq \frac{1}{3}$$
for any $n\geq 3$, $3\leq m \leq s$ and $s \geq \frac{q(n)}{2}$. 
Indeed, if $22n \leq m+1$ then for $s\geq \frac{q(n)}{2}$, since $\sqrt{2s-1}\geq (22n)^{22n} \geq 66^{66}$ we get $c \leq \frac{e}{8\sqrt{2\pi}66^{66}}$ for $m \geq 3$. 
If $22n \geq m+1$ then we write
$$ c \leq \frac{e}{\sqrt{2\pi}(m+1)^{3/2+m}(22n)^{22n-m}}\leq \frac{e}{512\sqrt{2\pi} \cdot 66}.$$
This ends the proof of (D1.6) for $\F_n$, $n\geq 3$. \fin

\renewcommand{\theequation}{\arabic{equation}}

\vskip5pt


\bibliographystyle{amsplain}

\hskip10pt

\noindent {\bf Marius Junge} \hfill \textbf{Carlos Palazuelos} \\ \noindent Department of Mathematics \hfill Instituto de Ciencias Matem{\'a}ticas \\ \noindent U. Illinois at Urbana-Champaign \hfill CSIC-UAM-UC3M-UCM \\ \noindent 1409 W. Green St. Urbana, IL. USA \hfill C/ \!\!\! Nicol\'as Cabrera 13-15. \!\!\! Madrid. \!\!\! Spain \\ \noindent \texttt{junge@math.uiuc.edu} \hfill \texttt{carlospalazuelos@icmat.es}

\vskip10pt

\noindent \textbf{Javier Parcet} \hfill \textbf{Mathilde Perrin} \\
\noindent Instituto de Ciencias Matem{\'a}ticas \hfill Instituto de Ciencias Matem{\'a}ticas 
\\ \noindent CSIC-UAM-UC3M-UCM \hfill CSIC-UAM-UC3M-UCM \\ C/ \!\!\! Nicol\'as Cabrera 13-15. \!\!\!
Madrid. \!\!\! Spain \hfill C/ \!\!\! Nicol\'as Cabrera 13-15. \!\!\! Madrid. \!\!\! Spain 
\\ \noindent \texttt{javier.parcet@icmat.es} \hfill \texttt{mathilde.perrin@icmat.es}

\end{document}